\documentclass[3p,times]{elsarticle}

\usepackage{graphicx}
\usepackage[dvipsnames]{xcolor}

\usepackage{mathrsfs}
\usepackage{amsfonts}
\usepackage{amssymb}
\usepackage{amsmath}
\usepackage{dsfont}
\usepackage{rotating}
\usepackage{multirow}
\usepackage{subfigure}

\DeclareMathAlphabet\mathbfcal{OMS}{cmsy}{b}{n}  % bold typing with \mathcal

\newtheorem{theorem}{Theorem}

\newtheorem{remark}[theorem]{Remark}

\newcommand{\diff}[1]{{\mathrm{d}{#1}}}

\newcommand{\RIcolor}[1] {{\leavevmode\color{black} #1}}
\newcommand{\RIIcolor}[1]{{\leavevmode\color{black} #1}}
\newcommand{\Rall}[1]    {{\leavevmode\color{black} #1}}

\begin{document}
	
\begin{frontmatter}
	
	\journal{Applied Mathematics and Computation}

%	\title{Polynomial Corrections for High Order Compressible Flows Simulations \\ on Curved Domains with Linear Meshes} 
	\title{{ Shifted boundary polynomial corrections for compressible flows:\\ high order on curved domains using linear meshes}}
	\cortext[mycorrespondingauthor]{Corresponding author}
	
	\author[Inria]{Mirco Ciallella\corref{mycorrespondingauthor}}
	\ead{mirco.ciallella@inria.fr}
	\author[Inria]{Elena Gaburro}
	\ead{elena.gaburro@inria.fr}
	\author[Inria]{Marco Lorini}
	\ead{marco.lorini@inria.fr}
	\author[Inria]{Mario Ricchiuto}
	\ead{mario.ricchiuto@inria.fr}

	\address[Inria]{Inria, Univ. Bordeaux, CNRS, Bordeaux INP, IMB, UMR 5251, 200 Avenue de la Vieille Tour, 33405 Talence cedex, France}

	%-------------------------------------------
	% ABSTRACT
	%
	\begin{abstract}
 In this work we propose a simple but effective high order polynomial correction allowing to enhance 
     the consistency of all kind of boundary conditions for the Euler equations (Dirichlet, characteristic far-field and slip-wall), both in 2D and 3D, preserving a high order of accuracy without the need of curved meshes.
      The method proposed is a simplified reformulation of  the Shifted Boundary Method (SBM) and relies on a \textit{correction} based on the \textit{extrapolated value} of the in cell polynomial \textit{to the true geometry}, thus not requiring  the explicit evaluation of high order Taylor series. Moreover, this strategy could be easily implemented into any already existing finite element and finite volume code. 
% 
%         has been developed and implemented to run compressible flow simulations characterized by curved domains
%        discretized using unstructured linear meshes. The formulation given by the Shifted Boundary Method (SBM) allows it to be fast to implement in any finite element framework
%        and easy to extend to three-dimensional configurations without any problem. The SBM, which usually would require cumbersome computations of high order partial derivatives,
%        has been replaced with a straightforward correction, which consists in an off-element evaluation of the polynomial, that takes into account all the high derivative terms in both 2D and 3D for polynomials of arbitrary order. 
%        This technique allows high order convergence rates without the need of generating high order curvilinear meshes to approximate the curved boundaries of the domain.
        Several validation tests are presented to prove the convergence properties up to order four for 2D and 3D simulations with curved boundaries, as well as an effective extension to flows with shocks.
	\end{abstract}	
	%------------------------------------------
	% KEY WORDS
	%
	%==========================================================================================
\end{frontmatter}

%%%%%%%%%%%%%%%%%%%%%%%%%%%%%%%%%%%%%%%%%%%%%%%%%%%%%%%%%%%%%%%%%%%%%%%%%%%%%%%%%%%%%%%%%%%%%%

%--------- SECTION ---------------------------------------------------------------------------
\section{Introduction} \label{sec.intro}
%---------------------------------------------------------------------------------------------

% \Rall{Please, write all the modifications concerning both the referees in this color.}

% \RIcolor{Please, write all the modifications concerning the first referee in this color.}

% \RIIcolor{Please write all the modifications concerning the second referee in this color.}

The potential of high order methods  in providing drastic error reductions in flow simulations with considerable savings in computational costs
is well established~\cite{Wang_et_al_ijnmf13}.  A key element in this respect is the proper treatment of \textit{boundary conditions}.
This involves two independent aspects: an appropriate geometrical representation of the boundaries and a high order approximation of the boundary condition itself.
When dealing with finite elements, the most classical approach is to work with an iso-parametric approximation in which the geometry as well as the flow solution are approximated by some high order polynomial~\cite{zienkiewic}. 
Standard approaches range from the use of various  maps  based on  some local interpolated or modal polynomial approximation of the curved geometry, to the more recent use of  rational B-spline or NURBS approximations used in the so called iso-geometric analysis (IGA)~\cite{nurbs0}.
In both cases, a crucial role is played by the availability of a high quality curved mesh. Progress has been made on methods to obtain such meshes, either via
curving  straight faced meshes~\cite{dey2001towards,luo2004automatic,sahni2010curved,FORTUNATO20161,MOXEY2016130,Puigt},
or by using some optimization or variational approach~\cite{https://doi.org/10.1002/nme.4888,TOULORGE20138,TURNER201873}. 
Despite the advances in this domain, obtaining easily a high quality curved mesh for a complex geometry remains a complex task,
still object of intense research both for dominant simplex meshes, and for coarse multiblock quad/hex meshes required in NURBS analysis (cf e.g.~\cite{doi:10.2514/6.2022-0389,doi:10.2514/6.2021-2991} and references therein).\\
 
In alternative to high order meshing, one can improve the boundary conditions by accounting, on a straight faced mesh, for the local features of the \textit{true} geometry.
Early work on curvature corrected wall boundary conditions can be found e.g. in~\cite{WangSun}, while the specific case of high order schemes and wall boundaries
in two space dimensions has been thoroughly treated in the well known paper by Krivodonova and Berger~\cite{krivodonova2006high}. 
The last reference in particular proposes an approach to correct the direction used when prescribing the slip-wall condition, which shows a recovery of third and, for some cases, 
fourth order of accuracy on 2D  geometries. However, this method can be formulated  only for slip-wall conditions, and the work is limited to 2D geometries.  
More recent developments  addressing this problem can be found in~\cite{costa2018very,costa2019very,fernandez2020very}.
The authors first presented a reconstruction off-site data (ROD) approach, in a high order finite volume framework, to apply Dirichlet boundary conditions~\cite{costa2018very}.
Afterwards, the work has also been extended to other boundary conditions~\cite{costa2019very}, which represents an important improvement for real applications.
Due to the finite volume framework, the approach is built on the least-squares method, which is used to handle several constraints from scattered mean values associated to the elements.
Due to the linear system arising from the constrained least-squares problem, a matrix should be inverted locally and this may introduce some issues for ill-conditioned problems.\\

In this work, we propose an alternative path to achieve a similar result. Our idea relies on  a simplified reformulation of the %This paper is going to focus on a new flux correction based on the
 Shifted Boundary Method (SBM)~\cite{Scovazzi1,Scovazzi2,Scovazzi3,Scovazzi4,lishifted,ciallella2022extrapolated,assonitis2022new,colomes2021weighted} that enhances
 the consistency of all kinds of boundary conditions for the Euler equations (Dirichlet, characteristic far-field and slip-wall), both in 2D and 3D.
% 
%to implement high order far-field and slip-wall boundary conditions for the 2D/3D Euler equations when the boundaries are approximated using conformal linear meshes.
%It should be noticed that with few adjustments the same approach can also be implemented for other hyperbolic systems.
Originally, %The 
SBM has been introduced to cope with embedded boundary problems and it consists in retaining the designed order of accuracy of the discretization method by using
modified boundary conditions based on truncated Taylor series expansions. Here we propose a \textit{simplified version} which allows to avoid the cumbersome explicit evaluation 
of all the Taylor development terms~\cite{Scovazzi4,atallah2022high}, and directly exploits the available (nodal or modal) polynomial expansion of a physical, linear element.
%Herein we exploit this capability to properly enforce weak boundary conditions when the physical boundary is not fitted by high order meshes.
The approach has been easily implemented for steady problems using indifferently nodal and modal bases on simplex meshes, 
and its extension to time dependent problems is based on the ADER space-time formulation.
The ADER methodology has been introduced by Toro and collaborators in~\cite{titarev2002ader,titarev2005ader,toro2006derivative} and it has been substantially simplified and generalized in~\cite{dumbser2008unified}. 
Then, the ADER approach has been widely used in the last 15 years: we recall here only some of its main developments, as the extension 
to staggered meshes~\cite{tavelli2017pressure,tavelli2018arbitrary}
to unstructured 2D and 3D moving meshes~\cite{loubere2014new, boscheri2015direct,gaburro2020high, gaburro2021unified},
to semi-implicit and implicit schemes~\cite{montecinos2014reformulations, busto2020efficient,han2021dec},
its coupling with the MOOD \textit{a posteriori} finite volume limiter~\cite{boscheri2018second,gaburro2021posteriori}
and its application in the modeling of quite complex PDE systems, as MHD in~\cite{balsara2009efficient,fambri2017space}, 
GRMHD and CCZ4 in~\cite{fambri2018ader,dumbser2018conformal,dumbser2020glm} and multiphase models and the recent first order hyperbolic unified model for continuum mechanics, known as GPR system, in~\cite{peshkov2021simulation,chiocchetti2021high,gabriel2021unified, boscheri2022cell}.
Although the extensive use of ADER schemes, the problems related to high accurate treatment of boundary conditions 
have been only marginally considered, see for example the work~\cite{kemm2020simple} that however treats the topic from a completely different approach employing indeed a diffuse interface technique.
Finally, for a formulation of ADER schemes employing a genuinely space-time \textit{modal expansion}, 
as the one adopted in this work, 
we refer to~\cite{gaburro2020high,gaburro2022high}.\\

The paper is organized as follows.
First, in Section~\ref{sec.model} we describe the physical model of interest in this work,
represented by the Euler equations of gasdynamic on curved domains. 
\RIcolor{Due to the generality of the approach, it should be noticed that the method proposed here can be easily applied to other systems of equations~\cite{Scovazzi1,Scovazzi2}.}
Next, in Section~\ref{sec.DG} we briefly recall the general framework of Discontinous Galerkin (DG) schemes,
used here to discretize the governing PDEs system with high order of accuracy in space and time on conformal meshes;
in particular, a brief overview on the standard treatment of boundary conditions,
involving the necessity of using curved meshes for curved domain, is given in Section~\ref{ssec.BCconformal}.
Then, the core of the paper is represented by Section~\ref{sec.BCpolynomial},
where we present the \textit{polynomial correction} to be applied to general high order numerical schemes
directly on \textit{linear} fitted meshes for modeling \textit{compressible} flows on \textit{curved} domains.
We would like to remark that the presented approach is simple, independent on the underlying type of high order discretization and on the space dimension
and allows to retrieve the formal order of accuracy of the original method for arbitrary domain and obstacles shapes.
Finally, a large set of numerical simulations, allowing to validate the proposed polynomial correction, is presented in Section~\ref{sec.test},
for both far-field and slip-wall boundary conditions, steady and unsteady problems, on two and three dimensional curved domains.
The paper is closed by some conclusive remarks and an outlook to future works is given in Section~\ref{sec.concl}.

%--------- SECTION ---------------------------------------------------------------------------
\section{{Governing equations}} \label{sec.model}
%---------------------------------------------------------------------------------------------

%\subsection{Governing equations} \label{ssec.model_pde}

We consider the numerical approximation of solutions of the Euler equations {in $d$ space dimensions} %in two dimensions 
reading:
\begin{equation}\label{eq:euler0}
\partial_t \mathbf{U}+\nabla\cdot\mathbf{F}\left(\mathbf{U}\right) =0, \quad \text{on}\quad \Omega_T=\Omega\times[0,T]\subset\mathbb{R}^d\times\mathbb{R}^+,
\end{equation}
{with  $\mathbf{U}$ the vector of conserved variables and $\mathbf{F}$ the non linear flux respectively defined as}
%where $\mathbf{U}$ is the solution of the system of conservation laws,
%$\mathbf{F}$ is the flux function and
%$\nabla:=\left(\partial_{x_1},\dots,\partial_{x_d}\right)$ is the gradient operator in the physical domain.\\
%Conserved variables and fluxes are given by:
\begin{equation}\label{eq:euler0a}
\mathbf{U}=\left(
\begin{array}{c}
\rho \\ \rho \mathbf{u}  \\ \rho E
\end{array}
\right)\;,\;\;
\mathbf{F}= \left(
\begin{array}{c}
\rho \mathbf{u}\\ \rho \mathbf{u} \otimes \mathbf{u} + p \mathbb{I} \\ \rho H \mathbf{u}
\end{array}
\right), 
\end{equation}
having denoted by $\rho$ the mass density, by $\mathbf{u}$ the velocity, by $p$ the pressure, and with $E=e+\mathbf{u}\cdot\mathbf{u}/2$ the specific total energy, $e$ being the specific internal energy.
Finally, the total specific enthalpy is  $H=h+\mathbf{u}\cdot\mathbf{u}/2$ with  $h=e+p/\rho$ the specific enthalpy. For simplicity in this paper we  work with the classical perfect gas
equation of state%:
\begin{equation}\label{eq:EOS}
p=(\gamma-1)\rho e,
\end{equation}
with $\gamma$ the constant % (for a perfect gas) 
ratio of specific heats.\\

\RIIcolor{In order to perform convergence analysis, we use the manufactured solution method which needs the discretization of an additional source term in Equation~\eqref{eq:euler0} which thus reads as follow 
\begin{equation}\label{eq:euler1}
\partial_t \mathbf{U}+\nabla\cdot\mathbf{F}\left(\mathbf{U}\right) = \mathbf{S}(\mathbf{U}).
\end{equation}
For this reason we are going to include $\mathbf{S}$ in the discretization presented in Section~\ref{sec.DG}.}

%\subsection{Boundary conditions} \label{ssec.model_bc}

%\subsection{Analytical solutions} \label{ssec.model_sol}	

\bigskip
%--------- SECTION ---------------------------------------------------------------------------
\section{High order Discontinuous Galerkin discretization} \label{sec.DG}
%---------------------------------------------------------------------------------------------

\subsection{Computational domain and data representation}

%\textcolor{red}{
%write down \\
%- Introduce the domain ... underline that the computational domain does not coincide with the physical domain		\\
%- Fix all but only the relevant notation\\
%- Conserved variables are discretized with high order polynomial in space (just one line to say that FV is a special case of that) .. cell centered approach}

The  spatial  domain $\Omega$ is discretized by means of a tessellation $\mathscr{T}$ composed of $\mathscr{N}$ non-overlapping simplicial  elements { (triangles in 2D, tetrahedra in 3D)}.
We denote by $K$ the generic element, {so that} %and set 
\RIIcolor{$\Omega_h=\bigcup_{j=1}^{\mathscr{N}} K_j$}. Note that  in general  $\Omega_h\ne \Omega$ {and in particular $\partial\Omega_h\ne\partial\Omega$ for  most approximations, even conformal, with the exception of very simple geometries or of iso-geometric approaches~\cite{nurbs0}.}

%% The numerical solution $\mathbf{U}$ is approximated {by $\mathbf{U}_h\in\pmb{\mathcal{V}}_h^p\times\mathbb{R}^{d+2}$,
%% where $\pmb{\mathcal{V}}_h^p$ is the space of piece-wise polynomials of degree $p$ within each triangle $K$  and discontinuous across faces. In particular, within each $K$ we set
%% $\mathbf{U}_h^n(\mathbf{x},t^n)\in(\pmb{\mathcal{V}}_h^p)_K=\text{span}\{\psi_i\}_{1\le i\le D}$, or in other words}
\RIIcolor{
The numerical solution $\mathbf{U}$ is approximated by $\mathbf{U}_h$, which belongs to a space of piece-wise polynomials within each triangle $K$ and discontinuous across faces, such that in each element $K_j$ we have
%
%
%% \begin{equation}\label{eq:basis psi}
%% \mathbf{U}_h(\mathbf{x},t) = \sum_{i=1}^D \mathbf{U}_i(t)\psi_i(\mathbf{x}).  %:= \mathbf{U}_i(t)\psi_i(\mathbf{x}) \,. %,\qquad \mathbf{v}_h(\mathbf{x}) = \sum_{i=1}^n \mathbf{V}_i(t)\psi_i(\mathbf{x}), \qquad \forall\mathbf{x}\in K.
%% \end{equation}
\begin{equation}\label{eq:basis psi}
\mathbf{U}_h(\mathbf{x},t)|_{K_j} = \sum_{i=1}^D \mathbf{U}_{k}(t)\psi_{k}(\mathbf{x}), %\quad \mathbf{x}\in K_j,\quad j=1,\ldots,\mathscr{N},
\end{equation}
where $\{\psi_{k},k=1,\ldots,D\}$ is a basis of polynomials of degree $p$.}
On simplex elements, the number of \textit{degrees of freedom} $D$ can be  shown to be $D=\prod_{l = 1}^d  (p + l)/l$. 
\Rall{The discontinuous finite element data representation~\eqref{eq:basis psi} leads to a Finite Volume (FV) scheme if $p=0$}.

%\subsection{High order discontinuous Galerkin method}
%
%These flux corrections have been embedded within two different implementations of discontinuous Galerkin (DG).
%For two and three dimensional steady state problems we use a  classical method of lines approach in which the semi-discrete equations obtained from the DG method are integrated by means of appropriate ODEs  methods.
%For time dependent problems, we have used an explicit ADER-DG method with local space-time predictors~\cite{dumbser2008unified,gaburro2020high}. When necessary,
%an \textit{a posteriori} sub-cell limiter provides control of the numerical solution in correspondence of discontinuities~\cite{loubere2014new,gaburro2021posteriori}. 
%
%The following subsections recall the main building blocks of these techniques, allowing to introduce the relevant notation used later in the paper.
%
\subsection{Discontinuous Galerkin discretization in space}
{The elemental semi-discrete discontinuous Galerkin (DG) weak formulation is classically written by projecting each component of~\RIIcolor{\eqref{eq:euler1}} 
on the relevant basis and integrating by parts~\cite{cockburn1998runge,bassi1997high}}. 
\begin{equation}\label{eq:DGweak semidiscrete 2}
\int_K \psi_i \, \frac{\diff{ \mathbf{U}_h}}{\diff{t}} \diff{\Rall{\mathbf{x}}} + \int_{\partial K} \psi_i \, \hat{F}(\mathbf{U}_h^-,\mathbf{U}_h^+)\cdot \mathbf{n}\diff{S}-\int_K \nabla\psi_i \, \cdot\mathbf{F}(\mathbf{U}_h)\diff{\Rall{\mathbf{x}}} = \RIIcolor{\int_K \psi_i \, \mathbf{S}(\mathbf{U}_h)\diff{\mathbf{x}}}, \qquad 1 \leq i \leq D,
\end{equation}
%
%The current work is performed within a discontinuous Galerkin (DG) framework~\cite{cockburn1998runge,bassi1997high}.
%To obtain semi-discrete equations, we replace $\mathbf{U}$ in~\eqref{eq:euler0} by a discrete approximation $\mathbf{U}_h\in\pmb{\mathcal{V}}_h^p\times\mathbb{R}^{d+2}$,
%where $\pmb{\mathcal{V}}_h$ is the  broken space of broken polynomials of degree $p$ within each triangle $K$,  and discontinuous across faces.
%The DG weak formulation is classically written as: find $\mathbf{U}_h\in \pmb{\mathcal{V}}_h^p\times\mathbb{R}^{d+2}$ such that $\forall K$, and  $\forall\mathbf{v}_h\in\pmb{\mathcal{V}}_h^p$ we have
%%
%\begin{equation}\label{eq:DGweak semidiscrete}
%\int_K \mathbf{v}_h\frac{\diff{ \mathbf{U}_h}}{\diff{t}} \diff{V} + \int_{\partial K} \mathbf{v}_h\hat{F}(\mathbf{U}_h^-,\mathbf{U}_h^+) \diff{S}-\int_K \nabla\mathbf{v}_h\cdot\mathbf{F}(\mathbf{U}_h)\diff{V}=0, \qquad \forall\mathbf{v}_h,
%\end{equation}
with $\hat{F}(\mathbf{U}_h^-,\mathbf{U}_h^+)$ a consistent numerical flux which depends on the face values of the internal  state $\mathbf{U}_h^-$, of the neighboring element  state 
$\mathbf{U}_h^+$ and on the face normal $\mathbf{n}$. We recall in particular that a consistent flux  is a Lipschitz continuous function of  each of its arguments which also veryfies
\begin{equation}\label{eq:constistencyDG}
\Rall{\hat{F}(\mathbf{u},\mathbf{u})\cdot\mathbf{n} = \mathbf{F}(\mathbf{u})\cdot\mathbf{n}} %,\quad \hat{F}(\mathbf{u},\mathbf{w})\cdot\mathbf{n} = - \hat{F}(\mathbf{w},\mathbf{u})\cdot\mathbf{n}. 
\end{equation}
In this work we  have used a simple and robust Rusanov-type flux:
\begin{equation}\label{eq:rusanov flux}
\RIcolor{\hat{F}}(\mathbf{U}_h^-,\mathbf{U}_h^+)\cdot \mathbf{n} = \frac{1}{2}\left(\mathbf{F}(\mathbf{U}_h^+) + \mathbf{F}(\mathbf{U}_h^-)\right)\cdot\mathbf{n}-\frac{1}{2}s_{max}\left(\mathbf{U}_h^+ \Rall{-} \mathbf{U}_h^-\right),
\end{equation}
where $s_{max}$ is the maximum eigenvalue of the Jacobians of the flux $A_{\mathbf{n}}\left(\mathbf{U}_h^+\right)$ and 
$A_{\mathbf{n}}\left(\mathbf{U}_h^-\right)$. %\\
%
%%
%The semi-discrete equations can be written as:  find $\mathbf{U}_h\in \pmb{\mathcal{V}}_h^p\times\mathbb{R}^{d+2}$ such that $\forall K$   we have
%%
%\begin{equation}\label{eq:DGweak semidiscrete 2}
%\int_K \psi_j\frac{\diff{ \mathbf{U}_h}}{\diff{t}} \diff{V} + \int_{\partial K} \psi_j \hat{F}(\mathbf{U}_h^-,\mathbf{U}_h^+)\cdot \mathbf{n}\diff{S}-\int_K \nabla\psi_j\cdot\mathbf{F}(\mathbf{U}_h)\diff{V}=0, \qquad 1 \leq j \leq D.
%\end{equation}
%
%On simplex elements, the number of \emph{degrees of freedom} $D$ can be  shown to be $D=\prod\limits_{l = 1}^d  (p + l)/l $.

\subsection{High order time integration}

In this Section we briefly recall two classical different strategies, both employed in this work, to achieve high order of accuracy during the time integration with a discontinous Galerkin scheme.

\subsubsection{Method of lines for steady problems} %Semi-discrete approach via Runge-Kutta ODE integrator}

For steady state computations, we use a classical method of lines. % approach.
Assembling  all the contributions~\eqref{eq:DGweak semidiscrete 2}, we obtain in each element $K$
a system of ODEs~\cite{bassi1997high} reading
\begin{equation}\label{eq:ODE semidiscrete}
\frac{\diff{\mathbf{U}}}{\diff{t}}+  M^{-1} \mathsf{R}(\RIIcolor{\mathbf{U}}) =0,
\end{equation}
where \RIIcolor{the array $\mathbf{U}$  contains the $d+2$ vector of unknowns corresponding to the degrees of freedom},
 $\mathsf{R}$ is the array  of size $D$, \RIIcolor{made up by} the second, third, and fourth integrals in~\eqref{eq:DGweak semidiscrete 2},
and   $M$ denotes the elemental mass matrix with  $D\times D$ block diagonal entries, \RIIcolor{ whose components are given by
%mass matrix
\begin{equation}\label{eq:massM}
[M]_{ik}= \int_K \psi_i\psi_k  \diff{\mathbf{x}}\,.
\end{equation}
}
%%
%The  $j$th component of the array $\mathbf{U}$  contains the $d+2$ vector of unknowns corresponding to the degree of freedom $\mathbf{U}_j$.
For steady problems, we simply integrate~\eqref{eq:ODE semidiscrete}  with explicit Euler or  classical Runge-Kutta methods with improved time step constraints~\cite{butcher2000numerical,gottlieb2001strong}.

\subsubsection{High order explicit ADER approximation with local  space-time predictors} % Fully discrete approach via ADER predictor corrector formulation}

For time accurate simulations, we have  used a high order explicit one-step \textit{predictor-corrector} ADER approach, \RIIcolor{that we are going to briefly describe in this section.
For further details we refer to ~\cite{dumbser2008unified,gaburro2020high}}.
The temporal domain is as usual discretized in temporal slabs  $[t^n, t^{n+1}]$.
With the same notation used before, on fixed meshes the ADER method  can be succinctly written as
\begin{equation}\label{eq:corrector2}
  \mathbf{U}^{n+1}  =    \mathbf{U}^{n} -  M^{-1} \int_{t^n}^{t^{n+1}} \mathsf{R}(\mathbf{q}_h) dt, % 
\end{equation}
where $\mathbf{q}_h=\mathbf{q}_h(\mathbf{x},t)$ is  a high order space-time predictor of the solution in the time slab $[t^n, t^{n+1}]$.
In practice,~\eqref{eq:corrector2} is replaced by some high order quadrature formula in time (with $\alpha_i$ and $\omega_i$, $1 \le i \le r$, the quadrature points and weights)
\begin{equation}\label{eq:corrector2a}
  \mathbf{U}^{n+1}  =    \mathbf{U}^{n} -  M^{-1}  \Delta t 
  \sum\limits_{i = 1}^{r}\omega_i \mathsf{R}\big(\mathbf{q}_h(t^{n+\alpha_i}) \big), % 
\end{equation}
and    $\mathbf{q}_h(\mathbf{x},t)$ is defined as a polynomial of degree $N$ in space \emph{and} time, with $N \ge p$.
This polynomial  is obtained by means of a genuinely local space-time procedure.
In the current implementation this local problem  is formulated  by means of a modal  expansion
\begin{equation}\label{eq:predictor0}
\mathbf{q}_h(\mathbf{x},t) = \sum_{\ell=0}^{\mathcal{Q}-1}\theta_\ell(\mathbf{x},t)\mathbf{q}^n_\ell, \qquad (\mathbf{x},t)\in \RIIcolor{K\times[t^n,t^{n+1}]}, \qquad \mathcal{Q}=\mathcal{L}(N,d+1),
\end{equation}
where  $\theta_\ell(\mathbf{x},t)$ \RIIcolor{being $\mathcal{L}(N,d)=\prod_{m=1}^d(N+m)/m$ space-time modal basis of the polynomials of degree $N$ in $d+1$ dimensions ($d$ space dimensions and time)}, defined as
\begin{equation}
\theta_\ell(x,y,t)\RIIcolor{|_{K\times[t^n,t^{n+1}]}} = \frac{(x-x_{K}^n)^{p_\ell}}{p_\ell ! h_\RIcolor{K}^{p_\ell}} \frac{(y-y_{K}^n)^{q_\ell}}{q_\ell ! h_K^{q_\ell}} \frac{(t-t^n)^{r_\ell}}{r_\ell ! h_K^{r_\ell}}, \quad \ell=0,\ldots,\mathcal{L}(N,d+1),\quad 0\leq p_\ell+q_\ell+r_\ell\leq N,
\end{equation}
with $(x_K,y_K)$ the gravity center of element $K$, \RIIcolor{and $(p_\ell,q_\ell,r_\ell)$ being the terms associated to the taylor series expansion}. 
The values \RIIcolor{of the predictor $\mathbf{q}^n_h$} are %obtained as the solution of the local space-time weak formulation
\RIIcolor{computed by means of an iterative procedure that seeks the solution for any space-time element $K\times[t^n,t^{n+1}]$} of the local space-time weak formulation  
\begin{align}
&\int_{K} \theta_k(\mathbf{x},t^{n+1}) \mathbf{q}^n_h(\mathbf{x},t^{n+1})\diff{\mathbf{x}} - \int_{K} \theta_k(\mathbf{x},t^{n}) \mathbf{U}^n_h(\mathbf{x},t^{n})\diff{\mathbf{x}} - \nonumber\\\qquad &\int_{t^n}^{t^{n+1}}\int_{K} \frac{\partial \theta_k}{\partial t}(\mathbf{x},t)\mathbf{q}^n_h(\mathbf{x},t)\diff{\mathbf{x}}\diff{t}  +\int_{t^n}^{t^{n+1}}\int_{K} \theta_k(\mathbf{x},t)\nabla\cdot\mathbf{F}(\mathbf{q}^n_h)\diff{\mathbf{x}}\diff{t}= \RIIcolor{\int_{t^n}^{t^{n+1}}\int_{K} \theta_k(\mathbf{x},t) \, \mathbf{S}(\mathbf{q}^n_h)\diff{\mathbf{x}}\diff{t}}, \label{eq:predictor2}
\end{align}
where $\mathbf{U}^n_h$ is the known initial condition at time $t^n$.\\

Equation~\eqref{eq:predictor2} is fully local, in the sense that it involves no
communication between $K$ and its neighbouring control volumes. The solution to this equation can be obtained \textit{independently within each element $K$} by means of some iterative procedure.
In this work  a simple discrete Picard   iteration for each space-time element is used.

\subsection{{\textit{A posteriori} sub-cell finite volume limiter and CFL constraint}}

%\begin{figure}
%\centering
%\subfigure[$N=1$]{\includegraphics[width=0.25\textwidth]{subcell1.pdf}}\quad
%\subfigure[$N=2$]{\includegraphics[width=0.25\textwidth]{subcell2.pdf}}\quad
%\subfigure[$N=3$]{\includegraphics[width=0.25\textwidth]{subcell3.pdf}}
%\caption{Sub-triangulation of $K$ for polynomials of degree $N$.}\label{fig:subcell}
%\end{figure}
%
%The DG scheme presented so far only consists in representing the solution through high order linear (Godunov theorem) polynomials. 
%However, linear polynomials are affected by the Gibbs phenomenon, i.e.\ oscillations are
%likely to appear where shock waves or other discontinuities occur.
%Therefore a limiting technique is required. The strategy employed herein is based on the 

{To handle discontinuous solutions in the time dependent case we use the }
MOOD approach~\cite{clain2011high,diot2012improved,diot2013multidimensional}, 
which has already been effectively applied in the ADER framework~\cite{loubere2014new,boscheri2017high,gaburro2021posteriori}.
 
The algorithm is based on an \textit{a posteriori} technique. 
{The solution is first evolved % meaning that the solution will be first evolved in time 
from  $t^n$ to $t^{n+1}$ using the high order ADER-DG method. Then several admissibility criteria are checked and 
the solution in all troubled cells (i.e. the cells not satisfying the admissibility criteria) is recomputed  \textit{a posteriori} using a MUSCL-Hancock TVD finite volume (FV) scheme, but working on a sub-triangulation of the initial grid in order to preserve the accuracy of the high order DG scheme also when passing to a lower order but more robust second order FV scheme.
All aspects of  the implementation of this technique are provided in~\cite{gaburro2021posteriori} to which we refer for details. }\\

{Concerning the choice of the time step, we have implemented the usual } % the Regarding the stability region, the method is stable as long as the time step size $\Delta t$ satisfies an 
explicit CFL condition %, given by
\begin{equation}
\label{eqn.cfl}
\Delta t < CFL \left( \frac{|h_{\text{min}}|}{(2N+1)|\lambda_{\text{max}}|}  \right),
\end{equation}
where $|h_{min}|$ is the minimum characteristic mesh-size and $|\lambda_{max}|$ is the maximum eigenvalue of the Jacobian of the flux.
{For DG on} unstructured meshes the CFL stability condition requires the inequality $CFL<1/d$ to be satisfied {(cf. discussion in~\cite{CHALMERS2020109095}).
We underline that the time step constraint does not need to be  modified in presence of troubled cells, because we subdivide each troubled triangle in exactly $(2N+1)^d$ sub-triangles and then we employ a FV scheme for which~\eqref{eqn.cfl} holds with $N=0$. 
We refer  again to~\cite{gaburro2021posteriori}  for details.}
% 
%It should be  also noticed that for the sub-cell FV scheme we have a different CFL stability condition
%%
%\begin{equation}
%\Delta t_{\text{FV}} < CFL_{\text{FV}} \;\frac{|h_{\text{min}}|}{N_k}\frac{1}{|\lambda_\text{max}|},
%\end{equation}
%%
%with $CFL_{\text{FV}}<1/d$ and $|h_{\text{min}}|$ the minimum cell size referred to $K$.

\subsection{{Boundary conditions on  conformal meshes}} \label{ssec.BCconformal}

When the boundary $\partial K$ of element $K$ belongs to $\partial\Omega_h$, the normal flux function $\mathbf{F}(\mathbf{U}_h)\cdot\mathbf{n}$ {must  account for the appropriate  boundary conditions.
The %numerical 
flux consistent with such conditions will be denoted by}
% consistent with the boundary condition that has to be imposed	
%on $\partial\Omega$ and will be referred to as 
$\mathbf{F}^{bc}_{\mathbf{n}}$. In this work, the boundary flux function $\mathbf{F}^{bc}_{\mathbf{n}}$ is {obtained}
%generally computed
 by {defining a ghost state $\mathbf{U}^{bc}$,} %adding a row of ghost cells behind $\partial\Omega_h$ 
 and {introducing a numerical flux $\mathbf{F}^{bc}_{\mathbf{n}} = \hat{F}(\mathbf{U}^-_h,\mathbf{U}^{bc})$ defined by some approximate Riemann solver based on the internal state }
%  solving
%a Riemann problem with a numerical flux function $\mathbf{F}^{bc}_{\mathbf{n}} = \hat{F}(\mathbf{U}^-_h,\mathbf{U}^{bc})$, where 
$\mathbf{U}^-_h$ {and on} % is the internal state and 
$\mathbf{U}^{bc}$. {Depending on the condition to be enforced, different definitions  of the ghost state are used: } % is chosen depending
%on 
%the condition that must be enforced. % at the boundary.
\begin{itemize}
\item for far-field, which can be seen as a Dirichlet-type BC enforced weakly through fluxes, {all the components of } $\mathbf{U}^{bc}$ {are set to prescribed values}; % is set equal to the condition one wants to impose;
\item at inflow/outflow boundaries %the state 
$\mathbf{U}^{bc}$ is obtained {by imposing    the} % through a local one-dimensional analysis of characteristic lines in the direction normal to the boundary by
%imposing 
Riemann invariants associated to {characteristics entering the domain  values obtained from prescribed reference values of density, pressure and Mach number}; %outgoing characteristics and the variables we want to impose.
\item {for slip walls}  we wish to set
\begin{equation}\label{eq:slip}
\mathbf{u}\cdot\mathbf{n}=w.
\end{equation}
In this case $\mathbf{U}^{bc}$ has the same density, internal energy and tangential velocity of $\mathbf{U}^-_h$, and the opposite normal relative velocity component $\mathbf{u}\cdot\mathbf{n}-w$.
%It should be noticed that, when solving the Riemann problem, only pressure contributes to the boundary flux. For this reason, 
{For this condition, an alternative way } % another way of implementing slip-wall boundary conditions 
consists in defining {directly}  the boundary flux function {by setting % such that 
$\mathbf{u}\cdot\mathbf{n}=w$, meaning that 
\begin{equation}\label{eq:f_wall}
\mathbf{F}^{bc}_{\mathbf{n}} = \mathbf{F}_{\mathbf{n}}(\mathbf{U}^{bc})=  
\left(\rho w,\ \rho\mathbf{u} \,w + p\,\mathbf{n},\ \rho H \, w \right)\,.
\end{equation}
When doing so, the value of the total enthalpy should be consistently  modified as 
\begin{equation}\label{eq:sbm enthalpy0}
 H = \frac{\gamma}{\gamma-1}\frac{p}{\rho}+\frac{1}{2}\left(w^2 + u_{\mathbf{t}}^2\right),
\end{equation}
with $u_\mathbf{t}=\mathbf{u}\cdot\mathbf{t}$ the tangent velocity.
For static   walls this reduces to  $w=0$ and $\mathbf{F}^{bc}_{\mathbf{n}}=\left(0,\  p\,\mathbf{n},\  0 \right)$.}
\end{itemize}%MARIUZ

\begin{figure}
	\centering
	\includegraphics[width=0.6\textwidth]{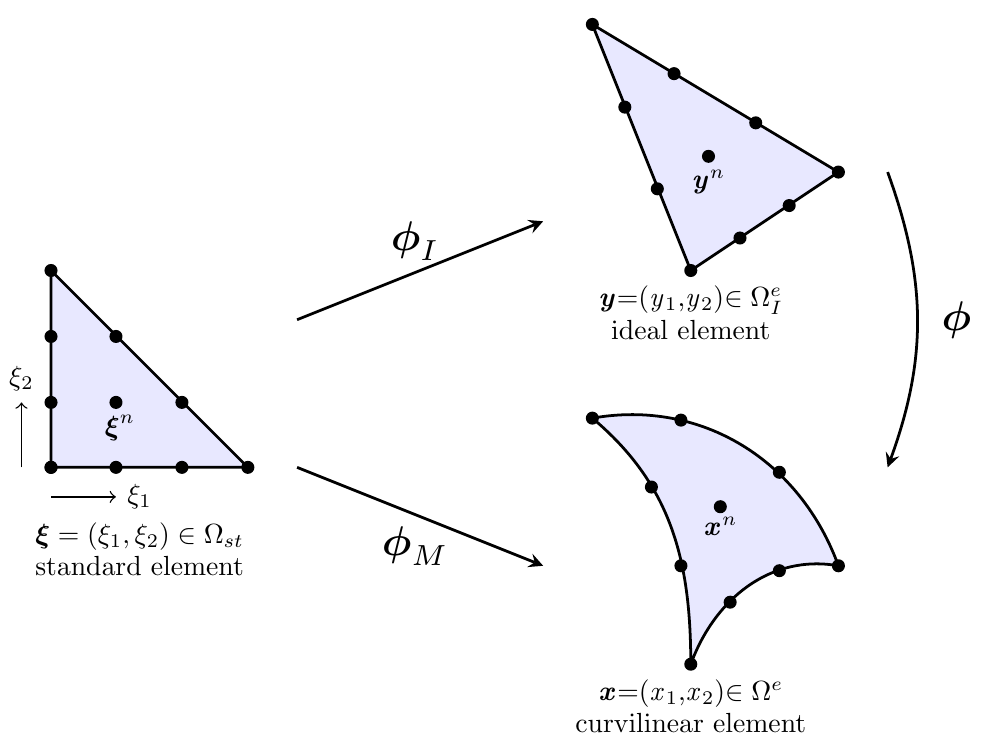}
	\caption{We give here a visual idea of the necessary transformations to be employed when dealing with curved meshes. 
		First, one needs the map from a standard reference element $\Omega_{st}$ (left) onto the straight-sided element $\Omega^e_{I}$ (top right), i.e. the mapping $\boldsymbol\phi_I:\Omega_{st}\rightarrow\Omega^e_I$, 
		and the one onto the curvilinear element (bottom right), i.e. $\boldsymbol\phi_M:\Omega_{st}\rightarrow\Omega^e$.
		Finally, the deformation mapping $\boldsymbol\phi:\Omega^e_I\rightarrow\Omega^e$ is defined through the composition $\boldsymbol\phi=\boldsymbol\phi_M\circ\boldsymbol\phi_I^{-1}$.}
	\label{fig:curvedElement}
\end{figure}

For high order methods, one of the key aspects in order to achieve a genuinely high order of accuracy is represented by the ability of simultaneously control the  error on the geometry and the flow variables. 
\RIcolor{This also includes approximating boundary integrals using consistent quadrature rules.}
To this end, an obvious, but also complex and expensive, solution consists in the use of a 
curved high order approximation of the boundary of the domain, 
which usually entails the use of some iso-parametric approximation of the boundary and the generation 
of a valid curved volume mesh~\cite{Wang_et_al_ijnmf13}.   
Curvilinear grids represent geometric boundaries with far superior accuracy 
allowing the use of larger elements than would be possible with linear elements.
Several approaches exist to obtain valid high order meshes, either based on  curving existing linear meshes~\cite{dey2001towards,luo2004automatic,sahni2010curved,FORTUNATO20161,MOXEY2016130},
or on some optimization or variational method~\cite{https://doi.org/10.1002/nme.4888,TOULORGE20138,TURNER201873}.
Moreover, these approaches always require the definition and evaluation of mappings between the curvilinear elements and the reference elements (see Figure~\ref{fig:curvedElement}), and despite the recent progress, 
while the generation of linear meshes for complex geometries has reached a very high level of maturity, 
the robust generation of curved meshes remains a relatively complex issue.

In this work we employ linear and curved meshes generated by the open source mesh generator {\tt Gmsh}~\cite{https://doi.org/10.1002/nme.2579}, 
and in the next Section, which represents the core of this work, 
we will show a novel technique to retrieve high order of accuracy for curved domain using linear meshes.

\section{{High order boundary conditions  on linear meshes via polynomial corrections}}\label{sec.BCpolynomial}
%---------------------------------------------------------------------------------------------

{In this work we aim at side stepping the need of generating curved meshes allowing to use directly conformal linear ones.
The idea is to use a simplified formulation of the }
Shifted Boundary Method (SBM), {originally introduced to handle non-conformal meshes within second order of accuracy for elliptic, parabolic and hyperbolic problems~\cite{Scovazzi1,Scovazzi2,Scovazzi3}, in order to 
compensate geometrical errors and retain the desired high order of accuracy.}
%the imposition of boundary conditions for curved geometries when high order finite element discretization are employed.
%The originality of the SBM lies in the idea of \textit{shifting} the location where the boundary conditions are applied.
%Originally, this method was proposed as a technique for solving problems with embedded geometries, however the latter and the problem
%discussed herein are much closer than one may think.
%In order to guarantee consistency, and retain the mesh convergence rates of the original method,
%the boundary conditions have to be modified.\\
%The main \MR{elements of the SBM} steps of the method are the following.

\begin{figure}[!tb]%
\centering
\subfigure{\includegraphics[width=0.45\textwidth]{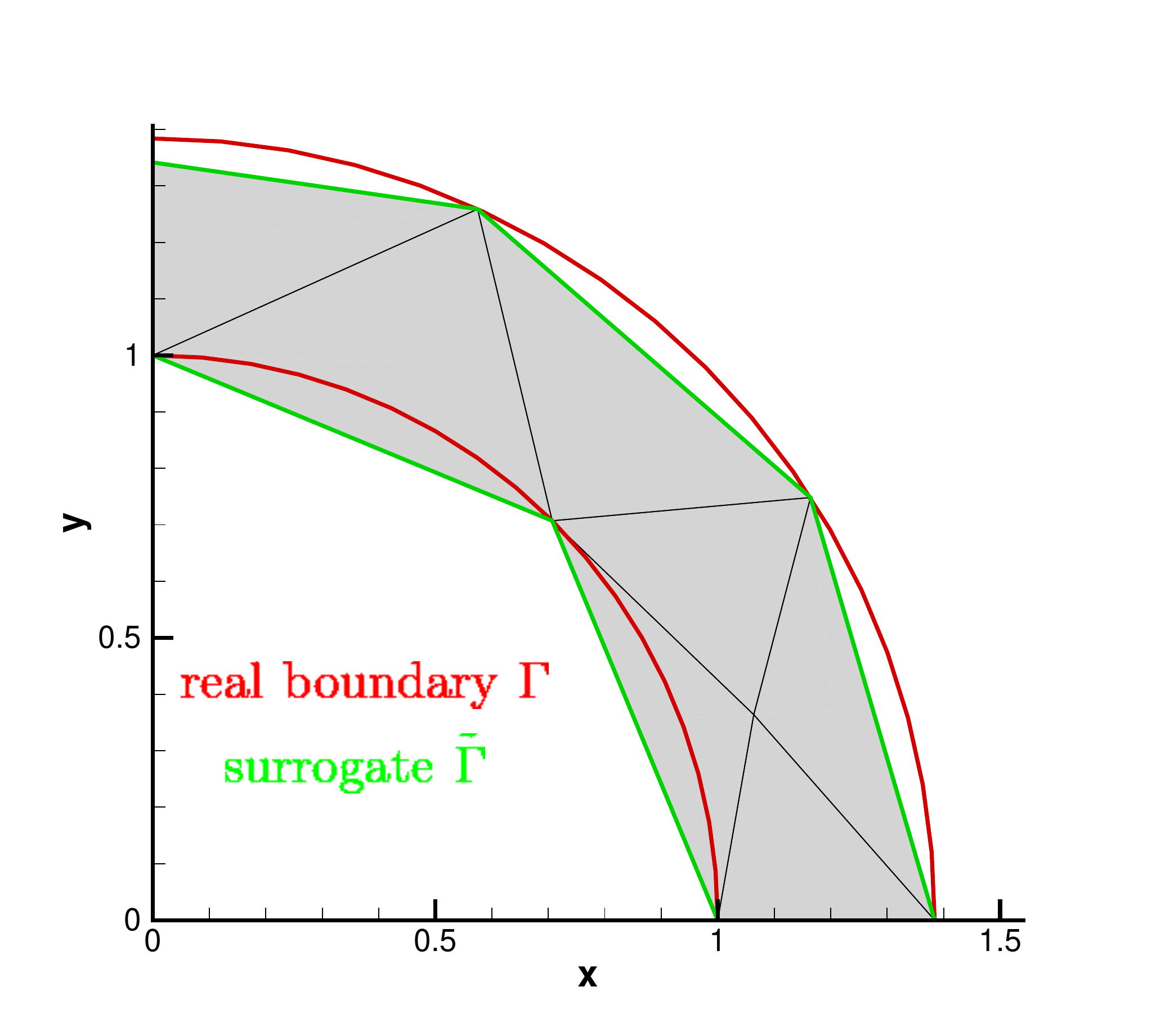}}\quad
\subfigure{\includegraphics[width=0.35\textwidth]{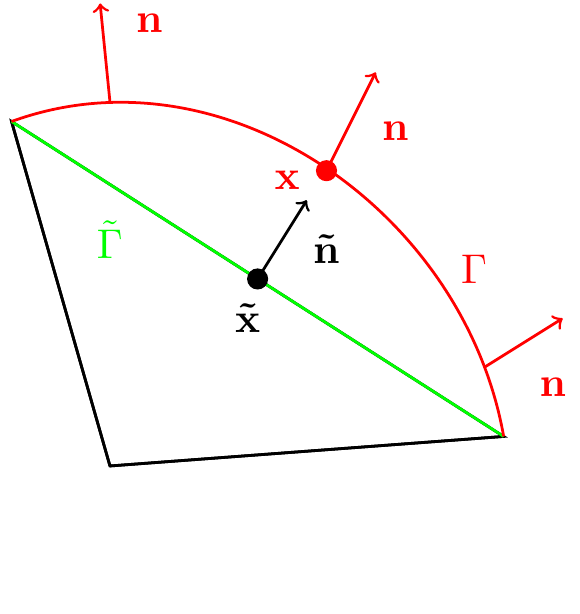}}
\caption{The SBM: the surrogate and actual boundaries with correspondent normals.}\label{fig:SBM1}%
\end{figure}

For the sake of clarity, we recall briefly the notation employed at boundaries. 
{Let  $\Omega_h$ be} a \textit{linear} conformal mesh discretizing the physical domain $\Omega$, {and $\tilde\Gamma:= \partial\Omega_h $ the \textit{linear} approximation of the}
%along with the 
curved boundary $\Gamma=\partial\Omega$. 
In particular, we refer to $\tilde\Gamma$ as to the 
\textit{surrogate boundary}. 
Moreover, {for any point on} % of the surrogate boundary 
$\tilde\Gamma$ {we assume to}
%one needs to 
be able to define a map to a unique point on the \textit{true boundary} $\Gamma$, s.t. 
\begin{align*}
\RIcolor{\mathcal{R}}\,:\,\tilde{\Gamma}\,&\rightarrow\,\Gamma \\
\mathbf{\tilde{x}}\,&\rightarrow\,\mathbf{x}.
\end{align*}
%
%which maps $\mathbf{\tilde{x}}\in\tilde{\Gamma}$  on the surrogate boundary to
% $\mathbf{x}\in\Gamma$  on the true boundary. 
The map $\RIcolor{\mathcal{R}}$ can be built in several ways, for example using a closest point projection, or  using level sets, or equivalently  using distances along directions normal to the true boundary $\Gamma$,
as shown in Figure~\ref{fig:SBM1}.
Since the gap between $\tilde{\Gamma}$ and $\Gamma$ is going to be of crucial importance, in terms of accuracy of the solution,
the map $\RIcolor{\mathcal{R}}$ will be characterized by a distance vector function:
\begin{equation}\label{dmap}
\mathbf{d(\tilde{x})}\,=\,\mathbf{x}\,-\,\mathbf{\tilde{x}}\,=\,[\RIcolor{\mathcal{R}}-\mathbf{I}](\mathbf{\tilde{x}}).
\end{equation}
\RIcolor{In our case, ${\mathcal{R}}$ is built using distances along  normals to   $\Gamma$, the vector  $\mathbf{d(\tilde{x})}$ is parallel to $\mathbf{n}(\mathbf{x})$ 
(the vector normal to $\Gamma$ in $\mathbf{x}=\RIcolor{\mathcal{R}}(\mathbf{\tilde{x}})$) such that
$$ \mathbf{x} = \mathbf{\tilde{x}} + \mathbf{d(\tilde{x})},\quad \mathbf{d(\tilde{x})}=\mathbf{n}(\mathbf{x})\|\mathbf{d(\tilde{x})}\|.$$
}

{Following the shifted boundary approach, we now modify the boundary conditions to retain the appropriate consistency order. 
We give a detailed description of the basic method for \textit{Dirichlet} conditions.}
Let $\rho_D$ be the prescribed value of the density (similar expressions can be written for all variables). The main idea is that a smooth exact solution of the problem  will verify the estimate
\begin{align}~\label{eq:taylor polynom}%\nonumber
\rho_D(\mathbf{\tilde x+d}) = \rho(\mathbf{\tilde x}) \;
 +& \;\|\mathbf{d}\|     \sum_{j=1}^d n_j \; \partial_{x_j} \rho(\mathbf{\tilde x}) % \text{ second order } \nonumber\\
 +\|\mathbf{d}\|^2   \sum_{j=1}^d \sum_{k=1}^d \frac{1}{2!}\; n_j n_k \; \partial^2_{x_j x_k} \rho(\mathbf{\tilde x}) \\ %\text{ third order }\\
 +& \;\|\mathbf{d}\|^3   \sum_{j=1}^d \sum_{k=1}^d \sum_{\ell=1}^d \frac{1}{3!}\; n_j n_k n_\ell \;\partial^3_{x_j x_k x_\ell} \rho(\mathbf{\tilde x})  + \dots \quad . \nonumber%  \mathcal{O}(\|\mathbf{d}\|^4)  \nonumber%\text{ fourth order }\nonumber
% &-\|\mathbf{d}\|^4   \sum_{j=1}^d \sum_{k=1}^d \sum_{\ell=1}^d \sum_{m=1}^d \frac{1}{4!}\; n_j n_k n_\ell n_m\; \partial^4_{x_j x_k x_\ell x_m} \rho &\text{ fifth order }
\end{align}
{The idea is thus to  modify the boundary condition  on $\tilde\Gamma$ to account for all the corrective terms, which boils down to use a modified prescribed value $\rho_{\text{SBM}}$ which, for different accuracy orders, is given by}
%has to be imposed as a Dirichlet boundary condition, the modified boundary condition would be
%
\begin{align}\label{eq:rho Taylor}
\rho_{\text{SBM}}(\mathbf{\tilde x}) = \rho_D(\mathbf{\tilde x+d})
 &-\|\mathbf{d}\|     \sum_{j=1}^d n_j \; \partial_{x_j} \rho(\mathbf{\tilde x}) & \text{ second order }\\
 &-\|\mathbf{d}\|^2   \sum_{j=1}^d \sum_{k=1}^d \frac{1}{2!}\; n_j n_k \; \partial^2_{x_j x_k} \rho(\mathbf{\tilde x}) &\text{ third order }\nonumber\\
 &-\|\mathbf{d}\|^3   \sum_{j=1}^d \sum_{k=1}^d \sum_{\ell=1}^d \frac{1}{3!}\; n_j n_k n_\ell \;\partial^3_{x_j x_k x_\ell} \rho(\mathbf{\tilde x}) &\text{ fourth order }\nonumber\\
 &\dots \ .  &\nonumber
% &-\|\mathbf{d}\|^4   \sum_{j=1}^d \sum_{k=1}^d \sum_{\ell=1}^d \sum_{m=1}^d \frac{1}{4!}\; n_j n_k n_\ell n_m\; \partial^4_{x_j x_k x_\ell x_m} \rho &\text{ fifth order }
\end{align}
{Note that for an explicit method all the terms involved in the right hand side of the last expression are known. For implicit schemes, they will modify the structure of the algebraic equations obtained.} 
\RIcolor{Also note that all the derivative terms  are evaluated  starting  from the available finite element approximation, and sampled at the appropriate quadrature points. Please refer to \cite{atallah2022high,Scovazzi4} for more details.}

%The last expression can be used as a modified boundary condition on the surrogate boundary of Dirichlet type $\tilde\Gamma_D$, which in this case preserves the accuracy up
%to fourh order.\\

%\subsection{Polynomial correction without gradient computation for Dirichlet boundary conditions} \label{ssec.Trick-Dir}
\subsection{{Derivative free formulation via polynomial corrections}} \label{ssec.Trick-Dir}

{The correction  terms in~\eqref{eq:rho Taylor}  become more and more cumbersome and costly  as the order of accuracy is increased,  especially in three space dimensions.}
%Moreover, the evaluation of these terms becomes even more involved when the finite element implementation systematically relies on a mapping in reference space, as shown in~\ref{appendix1}.}
%
% . Some details on the required operations to obtain
%the correction terms are }
%  The only implementation issue of~\eqref{eq:rho Taylor} would be of course the coding effort needed to compute all partial derivatives, which could really be cumbersome when dealing with very high order methods. Moreover, when the basis functions are defined in the physical space, the procedure consists in 
%computing the partial derivatives, which may already not be straightforward especially for 3D meshes. While, when the basis functions are defined in 
%the reference space, the procedure is even more cumbersome because of the different reference system. More information about this procedure can be
%found in~\ref{appendix1}.
{We propose here a different formulation that somehow \textit{simplifies} the evaluation of these terms, especially on straight-sided simplicial meshes for which both nodal and modal bases can be easily evaluated in physical space, without the needing of a map to the reference space, which would be instead required for curved elements.}

%We herein introduce a new formulation of the SBM that allows to overcome this.
\RIIcolor{Starting from a Taylor series expansion of arbitrary order of accuracy, for a generic variable $P$, centered in $\mathbf{\tilde{x}}$
$$ P(\mathbf{\tilde x+d}) = P(\mathbf{\tilde{x}}) + \nabla P(\mathbf{\tilde{x}})\cdot\mathbf{d} + \frac{1}{2}\mathbf{d^\intercal}\cdot\mathscr{H}(P(\mathbf{\tilde{x}})) \cdot\mathbf{d} + \dots \quad, $$
and moving $P(\mathbf{\tilde{x}})$ in the left-hand-side of the equation, we are left with
%of Eq.~\eqref{eq:taylor polynom}
%and recast it so as to have an alternative relation for all partial derivative terms:
%
$$P(\mathbf{\tilde x+d}) - P(\mathbf{\tilde{x}}) = \nabla P(\mathbf{\tilde{x}})\cdot\mathbf{d} + \frac{1}{2}\mathbf{d^\intercal}\cdot\mathscr{H}(P(\mathbf{\tilde{x}})) \cdot\mathbf{d} + \dots \quad. $$
Therefore, by simply evaluating the polynomial $P$ in $\mathbf{\tilde x+d}$ and by calculating the difference with that evaluated in $\mathbf{\tilde{x}}$
we get all {the necessary correction terms in \textit{only one} polynomial evaluation.} \Rall{A more interesting way to see our approach is to realize that 
the Taylor series expansion truncated to the same degree as the underlying elemental polynomial corresponds to a change of basis for the finite element space, 
passing to the local point wise Taylor basis. In this respect, the original SBM uses a different basis in different quadrature points, which is expensive and unnecessary.

In practice,  the data required in  all quadrature points can be evaluated using the unique  basis available in the implementation, whatever that may be.
Note also that the use of linear meshes, makes the finite element mapping fully linear, thus passing from the reference to the physical space can be done with no ambiguity for any
type of basis, and for all degrees of approxmation.}\\

For example, in practice the modified Dirichlet condition for the density $\rho$, introduced in~\eqref{eq:rho Taylor}, can be rewritten as
} 
%needed derivative terms up to the order of the accuracy of the polynomial used to discretize the solution.
%Finally we can recast the SBM flux correction~\eqref{eq:sbm example} of the polynomial to apply on $\tilde\Gamma$ as:
%
%\begin{equation}\label{eq:SBM trick}
%p_{SBM}(\mathbf{\tilde{x}}) = g(\mathbf{\tilde x+d}) - \left[ p(\mathbf{\tilde x+d}) - p(\mathbf{\tilde{x}}) \right]. 
%\end{equation}
%%
%When rewriting Eq.~\eqref{eq:rho Taylor} by using the trick of Eq.~\eqref{eq:SBM trick}, the density relation would simply be
%
\begin{equation}\label{eq:Dirichlet trick}
 \rho_{\text{SBM}}(\mathbf{\tilde x}) = \rho_D(\mathbf{\tilde x+d}) - \left[ \rho(\mathbf{\tilde x+d}) - \rho(\mathbf{\tilde x}) \right]
 {= \rho(\mathbf{\tilde x}) +  \left[ \rho_D(\mathbf{\tilde x+d}) - \rho(\mathbf{\tilde x+d})  \right],
 }
\end{equation}
which shows that the correction of the SBM method can be also seen as a direct \textit{shift} on the surrogate boundary
of the extrapolated polynomial error on the true boundary.  

This much simpler formulation only requires one extra polynomial evaluation
and thus it can be readily implemented on straight sides simplex elements
for which the basis functions are easily expressed in the physical space.
As we mentioned above, these extrapolated variables are then used to compose a \textit{ghost state} $\mathbf{U}^{bc}$ that will be used, along with
the internal state $\mathbf{U}_h^-$, as input for the numerical flux $\hat{F}(\mathbf{U}_h^-,\mathbf{U}^{bc})$ to obtain the consistent boundary flux $\mathbf{F}^{bc}_{\mathbf{n}}$.

\subsection{Treatment of slip wall boundary conditions} \label{ssec.Wall}

A similar approach can be applied in order to {impose the slip wall boundary condition~\eqref{eq:slip} \RIcolor{on the surrogate boundary $\tilde{\mathbf{x}}$}. An important issue to take into account in this case is that,
besides the position of the surrogate  wall boundary $\tilde\Gamma_w$, also its normal}
% compute the corresponding wall boundary condition that should be imposed on the real boundary $\Gamma_N$ but, instead,
%is going to be applied to the surrogate one $\tilde\Gamma_N$. In order to properly design the right boundary condition, it must be taken into account
%that the normal to the surrogate boundary 
$\tilde{\mathbf{n}}$ does not coincide with {$\mathbf{n}$, the normal to the true wall boundary $\Gamma_w$}. 
This difference {affects both the magnitude and the rate of the convergence of the error}
%plays an
%important role in the imposition of the boundary conditions. Indeed, when working with slip-wall boundary condition, for unviscous flows, the scalar product
%between the speed vector and the true normal should be set to zero:
%%
%\begin{equation}\label{eq:wall un}
%u_n=\mathbf{u}\cdot\mathbf{n}=0 \qquad\text{ on } \Gamma_N.
%\end{equation}
%%
%Hence, when mistaking the real boundary $\Gamma_N$ for the computational one $\tilde\Gamma_N$, the normal flux function reads $\mathbf{F}^{bc}_{\mathbf{\tilde n}} =  \left(0\quad p\,\mathbf{\tilde n}\quad 0 \right)$,
%which corresponds to set $u_{\tilde n}=\mathbf{u}\cdot\mathbf{\tilde n}=0$ rather than~\eqref{eq:wall un}. However, by looking at Figure~\ref{fig:SBM1}, it is clear than 
%this approximation yields to the introduction
%of spurious disturbances that will then affect both the discretization error and convergence rate
{as shown e.g.\ in~\cite{krivodonova2006high,bassi1997high,bassi1995accurate}. To overcome this issue, here we start from the formulation used in
~\cite{Scovazzi3}}.
%The procedure employed herein to overcome this limitation comes from the formulation of the SBM proposed in~\cite{Scovazzi3},
%in the shallow water framework, which will be recalled hereafter.
We start by decomposing the unit normal vector $\tilde{\mathbf{n}}$ at $\tilde{\mathbf{x}}$ as
\begin{equation}
\nonumber
\mathbf{\tilde n} = (\mathbf{\tilde n \cdot n})\,\mathbf{n}\,+\,\sum_{k=1}^{d-1}(\mathbf{\tilde n \cdot t}_k)\,\mathbf{t}_k,
\end{equation}
where $\mathbf{t}_k$ are the vector tangent to $\Gamma_w$.
By doing so, $\mathbf{F}^{bc}_{\mathbf{\tilde n}}$ can be recast as
\begin{equation}\label{eq:BC integral}
\mathbf{F}^{bc}_\mathbf{\tilde n} = (\mathbf{\tilde n \cdot n})\mathbf{F}_{\mathbf{n}}  + \sum_{k=1}^{d-1} (\mathbf{\tilde n \cdot t}_k)\mathbf{F}_{\mathbf{t}_k}.
\end{equation}
Then we can apply the Taylor expansion to {the normal velocity appearing in the flux terms, \RIcolor{while the other terms like $\rho$ and $p$ are taken from $\tilde{\mathbf{x}}$ without any corrections},
\begin{equation}\label{eq:SBM fluxes}
\mathbf{F}_{\mathbf{n}} = \left(\begin{array}{c} \rho w_{\text{SBM}} \\ \rho w_{\text{SBM}}\,\mathbf{u} + p\,\mathbf{n}\\ \rho w_{\text{SBM}}\,H \end{array}\right), \qquad 
   \mathbf{F}_{\mathbf{t}_k} = \left(\begin{array}{c} \rho\,u_{\mathbf{t}_k} \\ \rho\,u_{\mathbf{t}_k}\,\mathbf{u} + p\,\mathbf{t}_k\\\rho\,u_{\mathbf{t}_k}\,H \end{array}\right),
\end{equation}}
such that,
\begin{align}\label{eq:SBM un}
w_{\text{SBM}} = w% [u_n(\mathbf{\tilde x+d})]_N
 &-\|\mathbf{d}\|   \sum_{i=1}^d n_i \sum_{j=1}^d n_j \; \partial_{x_j} u_i (\mathbf{\tilde x})\quad &\text{ second order }\\
 &-\|\mathbf{d}\|^2 \sum_{i=1}^d n_i \sum_{j=1}^d \sum_{k=1}^d \frac{1}{2!}\; n_j n_k \; \partial^2_{x_j x_k} u_i (\mathbf{\tilde x})\quad &\text{ third order } \nonumber \\
 &-\|\mathbf{d}\|^3 \sum_{i=1}^d n_i  \sum_{j=1}^d \sum_{k=1}^d \sum_{\ell=1}^d \frac{1}{3!}\; n_j n_k n_\ell \;\partial^3_{x_j x_k x_\ell} u_i (\mathbf{\tilde x})\quad &\text{ fourth order }\nonumber\\ 
  &\dots \ ,  &\nonumber
% &-\|\mathbf{d}\|^4 \sum_{i=1}^d n_i  \sum_{j=1}^d \sum_{k=1}^d \sum_{\ell=1}^d \sum_{m=1}^d \frac{1}{4!}\; n_j n_k n_\ell n_m\; \partial^4_{x_j x_k x_\ell x_m} u_i, \quad &\text{ fifth order } 
\end{align}
{having set}
%where, for compactness, we define
 $\mathbf{n}=\{n_j\}_{j=1,\dots,d}$ and $\mathbf{u}=\{u_j\}_{j=1,\dots,d}$.\\

{As done before, we use the fact that the Taylor series development on the right hand side of~\eqref{eq:SBM un} is exact when applied
to a polynomial of degree lower or equal to the employed expansion, and thus we simplify it and we recast the correction as}   
%, and 
%Something similar to Eq.~\eqref{eq:Dirichlet trick} can be recovered as well for the wall boundary condition, given that $[u_n(\mathbf{\tilde x+d})]_N=0$:
%
\begin{equation}\label{eq:wall trick}
{
w_{\text{SBM}} =  
\mathbf{u}(\mathbf{\tilde x}) \cdot\mathbf{n} + \left[w -
 \mathbf{u}(\mathbf{\tilde x+d})  \cdot\mathbf{n}\right].
}
%(u_n)_{SBM} =  [\rho u_n(\mathbf{\tilde x+d})]_N - \left[ \rho\mathbf{u}(\mathbf{\tilde x+d}) - \rho\mathbf{u}(\mathbf{\tilde x}) \right]\cdot\mathbf{n} = 
%- \left[ \rho\mathbf{u}(\mathbf{\tilde x+d}) - \rho\mathbf{u}(\mathbf{\tilde x}) \right]\cdot\mathbf{n}. 
\end{equation}
It should be noticed that, without the formulation in Eq.~\eqref{eq:wall trick}, in order to perform the extrapolation given in Eq.~\eqref{eq:SBM un}, high order derivatives of the velocity components would have been needed.
However, since all derivatives are usually computed with respect to the conserved variables, either the \textit{chain rule} or some kind of linearization
would have been to be implemented to recover the necessary higher order derivatives. %This implementation issue is well deepen in Sect.~\ref{appendix1}. 
\Rall{
\begin{remark}[Boundary flux and penalty term]
{When using the numerical flux~\eqref{eq:f_wall} instead of a classical numerical flux $\hat F(\mathbf{U}^-_h,  \mathbf{U}^{bc}_{\tilde\Gamma})$
(as the Rusanov flux), for high (third, fourth, etc) order schemes, in order to obtain the correct convergence rates we had to 
include a penalty term similar to the diffusion term of the Rusanov flux. For slip walls this term reads
\begin{equation}
\mathcal{P}_w := \alpha_w \left( \mathbf{U} - \mathbf{U}^{bc}_{\tilde\Gamma} \right) = \alpha_w\,\rho \left(\begin{array}{c} 0 \\ \mathbf{u}\cdot\mathbf{n} - w_{\text{SBM}} \\[5pt] \dfrac{(\mathbf{u}\cdot\mathbf{n})^2}{2} - \dfrac{w_{\text{SBM}}^2}{2} \end{array}\right),
\end{equation}
where $\alpha_w = \|\mathbf{u}\| + \sqrt{\gamma p/\rho}$.}
\end{remark}
}

Finally, {for non moving walls with $w=0$ we consider the simpler strategy consisting in imposing directly $W:=(\rho  \mathbf{u}   \cdot\mathbf{n})|_{\Gamma_w}=0$. This allows  to  work with derivatives and variations of the momentum variable. 
We thus  replace~\eqref{eq:SBM un} by
\begin{align}\label{eq:SBM_run}
W_{SBM} =  
 &-\|\mathbf{d}\|   \sum_{i=1}^d n_i \sum_{j=1}^d n_j \; \partial_{x_j} (\rho u)_i (\mathbf{\tilde x})\quad &\text{ second order } \nonumber\\
 &-\|\mathbf{d}\|^2 \sum_{i=1}^d n_i \sum_{j=1}^d \sum_{k=1}^d \frac{1}{2!}\; n_j n_k \; \partial^2_{x_j x_k} (\rho u)_i (\mathbf{\tilde x})\quad &\text{ third order } \\
 &-\|\mathbf{d}\|^3 \sum_{i=1}^d n_i  \sum_{j=1}^d \sum_{k=1}^d \sum_{\ell=1}^d \frac{1}{3!}\; n_j n_k n_\ell \;\partial^3_{x_j x_k x_\ell} (\rho u)_i (\mathbf{\tilde x})\quad &\text{ fourth order } \nonumber\\ 
 &\dots \ ,  &\nonumber
% &-\|\mathbf{d}\|^4 \sum_{i=1}^d n_i  \sum_{j=1}^d \sum_{k=1}^d \sum_{\ell=1}^d \sum_{m=1}^d \frac{1}{4!}\; n_j n_k n_\ell n_m\; \partial^4_{x_j x_k x_\ell x_m} u_i, \quad &\text{ fifth order } 
\end{align}
and~\eqref{eq:wall trick} by}
\begin{equation}\label{eq:wall trick2}
{
W_{\text{SBM}} =   (\rho \mathbf{u})(\mathbf{\tilde x}) \cdot\mathbf{n}  - (\rho \mathbf{u})(\mathbf{\tilde x+d})  \cdot\mathbf{n}.
}
%w = - \left[ \rho\mathbf{u}(\mathbf{\tilde x+d}) - \rho\mathbf{u}(\mathbf{\tilde x}) \right]\cdot\mathbf{n}. 
\end{equation}

{As discussed  before a} fully consistent  definition of the flux {is} %can be 
obtained by {consistently correcting the value of the total enthalpy as in~\eqref{eq:sbm enthalpy0}, by replacing $w$ with $w_{SBM}$ or with $W_{SBM}/\rho$ depending on the case.}

%% noting that in all
%%expressions involving the velocity we should substitute 
%%%
%%\begin{equation}\label{eq:sbm speed}
%%\mathbf{\bar u}\,=\,w\,\mathbf{n}\,+\,u_t\,\mathbf{t}
%%\end{equation}
%%%
%%which involve a modification also in the computation of the total enthalpy which is going to be
%%
%\begin{equation}\label{eq:sbm enthalpy}
%\bar H = \frac{\gamma}{\gamma-1}\frac{p}{\rho}+\frac{1}{2}\left(w^2 + u_t^2\right)
%\end{equation}
%
%Although, when working with the polynomial correction we avoided these modifications in the speed and enthalpy in order
%not to introduce any linearization process. 

\subsection{{Other existing approaches}}

{The corrections proposed in the previous paragraph will be compared to the approach proposed by Krivodonova and Berger 
in~\cite{krivodonova2006high}, referred to as \textit{algorithm I} in the cited reference. 
\RIIcolor{The aforementioned approach, along with two more \textit{algorithms}, was introduced by the Krivodonova and Berger to provide improved solutions for two-dimensional slip-wall boundary conditions.
Even though, their corrections are limited to two-dimensional domains and slip-wall boundary conditions, 
we recall for completeness the one that we used herein to make comparisons.
More details can also be found in~\cite{krivodonova2006high}}.\\
 
We start by defining in each quadrature point a special state of primitive variables %is defined 
$\mathbf{U}^b$ and a corresponding numerical flux:}
\begin{equation}\label{eq:BergerFlux}
\mathbf{U}^b = \left(\begin{array}{c} \rho^b \\ \mathbf{u}^b \\ p^b  \end{array}\right) = \left(\begin{array}{c} \rho \\ u_t\,\mathbf{t} \\ p  \end{array}\right),\qquad 
\mathbf{F}^{bc}_{\tilde{\mathbf{n}}}(\mathbf{U}^b) = \left(\begin{array}{c} \rho^b\,(\mathbf{u}^b\cdot\tilde{\mathbf{n}})\\ \rho^b\,(\mathbf{u}^b\cdot\tilde{\mathbf{n}})\,\mathbf{u}^b+ p\,\tilde{\mathbf{n}} \\ \rho^b\,(\mathbf{u}^b\cdot\tilde{\mathbf{n}})\,H^b  \end{array}\right).
\end{equation}
%
%As mentioned above, 
%Flux~\eqref{eq:BergerFlux} can be easily recovered from Eqs.~\eqref{eq:BC integral} and~\eqref{eq:SBM fluxes}, by imposing $w=0$ in Eq.~\eqref{eq:sbm speed}:
%%
%\begin{align}
%\mathbf{F}^{bc}_{\tilde{\mathbf{n}}} &= \left(\begin{array}{c} \rho\,w \\ \rho\,w\,\mathbf{\bar u} + p\,\mathbf{n}\\\rho\,w\,\bar H \end{array}\right)(\mathbf{\tilde n \cdot n}) + \left(\begin{array}{c} \rho\,u_t \\ \rho\,u_t\,\mathbf{\bar u} + p\,\mathbf{t}\\\rho\,u_t\,\bar H \end{array}\right)(\mathbf{\tilde n \cdot t}) \nonumber \\
%&= \left(\begin{array}{c} \rho\,w \\ \rho\,w\,(w\,\mathbf{n}\,+\,u_t\,\mathbf{t}) + p\,\mathbf{n}\\\rho\,w\,\bar H \end{array}\right)(\mathbf{\tilde n \cdot n}) + \left(\begin{array}{c} \rho\,u_t \\ \rho\,u_t\,(w\,\mathbf{n}\,+\,u_t\,\mathbf{t}) + p\,\mathbf{t}\\\rho\,u_t\,\bar H \end{array}\right)(\mathbf{\tilde n \cdot t}) \nonumber \\
%&= \left(\begin{array}{c} 0 \\ p\,\mathbf{n}\\ 0 \end{array}\right)(\mathbf{\tilde n \cdot n}) + \left(\begin{array}{c} \rho\,u_t \\ \rho\,u_t\,(u_t\,\mathbf{t}) + p\,\mathbf{t}\\\rho\,u_t\,\bar H \end{array}\right)(\mathbf{\tilde n \cdot t}) \nonumber \\
%&= \left(\begin{array}{c} \rho\,u_t(\mathbf{t\cdot\tilde n}) \\ \rho\,u_t\,(\mathbf{t\cdot\tilde n})\,u_t\,\mathbf{t} + p\,\tilde{\mathbf{n}} \\ \rho\,u_t\,(\mathbf{t\cdot\tilde n})\,\bar H  \end{array}\right)
%\end{align}
%
{Note that~\cite{krivodonova2006high} contains a typo in the flux expression which does not include the pressure term.}
Even though the modification introduced in~\eqref{eq:BergerFlux} already allows a fair improvement of the
discretization error and convergence rate, {we will see that} this {shows some limitations when increasing the accuracy already beyond third order~\cite{costa2018very}.} 
\begin{remark}[Shock limiting]
\RIIcolor{For simulations with shocks,} when the \textit{a posteriori} limiter marks a cell on the wall boundary as troubled, the boundary conditions must be applied differently
{since the high order modes are discarded and the solution is updated with a} second order finite volume scheme. {In this case, we adopt  the flux modification~\eqref{eq:BergerFlux}.} 
\end{remark}

\bigskip
%--------- SECTION ---------------------------------------------------------------------------
\section{Numerical results} \label{sec.test}
%---------------------------------------------------------------------------------------------

In this Section, we test the new modified boundary treatments with several academic test-cases proving that the new method
is able to provide high-order convergence for both far-field and wall boundary conditions on 2D and 3D unstructured meshes. 
We also show the numerical results obtained on a problem that involves shocks, correctly captured in the framework of ADER-DG methods thanks to our \textit{a posteriori} subcell FV limiting technique. 
The results are provided with convergence analysis performed with classical and modified boundary conditions, and, when possible, with curvilinear meshes.

\subsection{2D tests with smooth solutions}

We start our suite of benchmarks with 2-dimensional tests involving smooth solution profiles; this easily allows to asses the claimed properties of the proposed SBM formulation.

\subsubsection{Manufactured solution on 2D curved domains: far-field BC} \label{ssec.Dir-2D}
\begin{figure}[b!]
	\centering
	%\subfigure[Coarse mesh]
	{\includegraphics[width=0.45\textwidth]{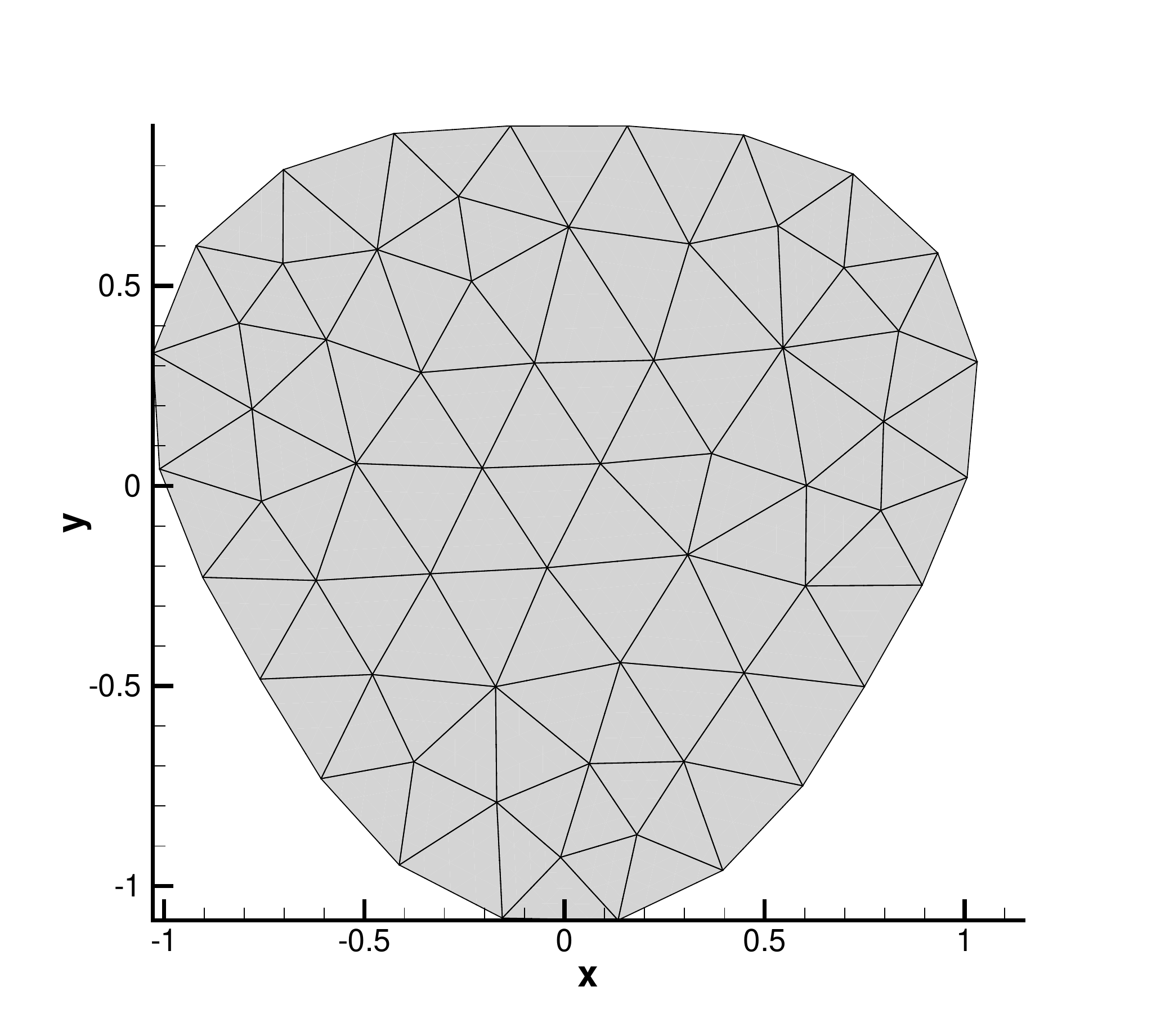}}
	%\subfigure[Initial condition]
	{\includegraphics[width=0.45\textwidth]{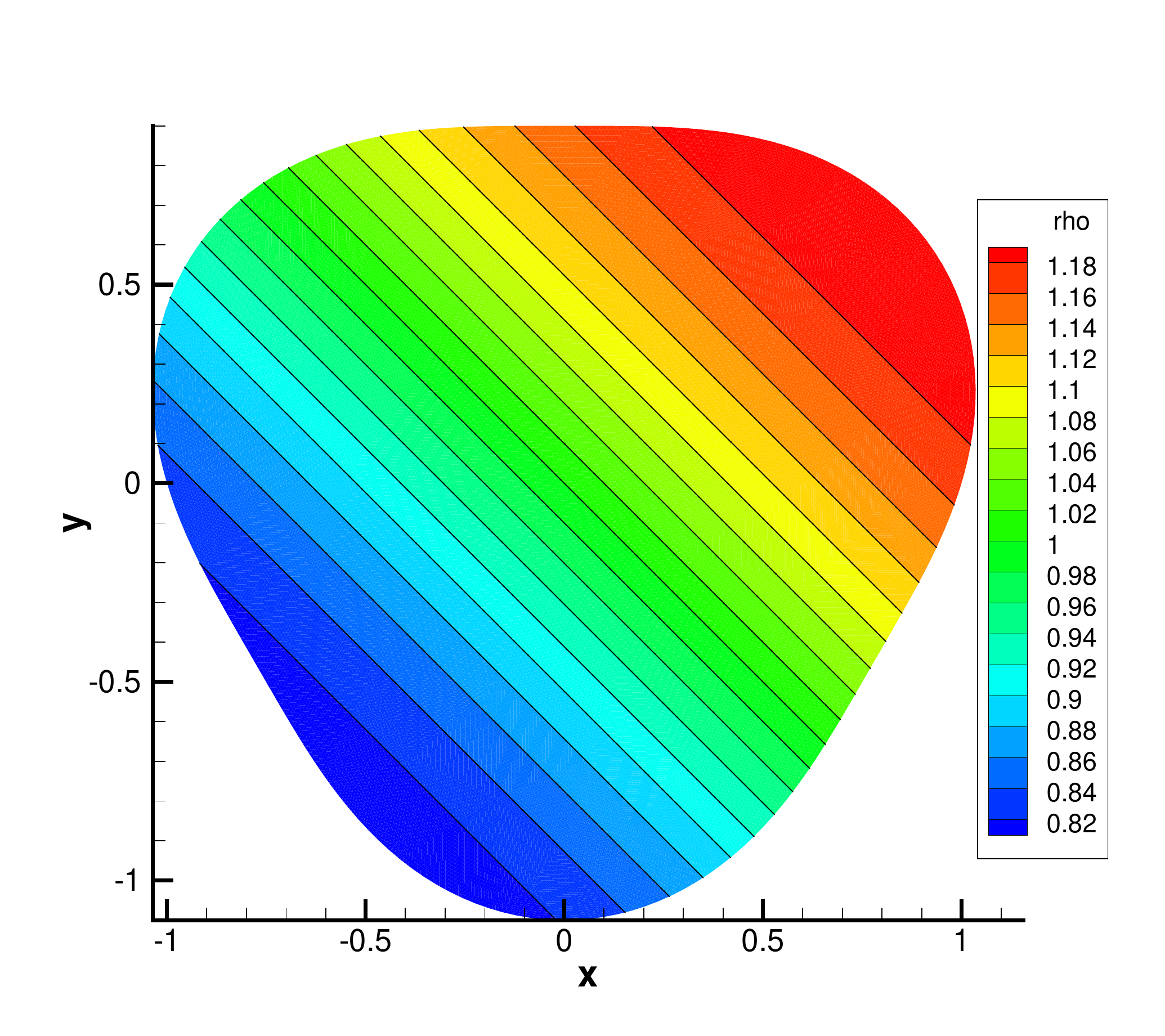}}
	\caption{Test case setup for the 2D manufactured solution test case presented in Section~\ref{ssec.Dir-2D}. We plot the coarsest employed mesh on the left and the initial density profile on the right.}
	\label{fig:Manuftest2}
\end{figure}

In order to assess the capability of our new flux correction, we start from a manufactured solution by considering the two-dimensional 
inhomogeneous Euler equations:
\begin{equation}
\mathbf{U}_t + \nabla\cdot \mathbf{F(U)} = \mathbf{S},\qquad \text{with}\qquad \mathbf{S} = \left(\begin{array}{c} 0.4\,\cos(x+y)\\ 0.6\,\cos(x+y)\\ 0.6\,\cos(x+y)\\ 1.8\,\cos(x+y) \end{array}\right).
\end{equation}
This system has the following exact steady state solution, as given in~\cite{mazaheri2019bounded},
\begin{equation}
\rho = 1+0.2\,\sin(x+y), \quad u = 1, \quad v = 1, \quad p = 1+0.2\,\sin(x+y),
\end{equation}
which is imposed on the real curved domain as far-field boundary conditions.
A very coarse mesh was generated and then refined by splitting. The four nested grids described in Table~\ref{tb:mesh manufactured2}
have been used to perform grid-convergence analysis.

The boundary where the far-field condition is applied will be referred to as $\Gamma_D$.
We tested the new SBM flux correction on a complicated geometry, taken from~\cite{costa2018very}, that can be described with the following equation written in polar coordinates:
\begin{equation}
\Gamma_D :\; \left(\begin{array}{c} x \\ y \end{array}\right) = r(\alpha,\theta) \left(\begin{array}{c} \cos\theta \\ \sin\theta \end{array}\right),\quad \text{ where }\quad r(\alpha,\theta)=r_0\left(1+\frac{1}{10}\sin(\alpha\theta)\right),\quad r_0,\alpha\in\mathbb{R},
\end{equation}
where $r_0=1$ and $\alpha=3$. We refer to \Rall{Figure~\ref{fig:Manuftest2}} for a visual representation.

Since the standard far-field boundary condition is enforced onto a curved boundary, discretized with a polygonal mesh, we expect to have second order
of accuracy at best no matter what the degree of the polynomial is and this is well-observed in the top part of Table~\ref{tb:manuf2D}.
This problem can be cured with the, however expensive, use of iso-parametric elements and curved meshes.
For the sake of completeness we provide throughout this paper also some results on curved meshes, see for example the central part of Table~\ref{tb:manuf2D}, but we underline that our technology for curved meshes allows to obtain third order accurate results with DG-$\mathbb{P}2$/$\mathbb{Q}2$ and slightly better results when using DG-$\mathbb{P}3$/$\mathbb{Q}3$, 
but without achieving the expected fourth order 
(however the improvement of curved meshes techniques is not the scope of this work).

Finally, the best results, in terms of error magnitude, convergence rates and reduced computational complexity, 
are obtained with the shifted boundary correction,
as one can see by comparing the previous results with those reported in bottom part of Table~\ref{tb:manuf2D}.
%A much better improvement is pointed out in Table~\ref{tb:manuf2D SBM}, where the results obtained with the modified Dirichlet boundary condition
%are shown. 
Convergence plots\footnote{\RIIcolor{For all convergence plots, we draw the dashed lines close to each result ($\mathbb{P}$1, $\mathbb{P}$2, $\mathbb{P}$3) to represent the optimal rate (2, 3, 4) for comparisons.}} 
until order four, for the conserved variable $\rho$, are presented in Figure~\ref{fig:dirichletConv}a.
\RIcolor{Although the paper is oriented towards polynomials of degree four, for this test case, we also presented the results obtained using polynomials of degree five to show that, 
the accuracy of the correction theoretically depends \textit{only} on the order of the polynomial itself.}
With the SBM correction, all convergence trends are correctly recovered until order five.

\begin{table}[p]
\caption{Characteristics of the employed meshes for the test case of Section~\ref{ssec.Dir-2D}, the 2D manufactured solution.}%\\[1pt]}
\label{tb:mesh manufactured2}
\scriptsize
\centering
\begin{tabular}{cccc} 
\hline %\hline
Grid level &Nodes  &Triangles  &$h$ \\[0.5mm]
\hline
0 & 61     &  98    & 1.5450E-1   \\
1 & 219    &  392   & 7.7252E-2   \\
2 & 829    & 1,568  & 3.8626E-2   \\
3 & 3,225  & 6,272  & 1.9313E-2   \\
\hline %\hline
\end{tabular}
\end{table}
\begin{table}[p]
	\caption{Convergence analysis for the test case of Section~\ref{ssec.Dir-2D}, the 2D manufactured solution. We provide the results obtained without the SBM correction on linear meshes (top part), with curved mesh (central part), and with SBM correction on linear meshes (bottom part). One can notice that the SBM correction allows to retrieve the expected high order convergence rate when working with methods of order greater than 2 on linear meshes.}%\\}
	\label{tb:manuf2D}
	\scriptsize
	\centering
	\begin{tabular}{ccccccccc} 
		\hline\\[1pt]
		\multicolumn{9}{c}{Convergence analysis \textit{without SBM} correction on \textit{linear} meshes}\\		
		\hline
		&\multicolumn{2}{c}{$\rho$} &\multicolumn{2}{c}{$\rho u$} &\multicolumn{2}{c}{$\rho v$} &\multicolumn{2}{c}{$\rho E$}\\[0.5mm]
		\cline{2-9}
		Grid level & $L_2$        & $\tilde{n}$ & $L_2$        & $\tilde{n}$ & $L_2$        & $\tilde{n}$ & $L_2$      & $\tilde{n}$ \\[0.5mm]\hline
		&\multicolumn{8}{c}{DG-$\mathbb{P}1$}\\[0.5mm]
		0          &  1.7293E-3   &   --        &  2.3514E-3   &  --         &  2.2344E-3   &   --        &  6.6466E-3 &   --    \\
		1          &  4.4790E-4   &    1.95     &  6.0925E-4   &    1.95     &  5.7869E-4   &    1.95     &  1.7138E-3 &  1.96   \\
		2          &  1.1276E-4   &    1.99     &  1.5358E-4   &    1.99     &  1.4579E-4   &    1.99     &  4.3195E-4 &  1.99   \\
		3          &  2.8241E-5   &    2.00     &  3.8473E-5   &    2.00     &  3.6512E-5   &    2.00     &  1.0820E-4 &  2.00   \\
		&\multicolumn{8}{c}{DG-$\mathbb{P}2$}\\[0.5mm]
		0          &  1.6583E-3   &   --        &  2.3344E-3   &  --         &  2.2040E-3   &   --        &  6.4258E-3 &   --    \\
		1          &  4.2195E-4   &   1.97      &  5.9156E-4   &    1.98     &  5.5908E-4   &    1.98     &  1.6308E-3 &  1.98   \\
		2          &  1.0522E-4   &   2.00      &  1.4777E-4   &    2.00     &  1.3975E-4   &    2.00     &  4.0791E-4 &  2.00   \\
		3          &  2.6187E-5   &   2.01      &  3.6855E-5   &    2.00     &  3.4867E-5   &    2.00     &  1.0170E-4 &  2.00   \\
		&\multicolumn{8}{c}{DG-$\mathbb{P}3$}\\[0.5mm]
		0          &  1.6817E-3   &   --        &  2.3809E-3   &  --         &  2.2562E-3   &   --        &  6.5431E-3 &   --    \\
		1          &  4.3095E-4   &   1.96      &  6.0642E-4   &   1.97      &  5.7430E-4   &   1.97      &  1.6640E-3 &  1.98   \\
		2          &  1.0768E-4   &   2.00      &  1.5092E-4   &   2.01      &  1.4280E-4   &   2.01      &  4.1580E-4 &  2.00   \\
		3          &  2.6792E-5   &   2.01      &  3.7467E-5   &   2.01      &  3.5428E-5   &   2.01      &  1.0345E-4 &  2.01   \\
                &\multicolumn{8}{c}{\RIcolor{DG-$\mathbb{P}4$}}\\[0.5mm]
                \RIcolor{0}          & \RIcolor{1.6950E-3}   &   \RIcolor{--  }      & \RIcolor{2.4059E-3}   &   \RIcolor{--  }      &  \RIcolor{2.2813E-3}   &   \RIcolor{--  }      &  \RIcolor{6.5915E-3} &  \RIcolor{ -- }    \\
                \RIcolor{1}          & \RIcolor{4.3318E-4}   &   \RIcolor{1.96}      & \RIcolor{6.1215E-4}   &   \RIcolor{1.97}      &  \RIcolor{5.8031E-4}   &   \RIcolor{1.98}      &  \RIcolor{1.6714E-3} &  \RIcolor{1.98}   \\
                \RIcolor{2}          & \RIcolor{1.0849E-4}   &   \RIcolor{2.00}      & \RIcolor{1.5293E-4}   &   \RIcolor{2.00}      &  \RIcolor{1.4497E-4}   &   \RIcolor{2.00}      &  \RIcolor{4.1861E-4} &  \RIcolor{2.00}   \\
                \RIcolor{3}          & \RIcolor{2.6976E-5}   &   \RIcolor{2.00}      & \RIcolor{3.7944E-5}   &   \RIcolor{2.01}      &  \RIcolor{3.5943E-5}   &   \RIcolor{2.01}      &  \RIcolor{1.0416E-4} &  \RIcolor{2.00}   \\
		\hline\\[1pt]
		\multicolumn{9}{c}{Convergence analysis on \textit{curved} meshes}\\
		\hline
		&\multicolumn{2}{c}{$\rho$} &\multicolumn{2}{c}{$\rho u$} &\multicolumn{2}{c}{$\rho v$} &\multicolumn{2}{c}{$\rho E$}\\[0.5mm]
		\cline{2-9}
		Grid level & $L_2$        & $\tilde{n}$ & $L_2$        & $\tilde{n}$ & $L_2$        & $\tilde{n}$ & $L_2$      & $\tilde{n}$ \\[0.5mm]\hline
		&\multicolumn{8}{c}{DG-$\mathbb{P}1$/$\mathbb{Q}1$}\\[0.5mm]
		0          &  1.7293E-3   &   --        &  2.3514E-3   &  --         &  2.2344E-3   &   --        &  6.6466E-3 &   --    \\
		1          &  4.4790E-4   &    1.95     &  6.0925E-4   &    1.95     &  5.7869E-4   &    1.95     &  1.7138E-3 &  1.96   \\
		2          &  1.1276E-4   &    1.99     &  1.5358E-4   &    1.99     &  1.4579E-4   &    1.99     &  4.3195E-4 &  1.99   \\
		3          &  2.8241E-5   &    2.00     &  3.8473E-5   &    2.00     &  3.6512E-5   &    2.00     &  1.0820E-4 &  2.00   \\
		&\multicolumn{8}{c}{DG-$\mathbb{P}2$/$\mathbb{Q}2$}\\[0.5mm]
		0          & 7.8113E-5    &   --        &  4.9287E-5   &  --         &  6.0091E-5   &   --        &  2.3914E-4 &   --    \\
		1          & 1.1494E-5    &    2.76     &  6.1144E-6   &    3.01     &  7.9719E-6   &    2.91     &  3.4731E-5 &  2.78   \\
		2          & 1.5902E-6    &    2.85     &  7.8242E-7   &    2.97     &  1.0271E-6   &    2.96     &  4.7828E-6 &  2.86   \\
		3          & 2.1096E-7    &    2.91     &  1.0184E-7   &    2.94     &  1.3078E-7   &    2.97     &  6.3447E-7 &  2.91   \\
		&\multicolumn{8}{c}{DG-$\mathbb{P}3$/$\mathbb{Q}3$}\\[0.5mm]
		0          & 5.2248E-6    &   --        &  4.7012E-6   &  --         &  4.6229E-6   &   --        &  1.5742E-5 &   --    \\
		1          & 5.6390E-7    &    3.21     &  4.7322E-7   &    3.31     &  4.7014E-7   &    3.29     &  1.6226E-6 &  3.28   \\
		2          & 6.1309E-8    &    3.20     &  4.7773E-8   &    3.31     &  4.8063E-8   &    3.29     &  1.7095E-7 &  3.25   \\
		3          & 1.0299E-8    &    2.57     &  9.2506E-9   &    2.36     &  1.0123E-8   &    2.25     &  3.6478E-8 &  2.23   \\
		\hline\\[1pt]
		\multicolumn{9}{c}{Convergence analysis \textit{with SBM} correction on \textit{linear} meshes}\\
		\hline
		&\multicolumn{2}{c}{$\rho$} &\multicolumn{2}{c}{$\rho u$} &\multicolumn{2}{c}{$\rho v$} &\multicolumn{2}{c}{$\rho E$}\\[0.5mm]
		\cline{2-9}
		Grid level & $L_2$        & $\tilde{n}$ & $L_2$        & $\tilde{n}$ & $L_2$        & $\tilde{n}$ & $L_2$      & $\tilde{n}$ \\[0.5mm]\hline
		&\multicolumn{8}{c}{DG-$\mathbb{P}1$/SBM-$\mathbb{P}1$}\\[0.5mm]
		0          &  1.0784E-3   &   --        &  1.2139E-3   &  --         &  1.2479E-3   &   --        &  3.3758E-3 &   --    \\
		1          &  2.3747E-4   &   2.18      &  2.4029E-4   &   2.34      &  2.4644E-4   &    2.34     &  7.2503E-4 &  2.22   \\
		2          &  5.5268E-5   &   2.10      &  5.3230E-5   &   2.17      &  5.4522E-5   &    2.18     &  1.6929E-4 &  2.10   \\
		3          &  1.3351E-5   &   2.05      &  1.2599E-5   &   2.08      &  1.2898E-5   &    2.08     &  4.1189E-5 &  2.04   \\
		&\multicolumn{8}{c}{DG-$\mathbb{P}2$/SBM-$\mathbb{P}2$}\\[0.5mm]
		0          &  7.5258E-5   &   --        &  5.0257E-5   &  --         &  5.8531E-5   &   --        & 2.2711E-4  &   --    \\
		1          &  1.1105E-5   &   2.76      &  5.6058E-6   &   3.16      &  7.1750E-6   &    3.03     & 3.2987E-5  &  2.78   \\
		2          &  1.5701E-6   &   2.82      &  7.2761E-7   &   2.95      &  9.2582E-7   &    2.95     & 4.6465E-6  &  2.83   \\
		3          &  2.1377E-7   &   2.88      &  1.0003E-7   &   2.86      &  1.1793E-7   &    2.97     & 6.3270E-7  &  2.88   \\
		&\multicolumn{8}{c}{DG-$\mathbb{P}3$/SBM-$\mathbb{P}3$}\\[0.5mm]
		0          &  1.7331E-6   &   --        &  2.8523E-6   &  --         &  2.9879E-6   &   --        & 6.3787E-6  &   --    \\
		1          &  7.0445E-8   &   4.62      &  9.8935E-8   &   4.85      &  1.0298E-7   &    4.86     & 2.3519E-7  &  4.76   \\
		2          &  3.2654E-9   &   4.43      &  3.9178E-9   &   4.66      &  4.0466E-9   &    4.67     & 1.0200E-8  &  4.53   \\
		3          &  1.7487E-10  &   4.22      &  1.8709E-10  &   4.39      &  1.9248E-10  &    4.39     & 5.3564E-10 &  4.25   \\
                &\multicolumn{8}{c}{\RIcolor{DG-$\mathbb{P}4$/SBM-$\mathbb{P}4$}}\\[0.5mm]
                \RIcolor{0}          &  \RIcolor{4.4808E-8 }  &   \RIcolor{--  }      &  \RIcolor{4.9888E-8 }  &    \RIcolor{ -- }     &  \RIcolor{5.3190E-8 }  &   \RIcolor{--  }      &  \RIcolor{1.4804E-7 }&  \RIcolor{ -- }   \\
                \RIcolor{1}          &  \RIcolor{1.2360E-9 }  &   \RIcolor{5.18}      &  \RIcolor{9.2420E-10}  &    \RIcolor{5.75}     &  \RIcolor{1.0080E-9 }  &   \RIcolor{5.72}      &  \RIcolor{3.8553E-9 }&  \RIcolor{5.26}   \\
                \RIcolor{2}          &  \RIcolor{4.0526E-11}  &   \RIcolor{4.93}      &  \RIcolor{2.5906E-11}  &    \RIcolor{5.16}     &  \RIcolor{2.7614E-11}  &   \RIcolor{5.19}      &  \RIcolor{1.2496E-10}&  \RIcolor{4.95}   \\
                \RIcolor{3}          &  \RIcolor{   --     }  &   \RIcolor{ -- }      &  \RIcolor{   --     }  &    \RIcolor{--  }     &  \RIcolor{  --      }  &   \RIcolor{   -}-     &  \RIcolor{ --       }&  \RIcolor{--  }   \\
		\hline
	\end{tabular}
\end{table}

\subsubsection{Supersonic vortex bounded by two circular walls: slip wall BC} \label{ssec.Wall-2D}
\begin{figure}%[b!]
	\centering
	{\includegraphics[width=0.45\textwidth]{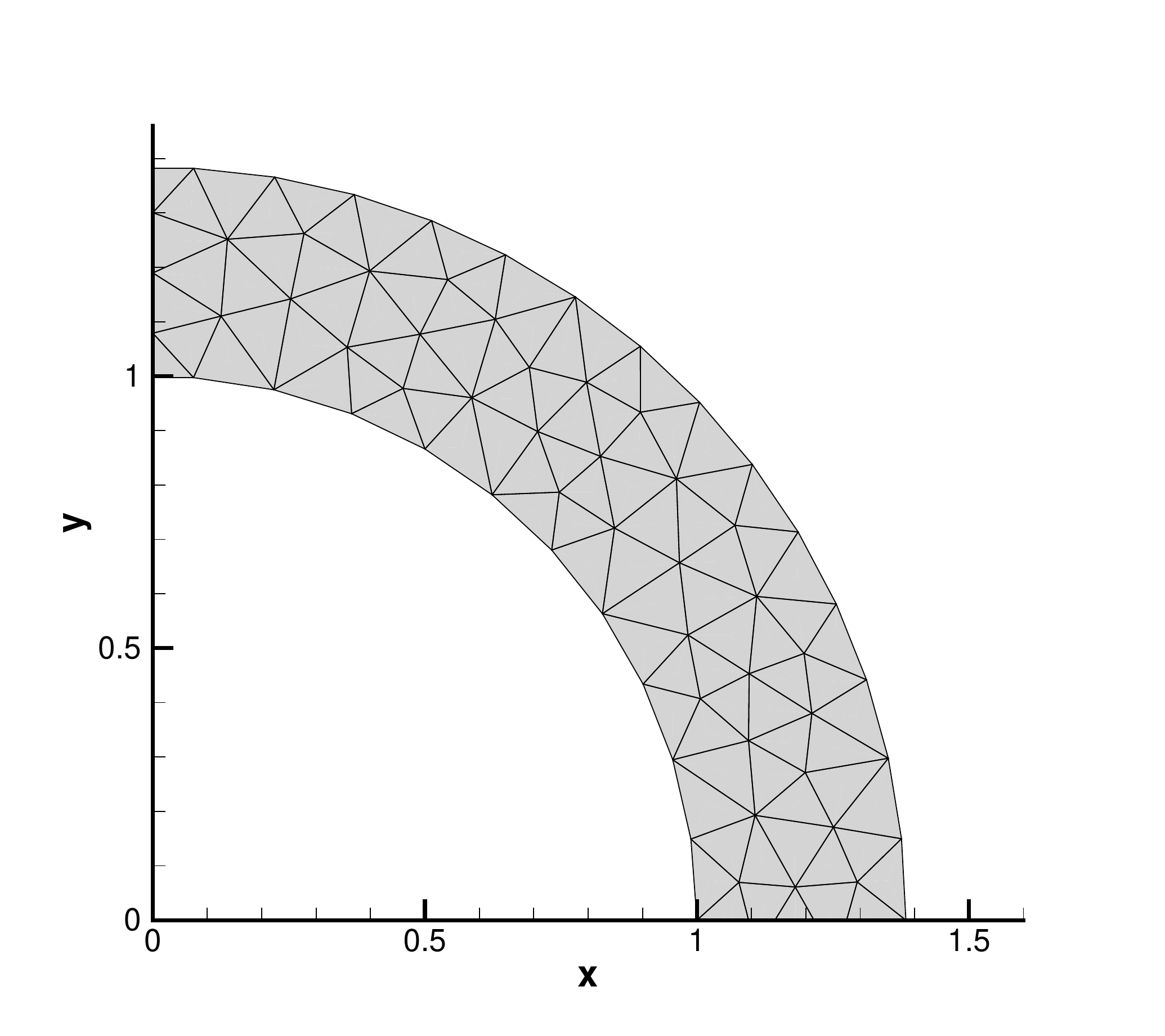}}
	{\includegraphics[width=0.45\textwidth]{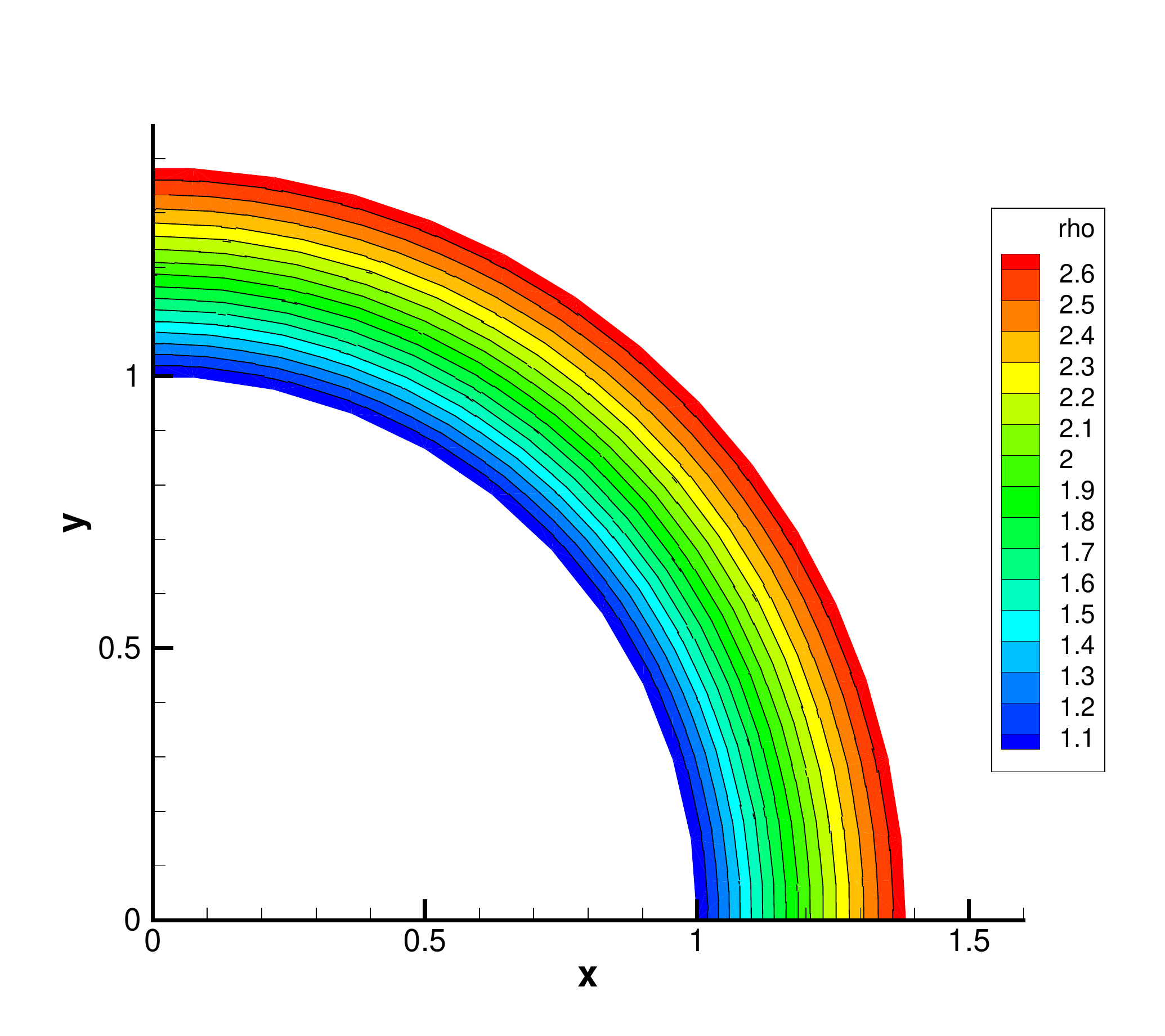}}
	\caption{Test case setup for the supersonic vortex bounded by two circular walls presented in Section~\ref{ssec.Wall-2D}. We plot the coarsest employed mesh on the left and the initial density profile on the right.}
	\label{fig:bergertest}
\end{figure}

In order to test the new flux corrections for wall boundary conditions, we now consider an isentropic supersonic flow between two concentric circular arcs of radii $r_i = 1$ and
$r_o = 1.384$. The exact density in terms of radius $r$ is given by
\begin{equation}\label{eq:vortex2D rho}
\rho = \rho_i \biggl( 1 + \frac{\gamma+1}{2}M_i^2  \biggl(1- \biggl(\frac{r_i}{r}\biggr)^2 \biggr) \biggr)^{\frac{1}{\gamma-1}}.
\end{equation}
The velocity and pressure are given by
\begin{equation}\label{eq:vortex2D u p}
\|\mathbf{u}\|=\frac{c_i M_i}{r}, \quad p = \frac{\rho^\gamma}{\gamma},
\end{equation}
where $c_i$ is the speed of sound on the inner circle. The Mach number on the inner circle $M_i$ is set to 2.25 and
the density $\rho_i$ to one.
The fluid's velocity vector components in $(x,y)$ can be computed as follows:
\begin{equation}
\left( \begin{array}{c} u \\ v\end{array} \right) = \|\mathbf{u}\| \left( \begin{array}{c} y/r \\ -x/r\end{array} \right),
\end{equation}
where $r=\sqrt{x^2+y^2}$.

A set of refined meshes is obtained by means of conformal refinement of an initial triangulation
shown in Figure~\ref{fig:bergertest}; the mesh characteristics are given in Table~\ref{tab.meshcarvortex2d}.
For this test case, we first performed the convergence test by enforcing the standard weak wall boundary condition \RIIcolor{(given by Equations~\eqref{eq:slip} and~\eqref{eq:f_wall})} 
on the polygonal boundary. Table~\ref{tb:supvortex2Dcurved} points out the grid-convergence analysis showing that the method is not even able to converge with
a full second order rate. It can also be noticed that, for each mesh, the higher is the degree of the polynomial, the larger the error is.
This is due to the fact that the scheme is trying to better approximate a solution that is indeed wrong because the boundary, which should be curved, is approximated by a polygon (see Figure~\ref{fig:SBM1}).

The same test, with the same boundary conditions, has also been run with high order curvilinear meshes showing that when increasing the order of the polynomial used to approximate the boundary
we are able to recover the formal order of accuracy of the method, {as visible on Table~\ref{tb:supvortex2Dcurved}.  Note that this result is only
possible when using iso-parametric elements, so} the {degree of the} polynomial $\mathbb{Q}$ that approximates the {geometry and the} mesh has to increase with the polynomial $\mathbb{P}$ of the scheme
(e.g.\ DG-$\mathbb{P}1$/$\mathbb{Q}1$, DG-$\mathbb{P}2$/$\mathbb{Q}2$, and DG-$\mathbb{P}3$/$\mathbb{Q}3$).
Convergence plots obtained with curvilinear element are shown in Figure~\ref{fig:slipwallConv}a for the density variable $\rho$.
%It should be recalled that, in general, generating high order meshes is not trivial and increasing the number of degree of freedom introduce a huge bottleneck in the simulation time.
%The grid-convergence analysis for this case is presented in Table~\ref{tb:berger2D curved} showing the expected results.

{Finally, the results obtained using the SBM corrections on straight sided meshes are reported in Table~\ref{tb:supvortex2DSBM} and shown in Figure~\ref{fig:slipwallConv}.
The top part of the Table shows the results obtained by using the truncated Taylor series development for the velocity; instead, 
the bottom part concerns the simplified polynomials correction~\eqref{eq:wall trick2} based on the momentum variables which does not require complex derivative evaluations.}
%Then, we tested our new SBM wall flux correction on simple linear meshes and, as shown in Table~\ref{tb:berger2D SBM}, all convergence trends are well recovered.
%To be complete, we also included in Table~\ref{tb:berger2D SBMtrick} the results obtained with the polynomials correction~\eqref{eq:wall trick2}, without any need of computing 
%high order partial derivatives. 
%It should be noticed that without introducing any linearization we are able to obtain almost superimposed results with the classical approach.
%They do not match only because the speed vector has not been updated as it was done previously. 
{We observe the expected rates and low error levels in both cases, 
with a slight improvement with the simplified approach which may be related}
to the fact that we are extrapolating 
the whole $\rho u_n$ term, rather than only $u_n$, {as} %it was 
done for the classical SBM formulation (with $\rho$ taken from the quadrature point).

For the sake of completeness, we present in Table~\ref{tb:berger2Dkrivo} and Figure~\ref{fig:slipwallConv}b the results obtained with the \Rall{correction~\eqref{eq:BergerFlux}}.
It should be noticed that even for this simple case, considering quadratic geometries, we have an order-of-accuracy degradation for finer meshes
when using polynomials of degree higher than two ($\mathbb{P}3$).

%Instead, better convergence trends are given also in this situation by the SBM correction in Figure~\ref{fig:slipwallConv}c.

\begin{table}[p]
	\caption{Characteristics of the employed meshes for the test case of Section~\ref{ssec.Wall-2D}.}
	\label{tab.meshcarvortex2d}
	\scriptsize
	\centering
	\begin{tabular}{cccc} \hline %\hline
		Grid level &Nodes  &Triangles  &$h$ \\[0.5mm]
		\hline
		0 & 238    & 376    & 1.024E-01 \\
		1 & 852    & 1,504  & 5.121E-02 \\
		2 & 3,208  & 6,016  & 2.560E-02 \\
		3 & 12,432 & 24,064 & 1.280E-02 \\
		\hline %\hline
	\end{tabular}
\end{table}
\begin{table}
\caption{
Convergence analysis for the test case of Section~\ref{ssec.Wall-2D}, the supersonic vortex bounded by two circular walls in 2D, performed \textit{with} linear and curved elements without the SBM correction.}\label{tb:supvortex2Dcurved}
\scriptsize
\centering
\begin{tabular}{ccccccccc} 
\hline\\[1pt]
\multicolumn{9}{c}{Convergence analysis \textit{without SBM} correction on \textit{linear} meshes.}\\
\hline
&\multicolumn{2}{c}{$\rho$} &\multicolumn{2}{c}{$\rho u$} &\multicolumn{2}{c}{$\rho v$} &\multicolumn{2}{c}{$\rho E$}\\[0.5mm]
\cline{2-9}
Grid level & $L_2$        & $\tilde{n}$ & $L_2$        & $\tilde{n}$ & $L_2$        & $\tilde{n}$ & $L_2$         & $\tilde{n}$ \\[0.5mm]\hline
&\multicolumn{8}{c}{DG-$\mathbb{P}1$}\\[0.5mm]
0        & 4.0507E-2 &   --        & 5.9773E-2 &  --         & 6.0006E-2 &   --        & 1.7033E-1  &   --        \\
1        & 1.3141E-2 &   1.62      & 1.9227E-2 &   1.64      & 1.9330E-2 &   1.63      & 5.6565E-2  &   1.59      \\
2        & 4.4645E-3 &   1.56      & 6.5585E-3 &   1.55      & 6.5822E-3 &   1.55      & 1.9499E-2  &   1.54      \\
3        & 1.5786E-3 &   1.50      & 2.3646E-3 &   1.47      & 2.3684E-3 &   1.47      & 6.9664E-3  &   1.48      \\
&\multicolumn{8}{c}{DG-$\mathbb{P}2$}\\[0.5mm]
0        & 7.5758E-2 &   --        & 1.3956E-1 &  --         & 1.4152E-1 &   --        & 3.3516E-1  &   --        \\
1        & 2.7584E-2 &   1.46      & 5.1868E-2 &   1.43      & 5.1990E-2 &   1.44      & 1.2308E-1  &   1.45       \\
2        & 9.6935E-3 &   1.51      & 1.8547E-2 &   1.48      & 1.8548E-2 &   1.49      & 4.3649E-2  &   1.50       \\
3        & 3.4687E-3 &   1.48      & 6.5731E-3 &   1.50      & 6.5656E-3 &   1.50      & 1.5548E-2  &   1.49       \\
&\multicolumn{8}{c}{DG-$\mathbb{P}3$}\\[0.5mm]
0        & 1.9104E-1 &   --        & 2.2847E-1 &  --         & 2.2992E-1 &   --        & 7.2469E-1  &   --        \\
1        & 8.3972E-2 &   1.19      & 9.9892E-2 &   1.19      & 1.0244E-1 &   1.17      & 3.2130E-1  &   1.17       \\
2        & 3.5508E-2 &   1.24      & 4.0370E-2 &   1.31      & 4.1040E-2 &   1.32      & 1.3453E-1  &   1.26      \\
3        & 1.5314E-2 &   1.21      & 1.6278E-2 &   1.31      & 1.6507E-2 &   1.31      & 5.7225E-2  &   1.23       \\
\hline\\[1pt]
\multicolumn{9}{c}{Convergence analysis on \textit{curved} meshes}\\
\hline
&\multicolumn{2}{c}{$\rho$} &\multicolumn{2}{c}{$\rho u$} &\multicolumn{2}{c}{$\rho v$} &\multicolumn{2}{c}{$\rho E$}\\[0.5mm]
\cline{2-9}
Grid level & $L_2$        & $\tilde{n}$ & $L_2$        & $\tilde{n}$ & $L_2$        & $\tilde{n}$ & $L_2$         & $\tilde{n}$ \\[0.5mm]\hline
&\multicolumn{8}{c}{DG-$\mathbb{P}1$/$\mathbb{Q}1$}\\[0.5mm]
0        & 4.0507E-2 &   --     & 5.9773E-2 &  --      & 6.0006E-2 &   --     & 1.7033E-1  &   --     \\
1        & 1.3141E-2 &   1.62   & 1.9227E-2 &   1.64   & 1.9330E-2 &   1.63   & 5.6565E-2  &   1.59   \\
2        & 4.4645E-3 &   1.56   & 6.5585E-3 &   1.55   & 6.5822E-3 &   1.55   & 1.9499E-2  &   1.54   \\
3        & 1.5786E-3 &   1.50   & 2.3646E-3 &   1.47   & 2.3684E-3 &   1.47   & 6.9664E-3  &   1.48   \\
&\multicolumn{8}{c}{DG-$\mathbb{P}2$/$\mathbb{Q}2$}\\[0.5mm]
0        & 1.2075E-3 &   --     & 1.7408E-3 &  --      & 1.7527E-3 &   --     & 4.8709E-3  &   --     \\
1        & 1.9156E-4 &   2.66   & 2.6389E-4 &   2.72   & 2.6721E-4 &   2.71   & 7.7666E-4  &   2.65   \\
2        & 2.5831E-5 &   2.89   & 3.3546E-5 &   2.98   & 3.3662E-5 &   2.99   & 1.0323E-4  &   2.91   \\
3        & 3.0727E-6 &   3.07   & 4.1031E-6 &   3.03   & 4.0812E-6 &   3.04   & 1.2467E-5  &   3.05   \\
&\multicolumn{8}{c}{DG-$\mathbb{P}3$/$\mathbb{Q}3$}\\[0.5mm]
0        & 7.0183E-5 &   --     & 8.8438E-5 &  --      & 8.9166E-5 &   --     & 2.7348E-4  &   --     \\
1        & 4.8983E-6 &  3.84    & 6.7326E-6 &   3.72   & 6.7616E-6 &    3.72  & 1.9264E-5  &   3.83   \\
2        & 4.0876E-7 &  3.58    & 5.9445E-7 &   3.50   & 5.9936E-7 &    3.50  & 1.6698E-6  &   3.50   \\
3        & 3.2822E-8 &  3.64    & 5.1193E-8 &   3.54   & 5.1593E-8 &    3.54  & 1.3795E-7  &   3.60   \\
\hline %\hline
\end{tabular}
\end{table}

\begin{table}
	\caption{Convergence analysis for the test case of Section~\ref{ssec.Wall-2D}, the supersonic vortex bounded by two circular walls in 2D, with the modified flux introduced in~\eqref{eq:BergerFlux}.}\label{tb:berger2Dkrivo}
	\scriptsize
	\centering
	\begin{tabular}{ccccccccc} \hline%\hline
		&\multicolumn{2}{c}{$\rho$} &\multicolumn{2}{c}{$\rho u$} &\multicolumn{2}{c}{$\rho v$} &\multicolumn{2}{c}{$\rho E$}\\[0.5mm]
		\cline{2-9}
		Grid level & $L_2$        & $\tilde{n}$ & $L_2$        & $\tilde{n}$ & $L_2$        & $\tilde{n}$ & $L_2$         & $\tilde{n}$ \\[0.5mm]\hline
		&\multicolumn{8}{c}{DG-$\mathbb{P}$1 with correction~\eqref{eq:BergerFlux}}\\[0.5mm]
		0              &  1.9580E-2   &   --        &  4.1681E-2   &  --         &  4.1868E-2   &   --        &  8.9637E-2    &   --        \\
		1              &  5.3355E-3   &   1.88      &  1.1620E-2   &    1.84     &  1.1676E-2   &   1.84      &  2.4676E-2    &   1.86       \\
		2              &  1.1690E-3   &   2.19      &  2.5781E-3   &    2.17     &  2.5983E-3   &   2.17      &  5.3291E-3    &   2.21       \\
		3              &  2.4711E-4   &   2.24      &  5.5186E-4   &    2.22     &  5.5738E-4   &   2.22      &  1.0988E-3    &   2.28       \\
		&\multicolumn{8}{c}{DG-$\mathbb{P}$2 with correction~\eqref{eq:BergerFlux}}\\[0.5mm]
		0              & 1.2764E-3    &   --        &  1.8543E-3   &  --         & 1.8656E-3    &   --        & 5.2077E-3     &   --        \\
		1              & 1.9749E-4    &    2.69     &  2.6604E-4   &    2.80     & 2.6928E-4    &   2.79      & 7.9744E-4     &    2.71      \\
		2              & 2.5994E-5    &    2.92     &  3.3359E-5   &    2.99     & 3.3458E-5    &   3.00      & 1.0360E-4     &    2.94      \\
		3              & 3.0959E-6    &    3.06     &  4.1008E-6   &    3.02     & 4.0776E-6    &   3.03      & 1.2531E-5     &    3.04      \\
		&\multicolumn{8}{c}{DG-$\mathbb{P}$3 with correction~\eqref{eq:BergerFlux}}\\[0.5mm]
		0              & 8.3269E-5    &   --        &  1.0748E-4   &  --         &  1.0489E-4   &   --        &  3.2138E-4    &   --        \\
		1              & 5.3483E-6    &   3.96      &  6.6779E-6   &   4.00      &  6.5981E-6   &    3.99     &  2.0358E-5    &   3.98       \\
		2              & 4.6348E-7    &   3.52      &  5.0053E-7   &   3.73      &  5.0276E-7   &    3.71     &  1.7297E-6    &   3.55       \\
		3              & 4.5383E-8    &   3.35      &  4.3567E-8   &   3.52      &  4.3975E-8   &    3.51     &  1.6679E-7    &   3.37       \\
		\hline%\hline
	\end{tabular}
\end{table}

\begin{table}
\caption{Convergence analysis for the test case of Section~\ref{ssec.Wall-2D}, the supersonic vortex bounded by two circular walls in 2D, performed \textit{with} the use of SBM corrections on linear meshes. We provide the results obtained with the entire Taylor series evaluation of Eq.~\eqref{eq:SBM un} (top part of the Table) and those obtained with the simplified formula of Eq.~\eqref{eq:wall trick2} (bottom part of the Table).}\label{tb:supvortex2DSBM}
\scriptsize
\centering
\begin{tabular}{ccccccccc} 
\hline\\[1pt]
\multicolumn{9}{c}{Convergence analysis \textit{with SBM} correction, as given in Eq.~\eqref{eq:SBM un}}\\
\hline
&\multicolumn{2}{c}{$\rho$} &\multicolumn{2}{c}{$\rho u$} &\multicolumn{2}{c}{$\rho v$} &\multicolumn{2}{c}{$\rho E$}\\[0.5mm]
\cline{2-9}
Grid level & $L_2$        & $\tilde{n}$ & $L_2$        & $\tilde{n}$ & $L_2$        & $\tilde{n}$ & $L_2$         & $\tilde{n}$ \\[0.5mm]\hline
&\multicolumn{8}{c}{DG-$\mathbb{P}$1/SBM-$\mathbb{P}$1}\\[0.5mm]
0          &  1.6446E-2   &   --        &  3.5594E-2   &  --         &  3.5415E-2   &   --        &  7.7895E-2    &   --        \\
1          &  4.5403E-3   &   1.86      &  9.9195E-3   &   1.84      &  9.9279E-3   &  1.84       &  2.1613E-2    &  1.85       \\
2          &  9.5300E-4   &   2.25      &  2.0713E-3   &   2.26      &  2.0770E-3   &  2.26       &  4.5047E-3    &  2.26       \\
3          &  1.9068E-4   &   2.32      &  4.1246E-4   &   2.33      &  4.1396E-4   &  2.33       &  8.9777E-4    &  2.33       \\
&\multicolumn{8}{c}{DG-$\mathbb{P}$2/SBM-$\mathbb{P}$2}\\[0.5mm]
0          & 1.2652E-3    &   --        &  1.8623E-3   &  --         &  1.8779E-3   &   --        & 5.2247E-3     &   --        \\
1          & 1.9622E-4    &   2.69      &  2.6306E-4   &   2.82      &  2.6695E-4   &  2.81       & 7.9217E-4     &  2.72       \\
2          & 2.5770E-5    &   2.93      &  3.3068E-5   &   2.99      &  3.3222E-5   &  3.01       & 1.0284E-4     &  2.95       \\
3          & 3.0587E-6    &   3.08      &  4.0533E-6   &   3.03      &  4.0349E-6   &  3.04       & 1.2390E-5     &  3.05       \\
&\multicolumn{8}{c}{DG-$\mathbb{P}$3/SBM-$\mathbb{P}$3}\\[0.5mm]
0          & 6.2321E-5    &   --        &  9.0618E-5   &  --         & 8.7148E-5    &   --        &  2.4363E-4    &   --        \\
1          & 2.9113E-6    &   4.42      &  4.9355E-6   &   4.20      & 4.7997E-6    &  4.18       &  1.1976E-5    &   4.35      \\
2          & 1.5293E-7    &   4.25      &  2.6814E-7   &   4.20      & 2.6463E-7    &  4.18       &  6.3695E-7    &   4.23      \\
3          & 8.5252E-9    &   4.17      &  1.4884E-8   &   4.17      & 1.4727E-8    &  4.17       &  3.5447E-8    &   4.17      \\
\hline\\[1pt]
\multicolumn{9}{c}{Convergence analysis \textit{with SBM} correction, as given in Eq.~\eqref{eq:wall trick2}}\\
\hline
&\multicolumn{2}{c}{$\rho$} &\multicolumn{2}{c}{$\rho u$} &\multicolumn{2}{c}{$\rho v$} &\multicolumn{2}{c}{$\rho E$}\\[0.5mm]
\cline{2-9}
Grid level & $L_2$        & $\tilde{n}$ & $L_2$        & $\tilde{n}$ & $L_2$        & $\tilde{n}$ & $L_2$         & $\tilde{n}$ \\[0.5mm]\hline
&\multicolumn{8}{c}{DG-$\mathbb{P}$1/SBM-$\mathbb{P}$1}\\[0.5mm]
0          & 1.6335E-2    &   --        &  3.5355E-2   &  --         &  3.5183E-2   &   --        &  7.7374E-2    &   --        \\
1          & 4.5259E-3    &   1.85      &  9.8900E-3   &    1.84     &  9.8996E-3   &  1.84       &  2.1548E-2    &   1.84      \\
2          & 9.5188E-4    &   2.25      &  2.0689E-3   &    2.26     &  2.0748E-3   &  2.26       &  4.4993E-3    &   2.26      \\
3          & 1.9058E-4    &   2.32      &  4.1225E-4   &    2.33     &  4.1376E-4   &  2.33       &  8.9728E-4    &   2.33      \\
&\multicolumn{8}{c}{DG-$\mathbb{P}$2/SBM-$\mathbb{P}$2}\\[0.5mm]
0          & 1.2604E-3    &   --        &  1.8419E-3   &  --         &  1.8568E-3   &   --        &  5.1834E-3    &   --        \\
1          & 1.9606E-4    &   2.69      &  2.6221E-4   &    2.82     &  2.6608E-4   &  2.81       &  7.9063E-4    &   2.72      \\
2          & 2.5765E-5    &   2.93      &  3.3039E-5   &    2.99     &  3.3192E-5   &  3.01       &  1.0279E-4    &   2.95      \\
3          & 3.0584E-6    &   3.08      &  4.0524E-6   &    3.03     &  4.0340E-6   &  3.04       &  1.2388E-5    &   3.05      \\
&\multicolumn{8}{c}{DG-$\mathbb{P}$3/SBM-$\mathbb{P}$3}\\[0.5mm] 
0          & 6.2289E-5    &   --        &  9.1157E-5   &  --         &  8.7853E-5   &   --        &  2.4332E-4    &   --        \\
1          & 2.8518E-6    &   4.45      &  4.9380E-6   &    4.20     &  4.7981E-6   &  4.19       &  1.1777E-5    &   4.37      \\
2          & 1.4569E-7    &   4.29      &  2.6571E-7   &    4.22     &  2.6206E-7   &  4.19       &  6.1275E-7    &   4.26      \\
3          & 7.9556E-9    &   4.20      &  1.4647E-8   &    4.18     &  1.4496E-8   &  4.18       &  3.3523E-8    &   4.19      \\
\hline %\hline
\end{tabular}
\end{table}

\RIcolor{
\subsubsection{Subsonic flow over a circular cylinder: slip wall BC} \label{ssec.SubsonicCyl-2D}

This test case has the same goal of that presented in~\cite{WangSun}: the validation of the new boundary treatment. 
The free-stream Mach number used in this simulation is 0.3. The comparison was carried by running the most accurate scheme used in this work, DG-$\mathbb{P}$3, with and without
the SBM correction to assess the influence of the present approach on this benchmark. 
We performed the computations with an unstructured Delaunay mesh made by 1,096 nodes and 2,087 triangles.
To impose consistent far-field conditions, we have a setup with an outer boundary 20 radiuses away from the cylinder.
In order to test the new boundary conditions we studied the symmetry of iso-contours for both density, in Figure~\ref{fig:SubCylrho}, and Mach number, in Figure~\ref{fig:SubCylMach}.
For both variables, a clear improvement of the solution is pointed out by the resuls performed with the SBM correction. Note that, when slip-wall conditions without corrections are considered, the 
flow is characterized by separation behind the body (see Figure~\ref{fig:SubCylMach}a). Instead, a much better prediction of the field is given by the results obtained with the polynomial correction.
Finally, in Figure~\ref{fig:SubCylEntropy}, we plot the spurious entropy production generated by the two boundary conditions, indicating a remarkable result also in this case. 

\begin{figure}
        \centering
        {\includegraphics[width=0.33\textwidth]{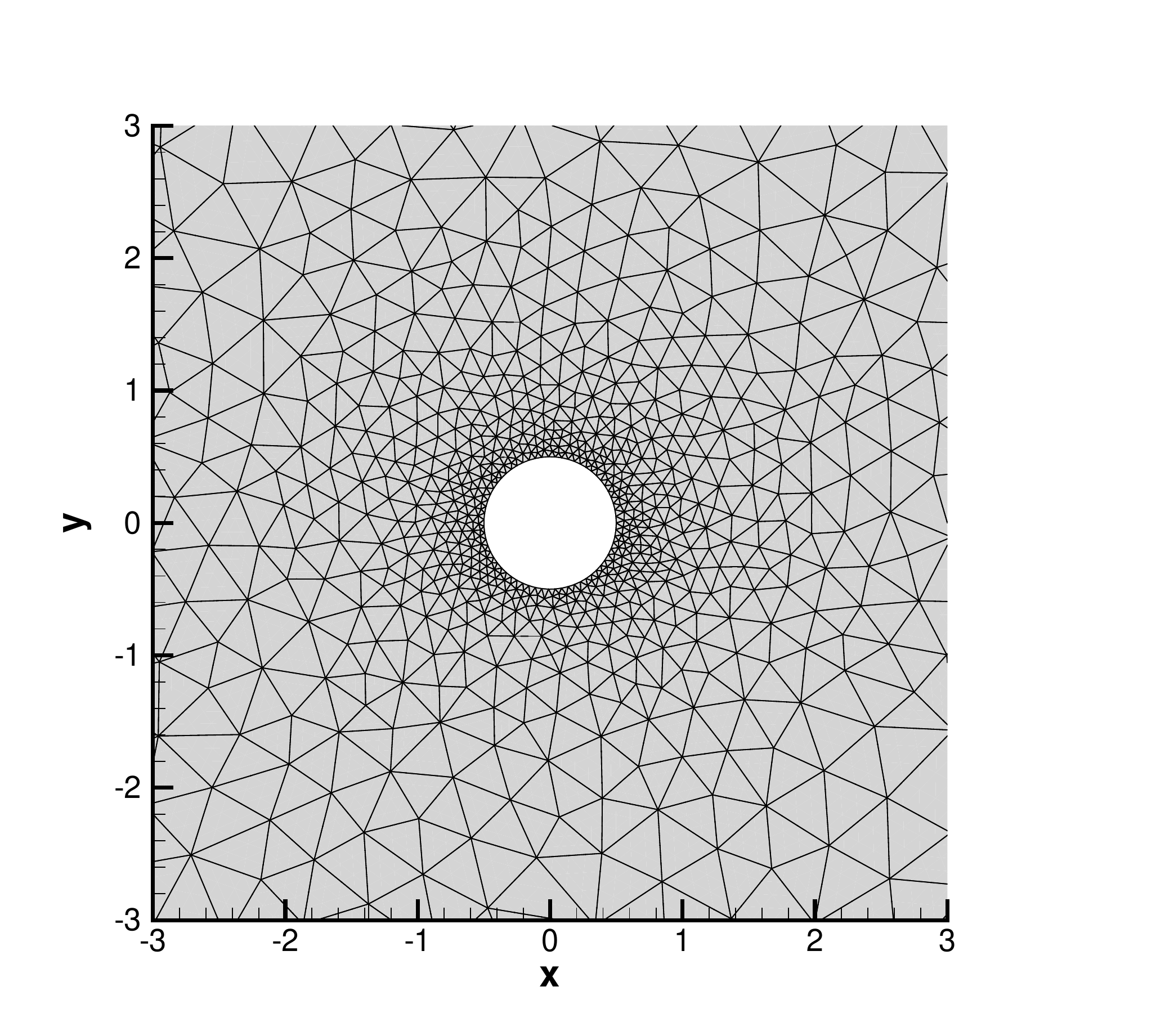}}
        {\includegraphics[width=0.33\textwidth]{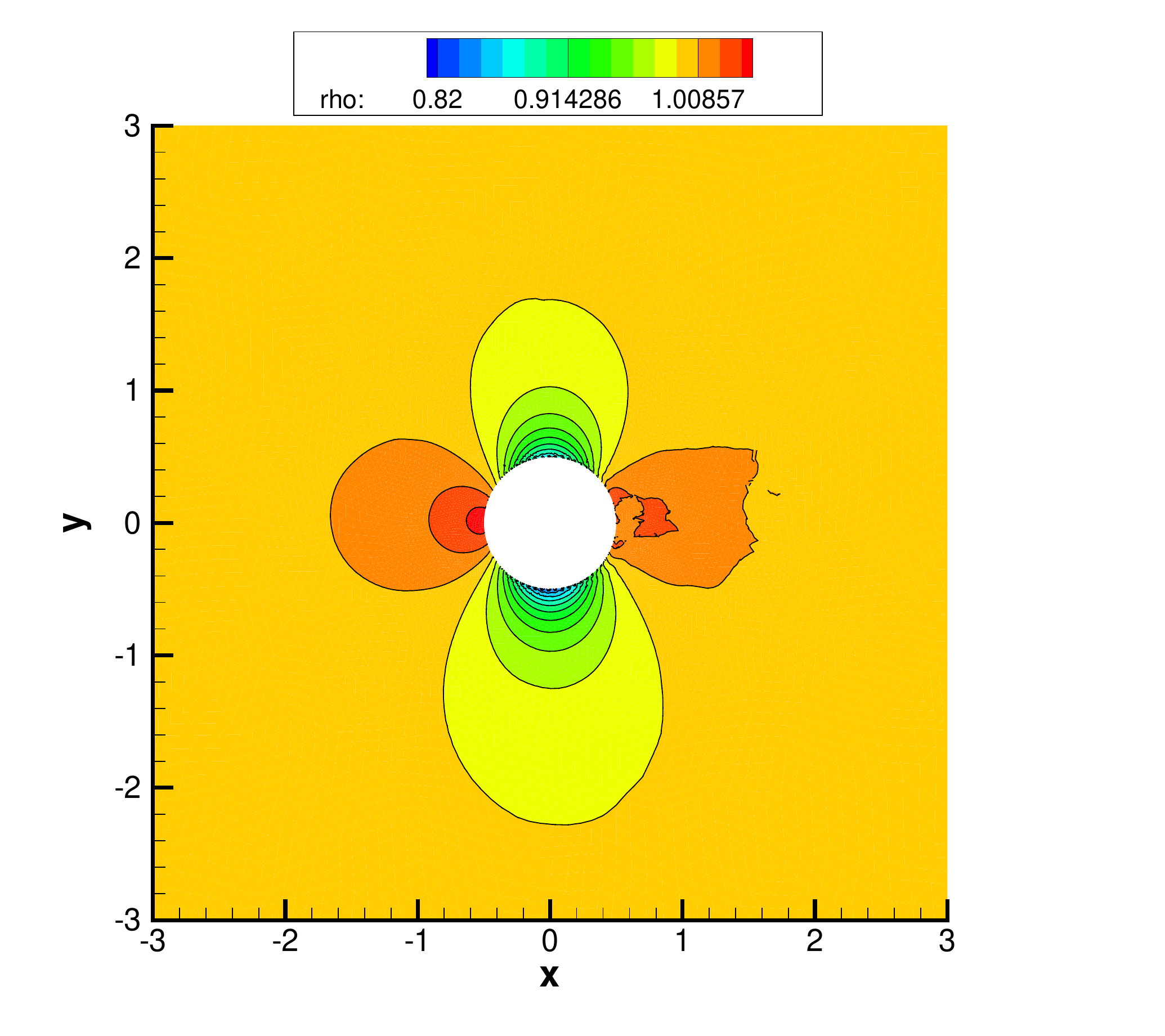}}
        {\includegraphics[width=0.33\textwidth]{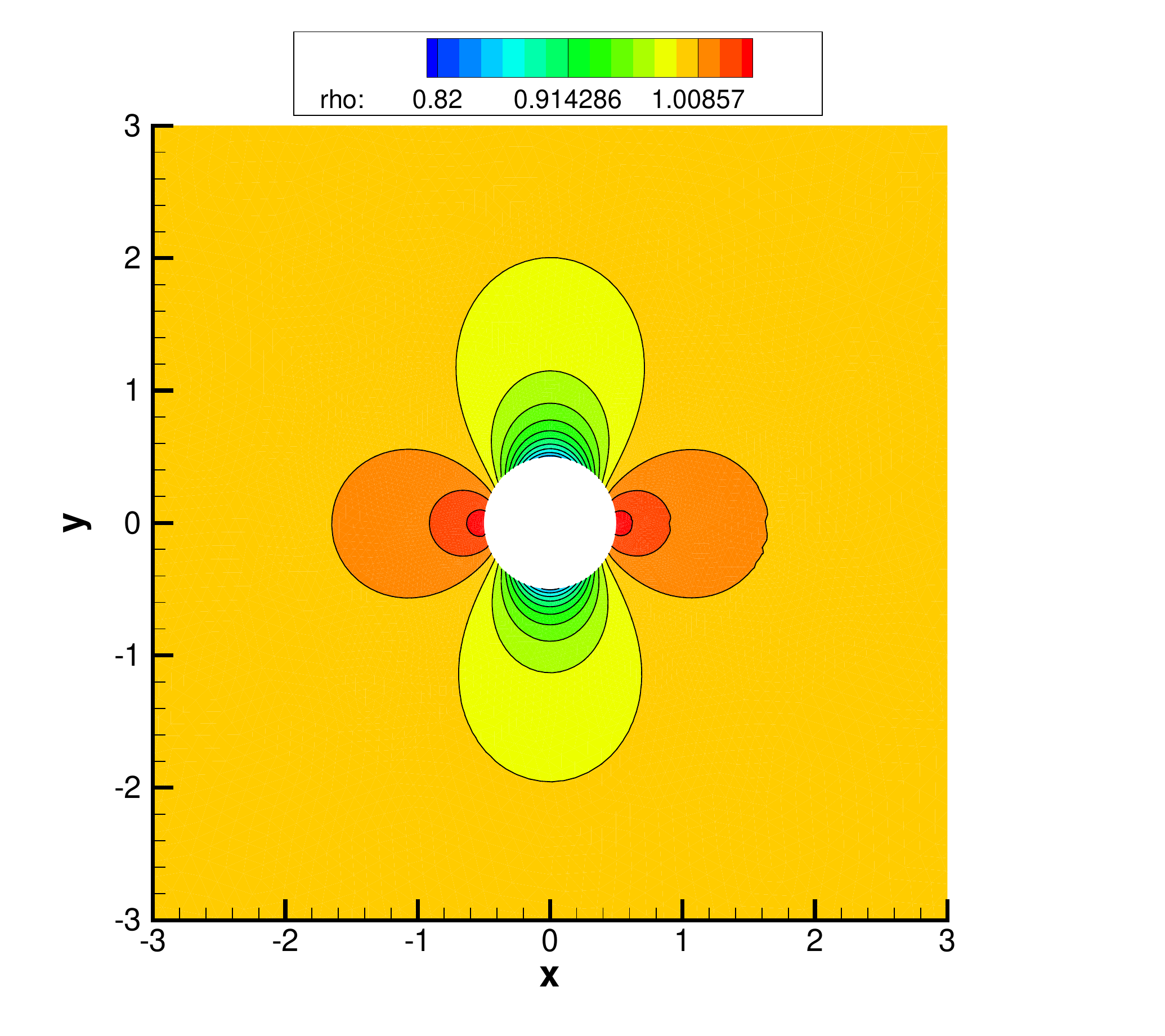}}
        \caption{Test case setup for the subsonic flow over a circular cylinder presented in Section~\ref{ssec.SubsonicCyl-2D}. We plot the mesh on the left, the density profile computed with DG-$\mathbb{P}3$ in the middle, and the density profile computed with DG-$\mathbb{P}3$/SBM-$\mathbb{P}3$ on the right.}
        \label{fig:SubCylrho}
\end{figure}

\begin{figure}
        \centering
        {\includegraphics[width=0.48\textwidth,trim={0cm 0cm 0cm 5.3cm},clip]{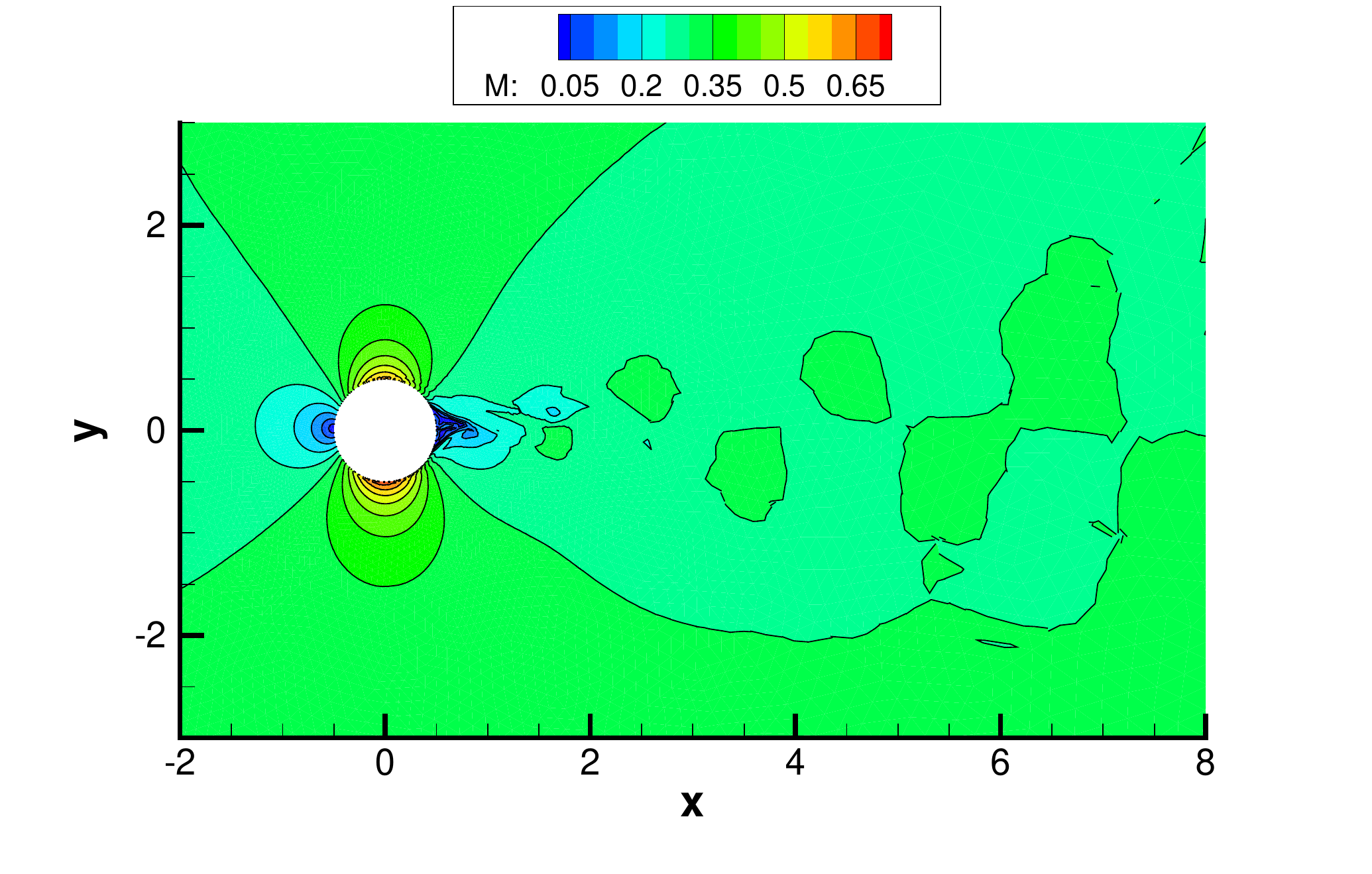}}
        {\includegraphics[width=0.48\textwidth,trim={0cm 0cm 0cm 5.3cm},clip]{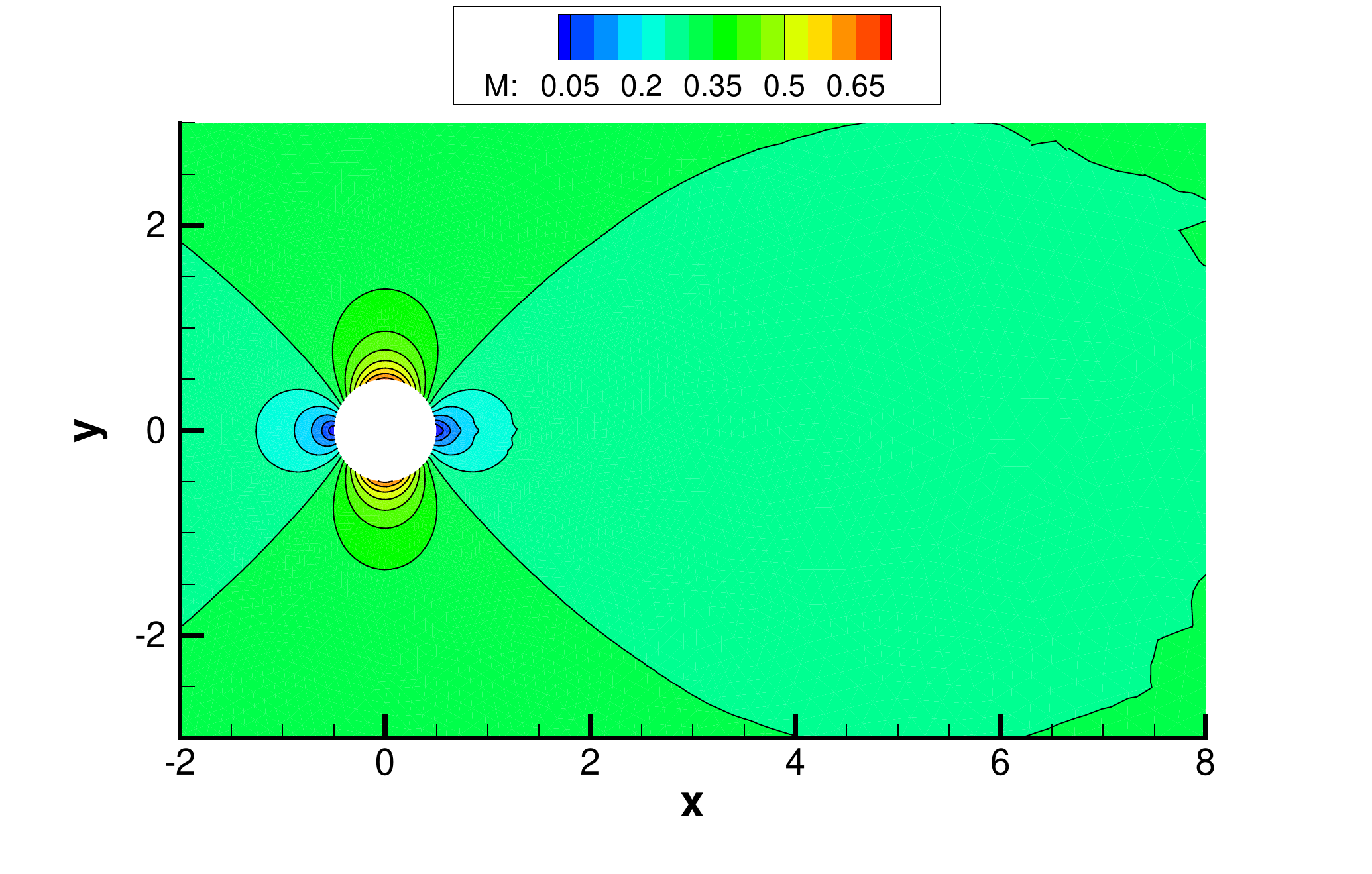}}
        \caption{Test case setup for the subsonic flow over a circular cylinder presented in Section~\ref{ssec.SubsonicCyl-2D}. We plot the Mach number profile computed with DG-$\mathbb{P}3$ on the left, and the Mach number profile computed with DG-$\mathbb{P}3$/SBM-$\mathbb{P}3$ on the right.}
        \label{fig:SubCylMach}
\end{figure}

\begin{figure}
        \centering
        {\includegraphics[width=0.48\textwidth,trim={0cm 0cm 0cm 0cm},clip]{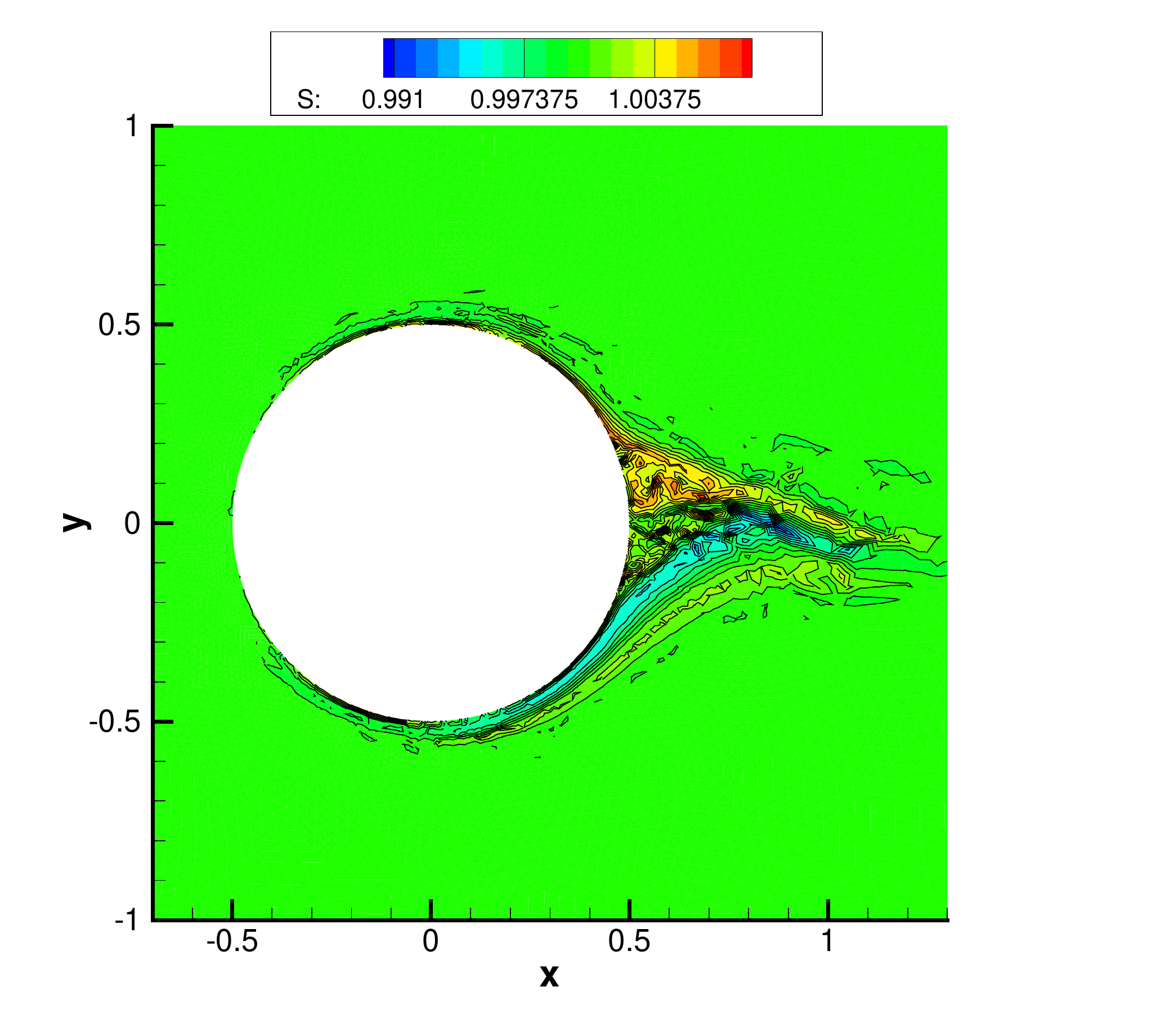}}
        {\includegraphics[width=0.48\textwidth,trim={0cm 0cm 0cm 0cm},clip]{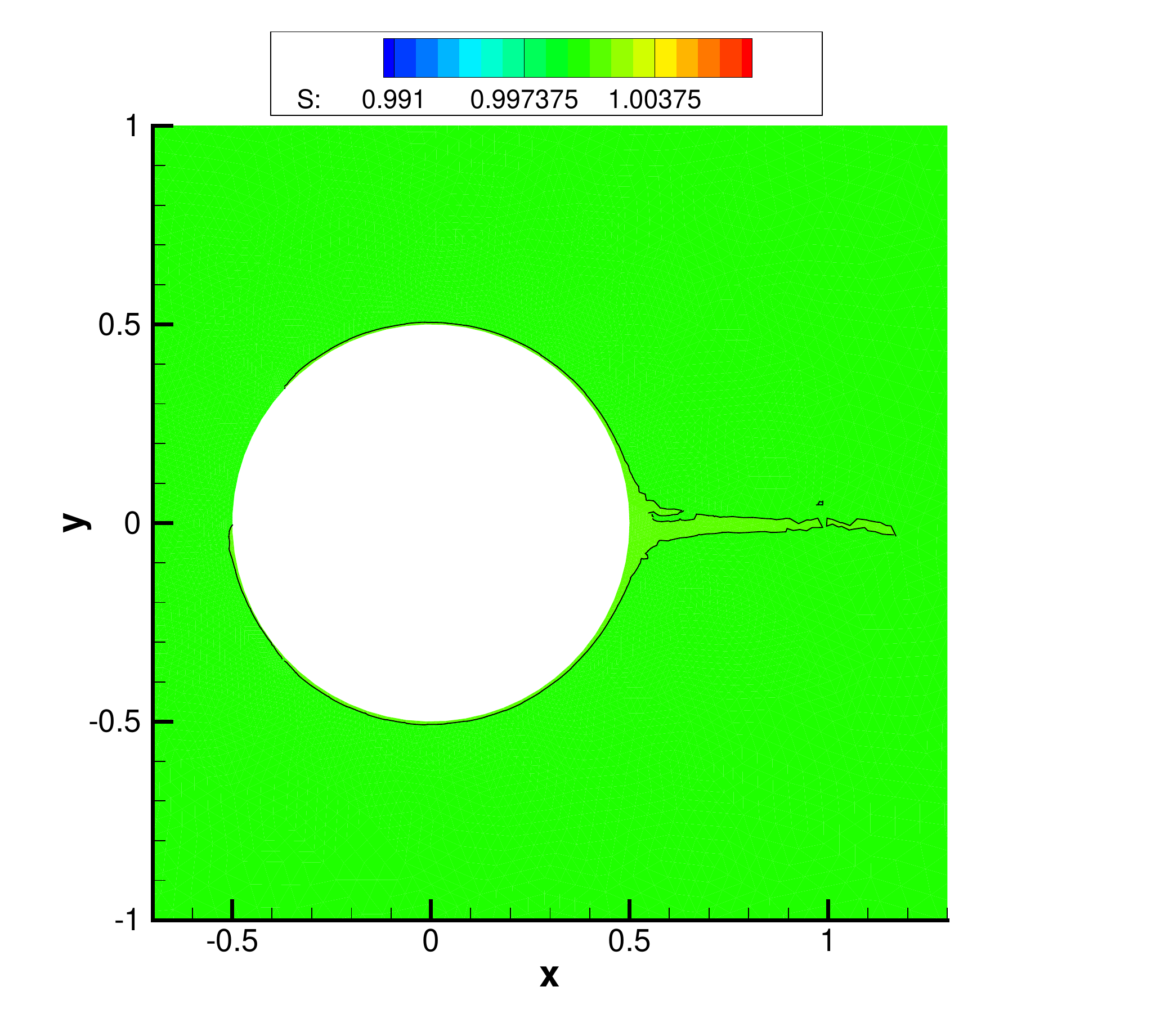}}
        \caption{Test case setup for the subsonic flow over a circular cylinder presented in Section~\ref{ssec.SubsonicCyl-2D}. We plot the entropy levels computed with DG-$\mathbb{P}3$ on the left, and the entropy levels computed with DG-$\mathbb{P}3$/SBM-$\mathbb{P}3$ on the right.}
        \label{fig:SubCylEntropy}
\end{figure}

}

\subsection{3D tests with smooth solutions}

We repeat the same study of the previous Section in the three-dimensional case.

\subsubsection{Manufactured solution on 3D curved domains: far-field BC} \label{ssec.Dir-3D}

\begin{table}
	\caption{Characteristics of the employed meshes for the test case of Section~\ref{ssec.Dir-3D}, the 3D manufactured solution.}
	\scriptsize
	\centering
	\begin{tabular}{cccc} \hline %\hline
		Grid level &Nodes  &Tetrahedra &$h$ \\[0.5mm]
		\hline
		0 &  169     &     778   &  2.8274E-01   \\
		1 &  971     &    5,072  &  1.4137E-01   \\
		2 &  6,437   &   35,968  &  7.0686E-02   \\
		\hline%\hline
	\end{tabular}
\end{table}

First, we consider the three-dimensional inhomogeneous Euler equations:
\begin{equation}
\mathbf{U}_t + \nabla\cdot \mathbf{F(U)} = \mathbf{S},\qquad \text{with}\qquad \mathbf{S} = \left(\begin{array}{c} 0.6\,\cos(x+y+z)\\ 0.8\,\cos(x+y+z)\\ 0.8\,\cos(x+y+z)\\ 0.8\,\cos(x+y+z)\\ 3.0\,\cos(x+y+z) \end{array}\right).
\end{equation}
This system has the following exact steady state solution, recovered using the \textit{manufactured solution technique},
\begin{equation}
\rho = 1+0.2\,\sin(x+y+z), \quad u = 1, \quad v = 1, \quad \omega = 1, \quad p = 1+0.2\,\sin(x+y+z),
\end{equation}
which is imposed on the domain boundaries as far-field boundary conditions.

\begin{figure}[t]
\centering
\subfigure[Grid level 0]{\includegraphics[width=0.32\textwidth,trim={2cm 2cm 2cm 2cm},clip]{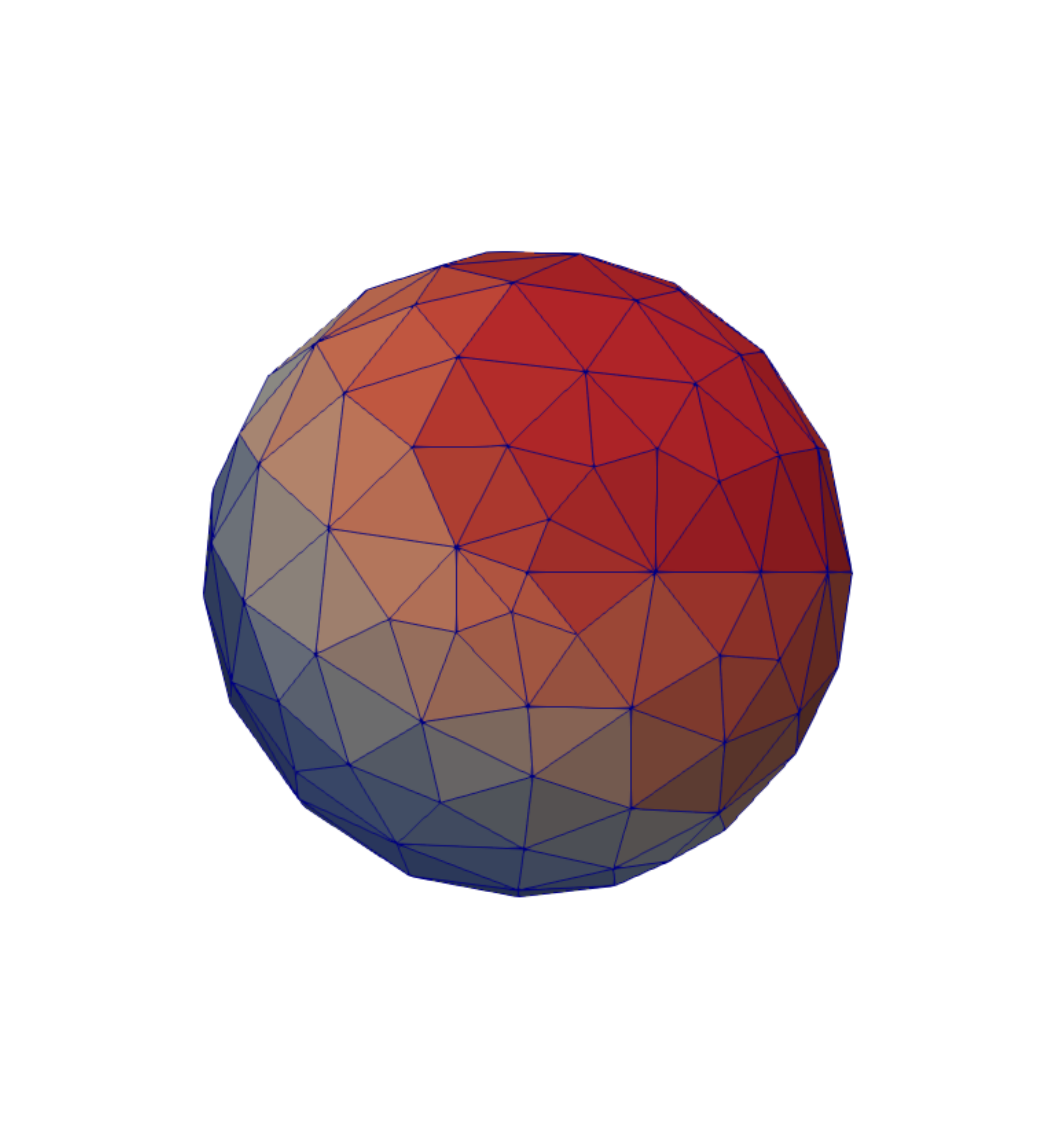}}
\subfigure[Grid level 1]{\includegraphics[width=0.32\textwidth,trim={2cm 2cm 2cm 2cm},clip]{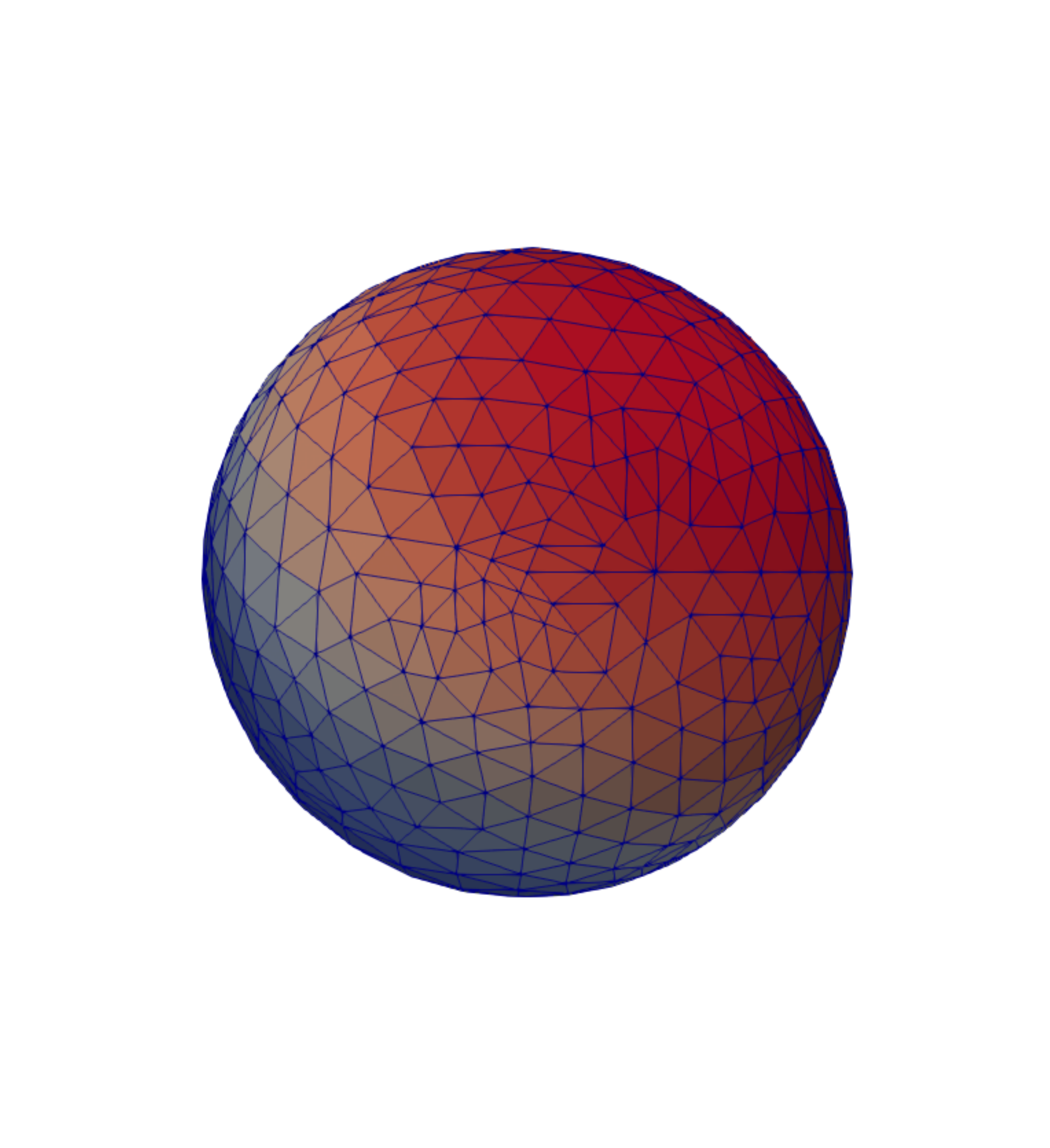}}
\subfigure[Grid level 2]{\includegraphics[width=0.32\textwidth,trim={2cm 2cm 2cm 2cm},clip]{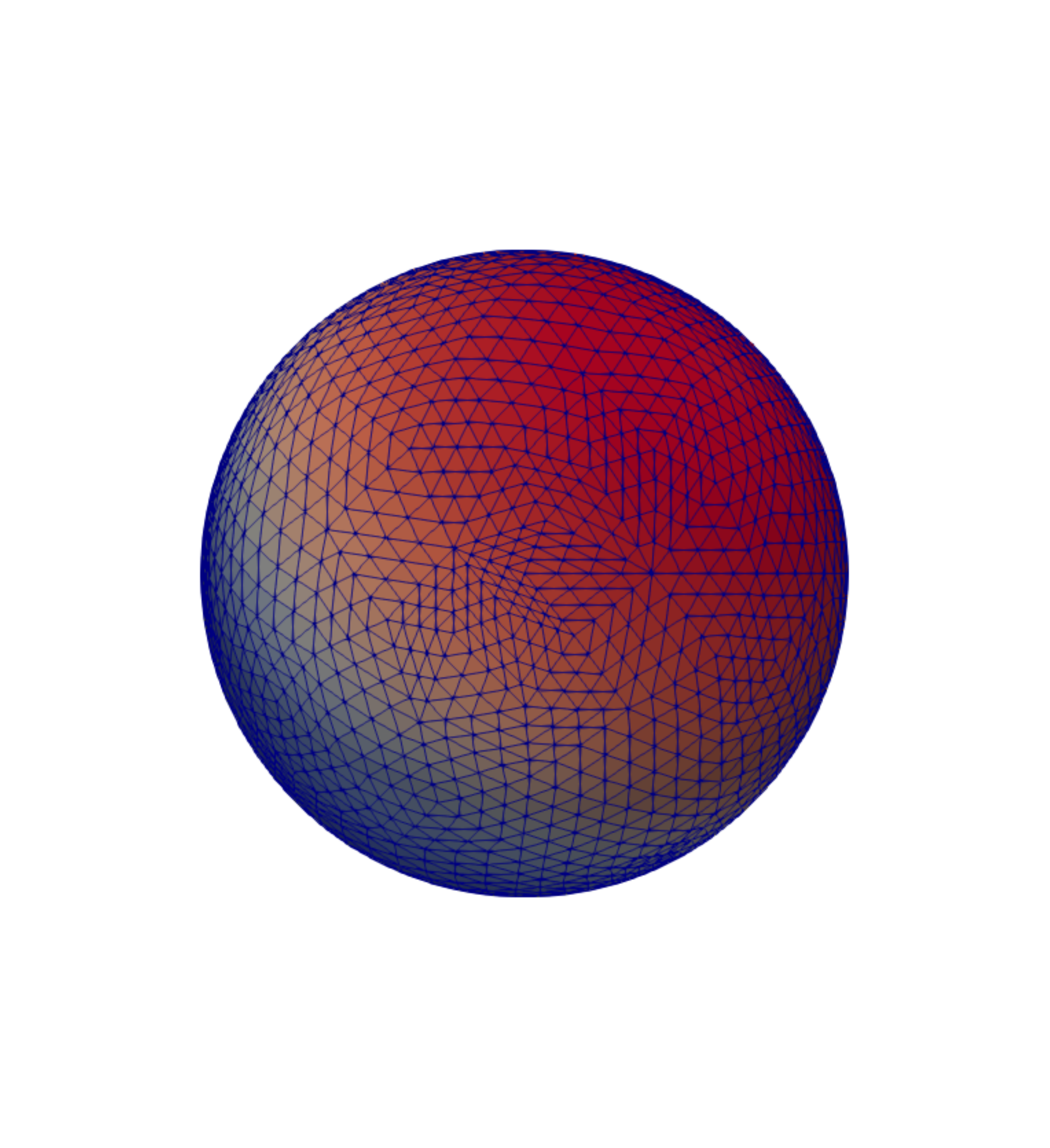}}
\caption{Test case setup for the 3D manufactured solution test case presented in Section~\ref{ssec.Dir-3D}.}\label{fig:manuf3Dsquare mesh}
\end{figure}

\begin{figure}[p]
\centering
\subfigure[]{\includegraphics[width=0.48\textwidth]{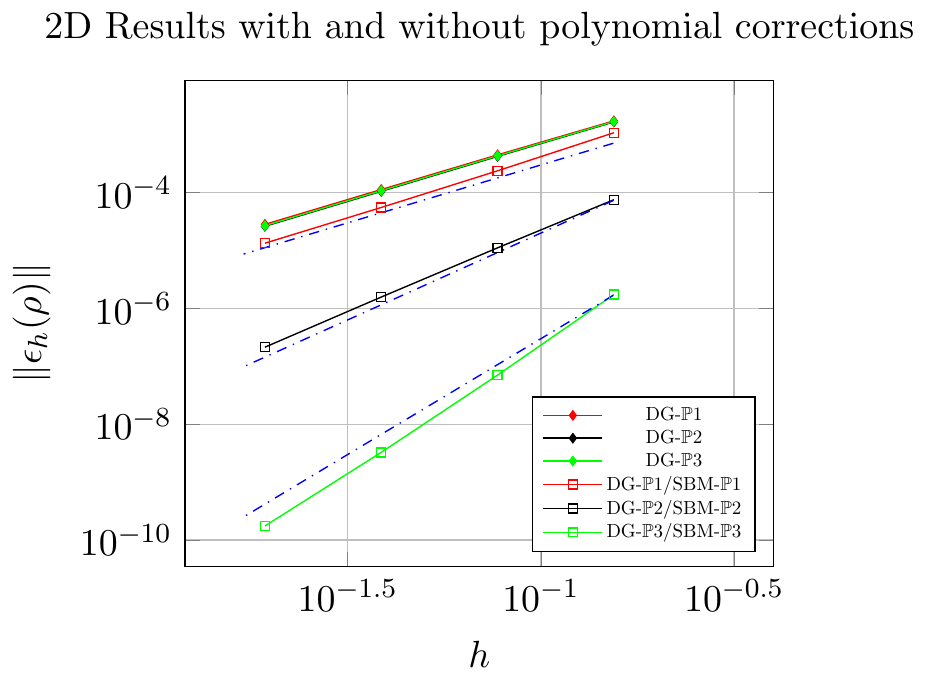}}
\subfigure[]{\includegraphics[width=0.47\textwidth]{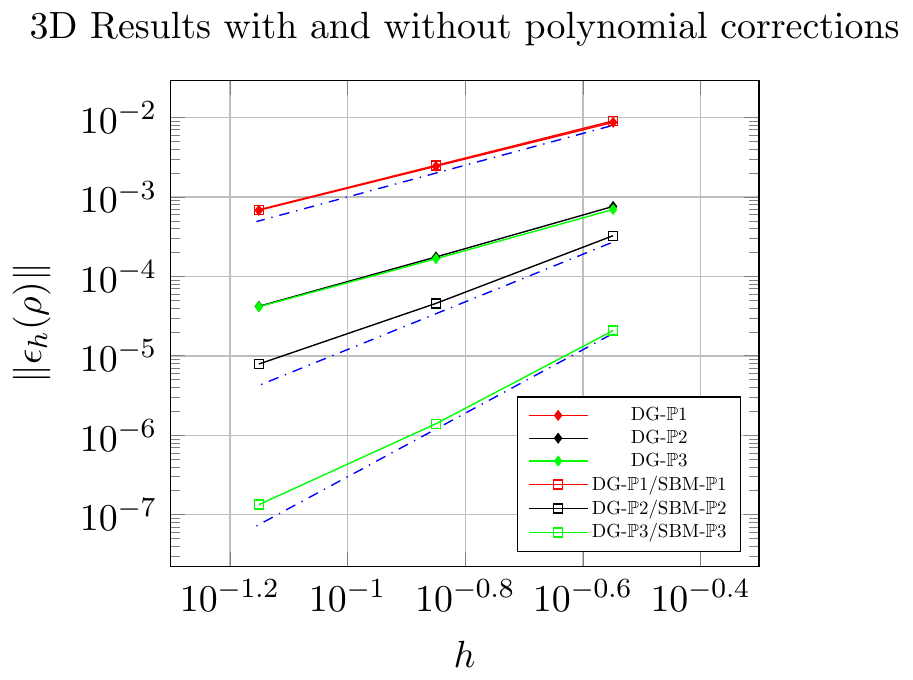}}
\caption{Manufactured solution (2D/3D): results obtained with and without the SBM correction.}\label{fig:dirichletConv}
\end{figure}

A sphere is now considered as boundary of the domain which introduces an error when using linear meshes due to the curvature. 
The employed linear meshes are shown in Figure~\ref{fig:manuf3Dsquare mesh}.
As expected when applying the far-field condition on the real curved boundary, the geometrical error given by the linear mesh
overcomes by far that given by the discretization technique with the outcome that no better than second order of accuracy can be achieved
(see Table~\ref{tb:manuf3D}).
{As before better results are obtained when using the SBM correction, indeed high order convergence is well recovered, 
as shown in Table~\ref{tb:manuf3D}.}
Again, convergence plots for the conserved variable $\rho$ are presented in Figure~\ref{fig:dirichletConv}b.

\begin{table}[p]
	\caption{Convergence analysis for the test case of Section~\ref{ssec.Dir-3D}, the 3D manufactured solution. We provide the results obtained without the SBM correction on linear meshes (top part) and with SBM correction on linear meshes (bottom part).}%\\}
	\label{tb:manuf3D}
	\scriptsize
	\centering
	\begin{tabular}{ccccccccccc} 
		\hline\\[1pt]
		\multicolumn{11}{c}{Convergence analysis \textit{without SBM} correction on \textit{linear} meshes}\\	
		\hline
		&\multicolumn{2}{c}{$\rho$} &\multicolumn{2}{c}{$\rho u$} &\multicolumn{2}{c}{$\rho v$} &\multicolumn{2}{c}{$\rho \omega$} &\multicolumn{2}{c}{$\rho E$}\\[0.5mm]
		\cline{2-11}
		Grid level & $L_2$        & $\tilde{n}$ & $L_2$        & $\tilde{n}$ & $L_2$        & $\tilde{n}$ & $L_2$      & $\tilde{n}$ \\[0.5mm]\hline
		&\multicolumn{10}{c}{DG-$\mathbb{P}1$}\\[0.5mm]
		0              & 8.7027E-3 &   --        & 8.8173E-3 &  --         & 8.8456E-3 &   --        & 8.7977E-3  &   --         & 3.3551E-2  &   --        \\
		1              & 2.4345E-3 &   1.84      & 2.4372E-3 &    1.86     & 2.4416E-3 &   1.86      & 2.4330E-3  &    1.85      & 9.3166E-3  &   1.85      \\
		2              & 6.7942E-4 &   1.84      & 6.7824E-4 &    1.85     & 6.8014E-4 &   1.84      & 6.7829E-4  &    1.84      & 2.6093E-3  &   1.84      \\
		&\multicolumn{10}{c}{DG-$\mathbb{P}2$}\\[0.5mm]
		0              & 7.6027E-4 &   --        & 9.1039E-4 &  --         & 8.9306E-4 &   --        & 8.4296E-4  &   --         & 3.2308E-3  &   --        \\
		1              & 1.7566E-4 &   2.11      & 2.1977E-4 &    2.05     & 2.1528E-4 &   2.05      & 1.9975E-4  &    2.08      & 7.8918E-4  &   2.03      \\
		2              & 4.2006E-5 &   2.06      & 5.4457E-5 &    2.01     & 5.3289E-5 &   2.01      & 4.8978E-5  &    2.03      & 1.9422E-4  &   2.02      \\
		&\multicolumn{10}{c}{DG-$\mathbb{P}3$}\\[0.5mm]
		0              & 6.9598E-4 &   --        & 8.7766E-4 &  --         & 8.5810E-4 &   --        & 7.9330E-4  &   --         & 3.1103E-3  &   --        \\
		1              & 1.6687E-4 &   2.06      & 2.1777E-4 &    2.01     & 2.1300E-4 &   2.01      & 1.9540E-4  &    2.02      & 7.7173E-4  &   2.01      \\
		2              & 4.1490E-5 &   2.01      & 5.4856E-5 &    1.99     & 5.3567E-5 &   1.99      & 4.9074E-5  &    1.99      & 1.9345E-4  &   2.00      \\
		\hline\\[1pt]
		\multicolumn{11}{c}{Convergence analysis \textit{with SBM} correction on \textit{linear} meshes}\\
		\hline
		&\multicolumn{2}{c}{$\rho$} &\multicolumn{2}{c}{$\rho u$} &\multicolumn{2}{c}{$\rho v$} &\multicolumn{2}{c}{$\rho \omega$} &\multicolumn{2}{c}{$\rho E$}\\[0.5mm]
		\cline{2-11}
		Grid level & $L_2$        & $\tilde{n}$ & $L_2$        & $\tilde{n}$ & $L_2$        & $\tilde{n}$ & $L_2$      & $\tilde{n}$ \\[0.5mm]\hline
		&\multicolumn{10}{c}{DG-$\mathbb{P}1$/SBM-$\mathbb{P}1$}\\[0.5mm]
		0              & 9.0767E-3 &   --        & 9.1066E-3 &  --         & 9.1410E-3 &   --        & 9.1148E-3  &   --         & 3.5055E-2  &   --        \\
		1              & 2.4981E-3 &   1.86      & 2.4841E-3 &   1.87      & 2.4895E-3 &   1.88      & 2.4859E-3  &   1.87       & 9.5724E-3  &   1.87      \\
		2              & 6.9067E-4 &   1.85      & 6.8601E-4 &   1.86      & 6.8798E-4 &   1.86      & 6.8738E-4  &   1.85       & 2.6541E-3  &   1.85      \\
		&\multicolumn{10}{c}{DG-$\mathbb{P}2$/SBM-$\mathbb{P}2$}\\[0.5mm]
		0              & 3.2537E-4 &   --        & 3.3772E-4 &  --         & 3.4049E-4 &   --        & 3.3983E-4  &   --         & 1.2000E-3  &   --        \\
		1              & 4.5765E-5 &   2.83      & 4.5417E-5 &   2.89      & 4.5165E-5 &   2.91      & 4.5562E-5  &   2.90       & 1.8400E-4  &   2.71      \\
		2              & 7.8831E-6 &   2.54      & 7.4985E-6 &   2.60      & 7.4576E-6 &   2.60      & 7.5288E-6  &   2.60       & 2.9674E-5  &   2.63      \\
		&\multicolumn{10}{c}{DG-$\mathbb{P}3$/SBM-$\mathbb{P}3$}\\[0.5mm]
		0              & 2.0810E-5 &   --        & 2.3072E-5 &  --         & 2.3226E-5 &   --        & 2.2280E-5  &   --         & 8.4229E-5  &   --        \\
		1              & 1.3956E-6 &   3.90      & 1.4671E-6 &   3.98      & 1.4566E-6 &   4.00      & 1.4633E-6  &   3.93       & 5.5552E-6  &   3.92      \\
		2              & 1.3372E-7 &   3.38      & 1.3562E-7 &   3.44      & 1.3579E-7 &   3.42      & 1.3562E-7  &   3.43       & 5.1363E-7  &   3.44      \\
		\hline
	\end{tabular}
\end{table}

\subsubsection{Supersonic vortex bounded by two cylindrical walls: slip wall BC} \label{ssec.Wall-3D}

We consider a 3D isentropic supersonic flow between two concentric cylindrical surfaces of radii $r_i = 1$ and
$r_o = 1.384$. The exact density, velocity and pressure in terms of radius $r$ are given by Eqs.~\eqref{eq:vortex2D rho} and~\eqref{eq:vortex2D u p} and the solution in the inner surface is taken as that in Section~\ref{ssec.Wall-2D}.
The fluid's velocity vector components in $(x,y,z)$ can be computed as follows:
\begin{equation}
\left( \begin{array}{c} u \\ v \\ \omega \end{array} \right) = \|\mathbf{u}\| \left( \begin{array}{c} y/r \\ -x/r \\ 0\end{array} \right),
\end{equation}
where $r=\sqrt{x^2+y^2}$ because the cylindrical surfaces are developed along the z--axis.

\begin{figure}
\centering
\subfigure[]{\includegraphics[width=0.48\textwidth]{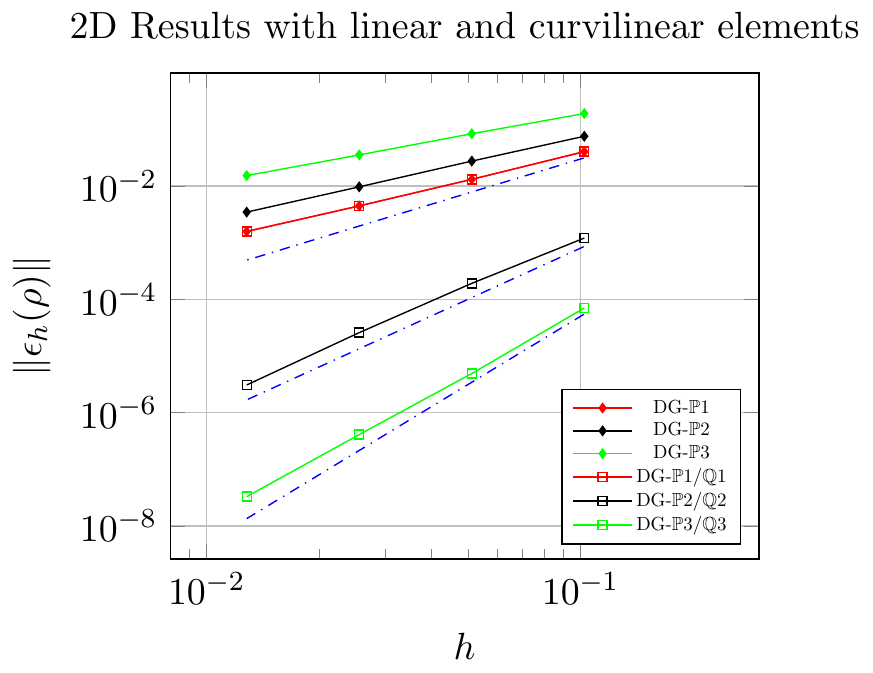}}
\subfigure[]{\includegraphics[width=0.46\textwidth]{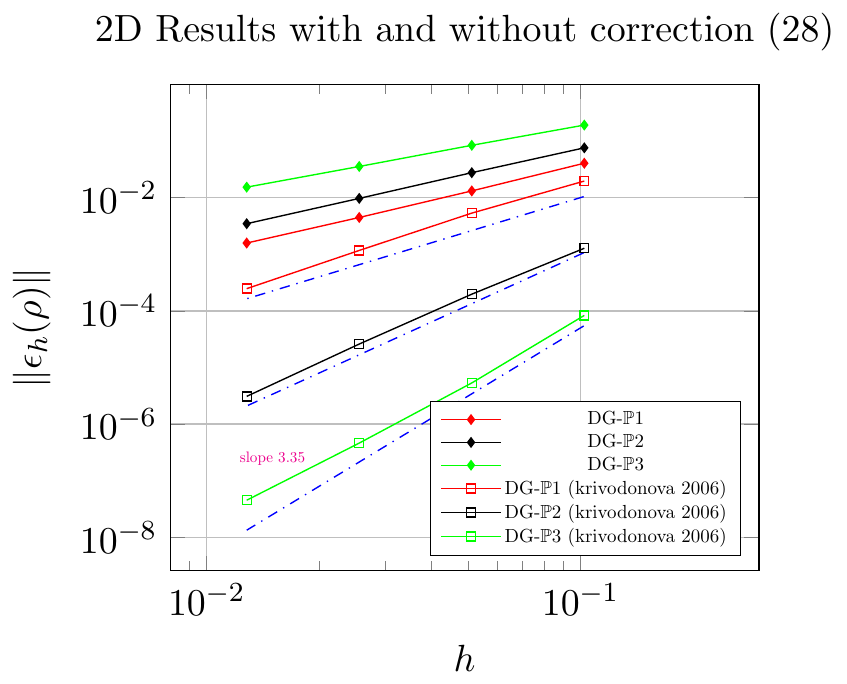}}
\subfigure[]{\includegraphics[width=0.48\textwidth]{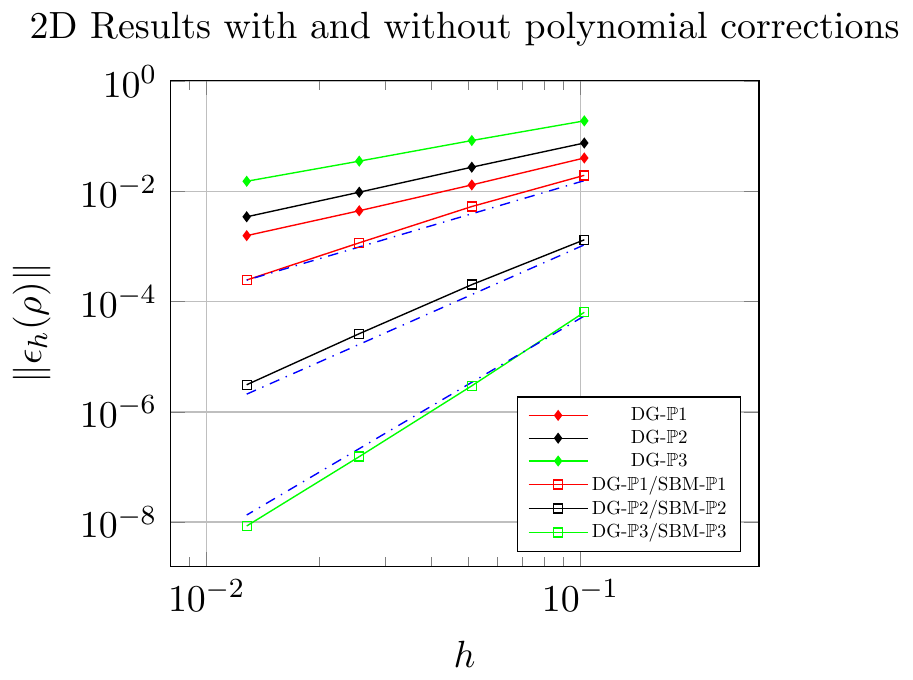}}
\subfigure[]{\includegraphics[width=0.48\textwidth]{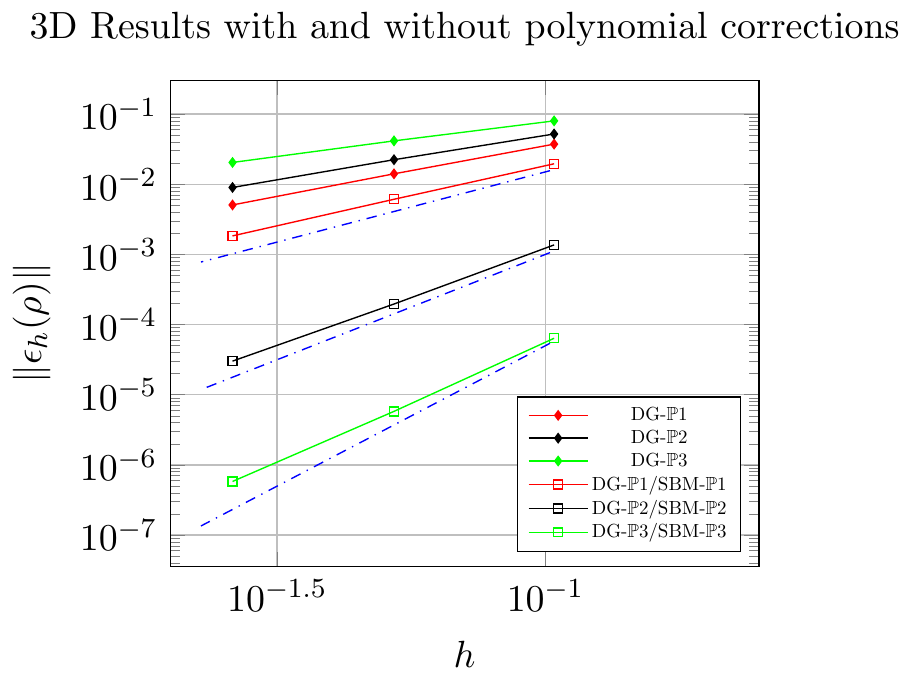}}
\caption{Supersonic vortex bounded by two walls (2D/3D): convergence tests performed with (a) linear and curvilinear 2D grids; (b) wall correction provided in~\eqref{eq:BergerFlux};
(c) polynomial corrections for 2D grids; (d) polynomial corrections for 3D grids.}\label{fig:slipwallConv}
\end{figure}

Simulations are first run with classical reflecting wall boundary conditions applied on the approximated boundary.
The results presented in Table~\ref{tb:SlipWall3D} show convergence trends of rates between 1.5 and 1 for all high order polynomials.
Finally, Table~\ref{tb:SlipWall3D} presents the results obtained using the SBM wall flux correction:
already a little improvement is shown for $\mathbb{P}1$ and for higher order all convergence trends are properly recovered 
(see Figure~\ref{fig:slipwallConv}d for the convergences plot for the density variable $\rho$).

\begin{figure}
\centering
\subfigure[Grid level 0]{\includegraphics[width=0.32\textwidth,trim={0cm 0cm 0cm 0cm},clip]{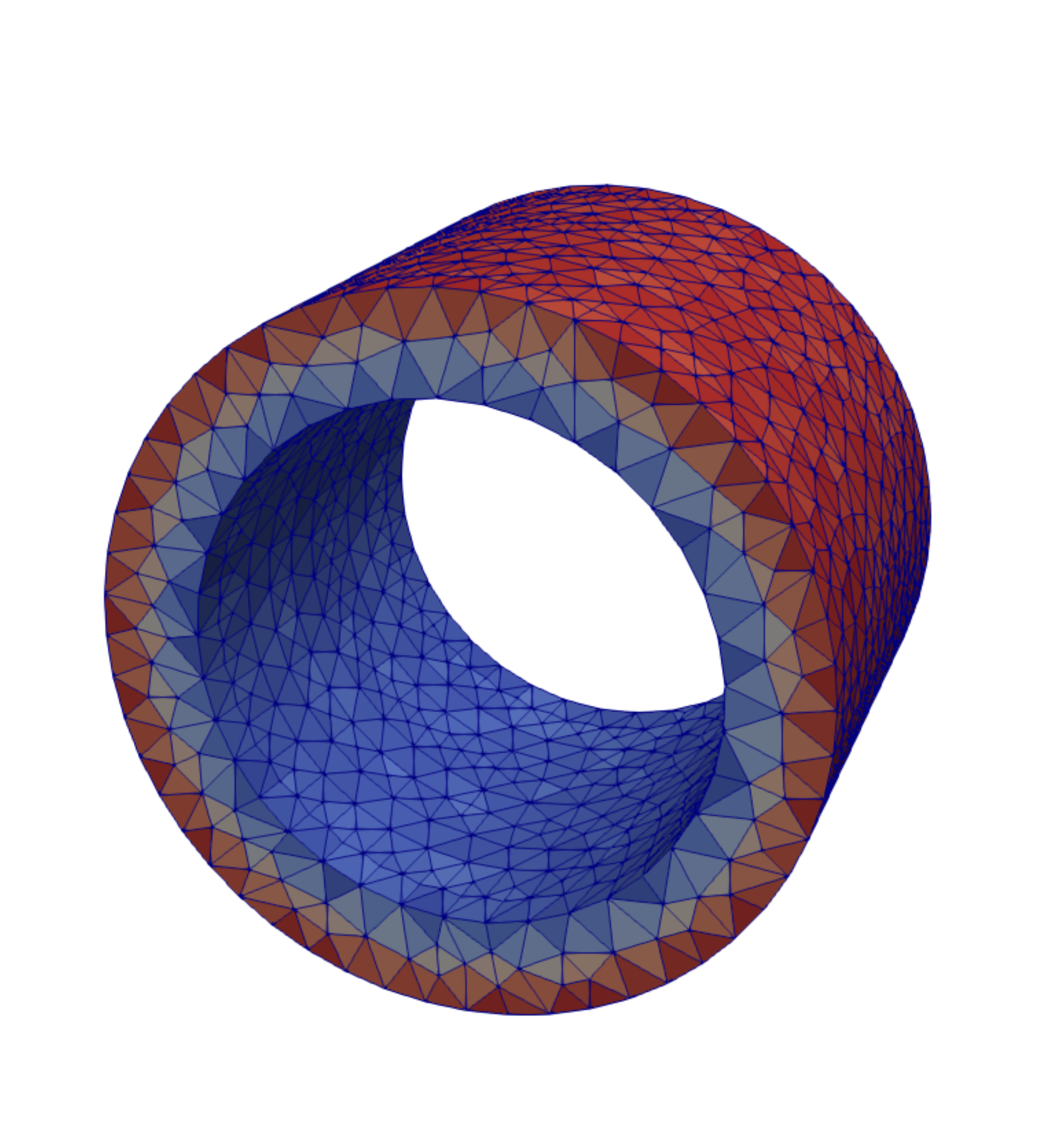}}
\subfigure[Grid level 1]{\includegraphics[width=0.32\textwidth,trim={0cm 0cm 0cm 0cm},clip]{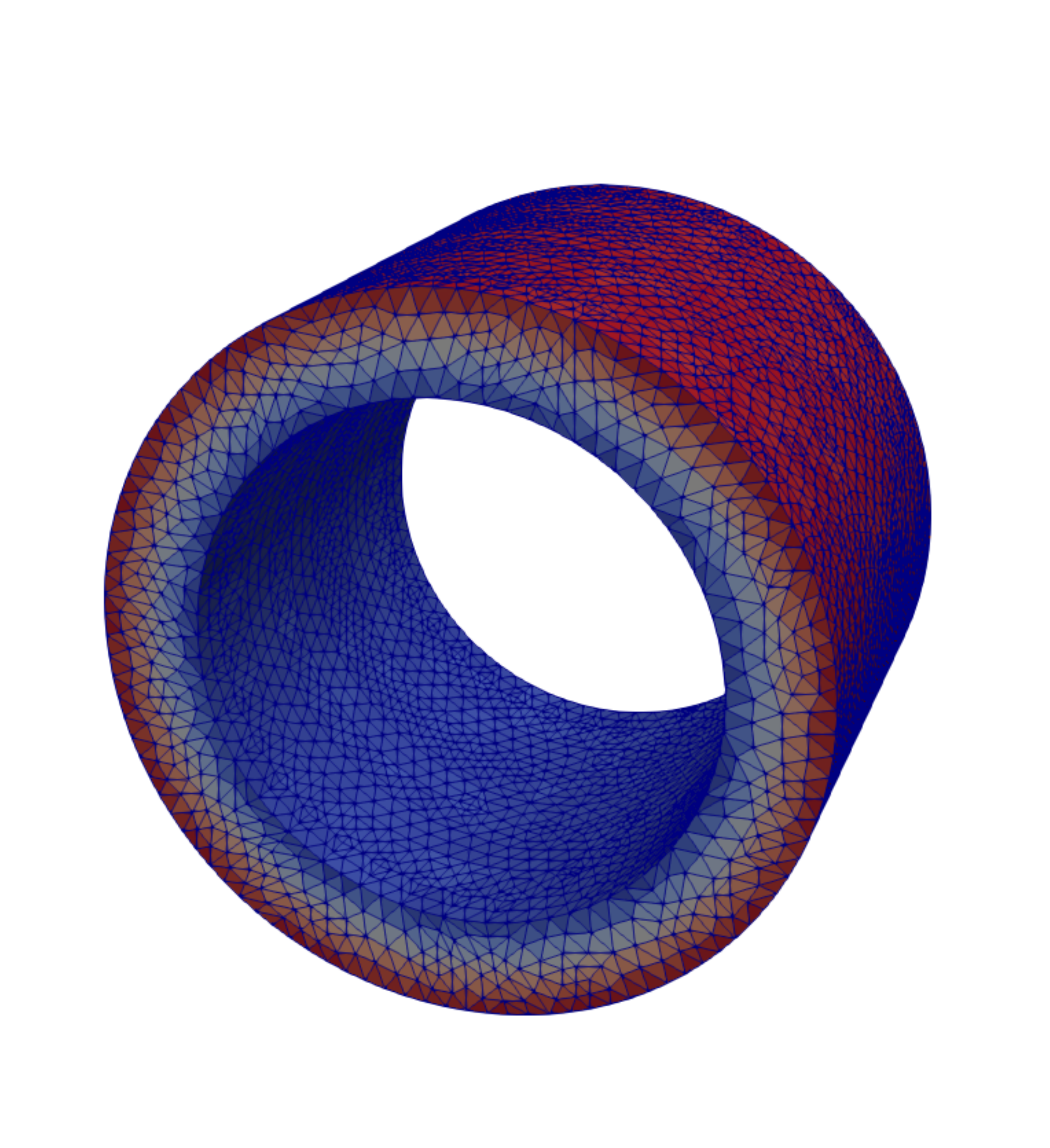}}
\subfigure[Grid level 2]{\includegraphics[width=0.32\textwidth,trim={0cm 0cm 0cm 0cm},clip]{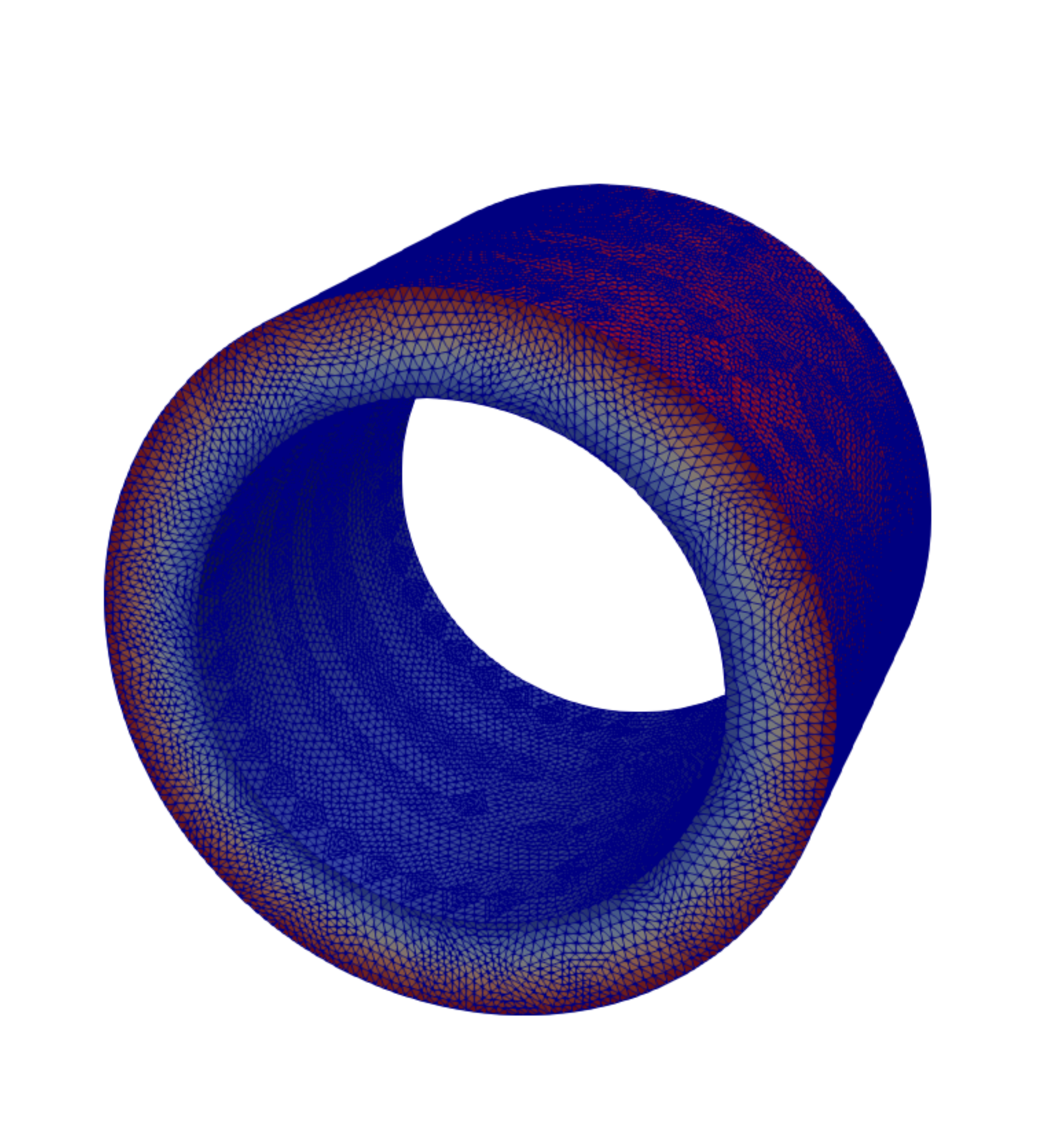}}
\caption{Supersonic vortex bounded by two cylindrical walls (3D): density contours.}\label{fig:bergerVortex3D mesh}
\end{figure}

\begin{table}[h!]
\caption{Characteristics of the employed meshes for the test case of Section~\ref{ssec.Wall-3D}, the supersonic vortex bounded by two cylindrical walls in 3D.}
\scriptsize
\centering
\begin{tabular}{cccc} \hline%\hline
Grid level &Nodes  &Tetrahedra &$h$ \\[0.5mm]
\hline
0 &   3,071  &   13,059  &   1.0382E-01 \\
1 &  20,929  &  104,472  &   5.2221E-02 \\
2 & 153,242  &  835,776  &   2.6112E-02 \\
\hline%\hline
\end{tabular}
\end{table}

\begin{table}[h!]
	\caption{Convergence analysis for the test case of Section~\ref{ssec.Wall-3D}, the supersonic vortex bounded by two cylindrical walls in 3D. We provide the results obtained without the SBM correction on linear meshes (top part) and with SBM correction on linear meshes (bottom part).}%\\}
	\label{tb:SlipWall3D}
	\scriptsize
	\centering
	\begin{tabular}{ccccccccccc} 
		\hline\\[1pt]
		\multicolumn{11}{c}{Convergence analysis \textit{without SBM} correction on \textit{linear} meshes}\\	
		\hline
		&\multicolumn{2}{c}{$\rho$} &\multicolumn{2}{c}{$\rho u$} &\multicolumn{2}{c}{$\rho v$} &\multicolumn{2}{c}{$\rho \omega$} &\multicolumn{2}{c}{$\rho E$}\\[0.5mm]
		\cline{2-11}
		Grid level & $L_2$        & $\tilde{n}$ & $L_2$        & $\tilde{n}$ & $L_2$        & $\tilde{n}$ & $L_2$      & $\tilde{n}$ \\[0.5mm]\hline
		&\multicolumn{10}{c}{DG-$\mathbb{P}1$}\\[0.5mm]
		0              & 3.7322E-2 &   --        & 6.4933E-2 &  --         & 6.5061E-2 &   --        & 1.1475E-2  &   --         & 1.6600E-1  &   --        \\
		1              & 1.4081E-2 &   1.41      & 2.0592E-2 &   1.65      & 2.0796E-2 &    1.65     & 5.0424E-3  &   1.18       & 6.0193E-2  &   1.46      \\
		2              & 5.0733E-3 &   1.47      & 6.6189E-3 &   1.64      & 6.7225E-3 &    1.63     & 1.9803E-3  &   1.35       & 2.1242E-2  &   1.50      \\
		&\multicolumn{10}{c}{DG-$\mathbb{P}2$}\\[0.5mm]
		0              & 5.2040E-2 &   --        & 5.9910E-2 &  --         & 6.1607E-2 &   --        & 1.6668E-2  &   --         & 2.0989E-1  &   --        \\
		1              & 2.2376E-2 &    1.22     & 2.3548E-2 &    1.35     & 2.4191E-2 &   1.35      & 8.6415E-3  &   0.95       & 8.6540E-2  &   1.28      \\
		2              & 9.0045E-3 &    1.31     & 8.8598E-3 &    1.41     & 9.1564E-3 &   1.40      & 3.8322E-3  &   1.17       & 3.4325E-2  &   1.33      \\
		&\multicolumn{10}{c}{DG-$\mathbb{P}3$}\\[0.5mm]
		0              & 8.0027E-2 &   --        & 8.787E-2  &  --         & 8.9290E-2 &   --        & 2.5696E-2  &   --         & 3.0165E-1  &   --        \\
		1              & 4.1512E-2 &   0.94      & 3.931E-2  &   1.16      & 4.0139E-2 &   1.15      & 1.5240E-2  &   0.75       & 1.5087E-1  &   1.00      \\
		2              & 2.0451E-2 &   1.02      & 1.729E-2  &   1.18      & 1.7731E-2 &   1.18      & 8.3006E-3  &   0.87       & 7.3682E-2  &   1.03      \\
		\hline\\[1pt]
		\multicolumn{11}{c}{Convergence analysis \textit{with SBM} correction on \textit{linear} meshes}\\
		\hline
		&\multicolumn{2}{c}{$\rho$} &\multicolumn{2}{c}{$\rho u$} &\multicolumn{2}{c}{$\rho v$} &\multicolumn{2}{c}{$\rho \omega$} &\multicolumn{2}{c}{$\rho E$}\\[0.5mm]
		\cline{2-11}
		Grid level & $L_2$        & $\tilde{n}$ & $L_2$        & $\tilde{n}$ & $L_2$        & $\tilde{n}$ & $L_2$      & $\tilde{n}$ \\[0.5mm]\hline
		&\multicolumn{10}{c}{DG-$\mathbb{P}1$/SBM-$\mathbb{P}1$}\\[0.5mm]
		0              & 1.9580E-2 &   --        & 5.4668E-2 &  --         & 5.3611E-2 &   --        & 5.8383E-3  &   --         & 9.1443E-2  &   --        \\
		1              & 6.1094E-3 &    1.68     & 1.5653E-2 &    1.80     & 1.5334E-2 &    1.80     & 2.0797E-3  &    1.49      & 2.7820E-2  &   1.72      \\
		2              & 1.8376E-3 &    1.73     & 4.4901E-3 &    1.80     & 4.4063E-3 &    1.80     & 7.0284E-4  &    1.57      & 8.2183E-3  &   1.76      \\
		&\multicolumn{10}{c}{DG-$\mathbb{P}2$/SBM-$\mathbb{P}2$}\\[0.5mm]
		0              & 1.3668E-3 &   --        & 2.4856E-3 &  --         & 2.4388E-3 &   --        & 6.4125E-4  &   --         & 6.3670E-3  &   --        \\
		1              & 1.9656E-4 &   2.80      & 3.6385E-4 &    2.77     & 3.5349E-4 &   2.79      & 9.7458E-5  &   2.72       & 9.2034E-4  &   2.79      \\
		2              & 3.0200E-5 &   2.70      & 5.3739E-5 &    2.76     & 5.2281E-5 &   2.75      & 1.5670E-5  &   2.64       & 1.3566E-4  &   2.76      \\
		&\multicolumn{10}{c}{DG-$\mathbb{P}3$/SBM-$\mathbb{P}3$}\\[0.5mm]
		0              & 6.3903E-5 &   --        & 1.0992E-4 &  --         & 1.0858E-4 &   --        &  2.7024E-5 &   --         & 2.6291E-4 &   --        \\
		1              & 5.8075E-6 &   3.46      & 1.0098E-5 &   3.44      & 9.9331E-6 &    3.45     &  2.3153E-6 &    3.54      & 2.3666E-5 &   3.47      \\
		2              & 5.8223E-7 &   3.32      & 1.0373E-6 &   3.28      & 1.0171E-6 &    3.29     &  2.3309E-7 &    3.31      & 2.3781E-6 &   3.31      \\
		\hline
	\end{tabular}
\end{table}

%\subsection{Wall boundary conditions in 2D: shock-cylinder interaction} \label{ssec.shock-cyl2D interaction}
\subsection{{Shock-cylinder interaction}} \label{ssec.shock-cyl2D interaction}

\begin{figure}
\centering
{\includegraphics[width=0.5\textwidth,trim={0cm 0cm 0cm 4.5cm},clip]{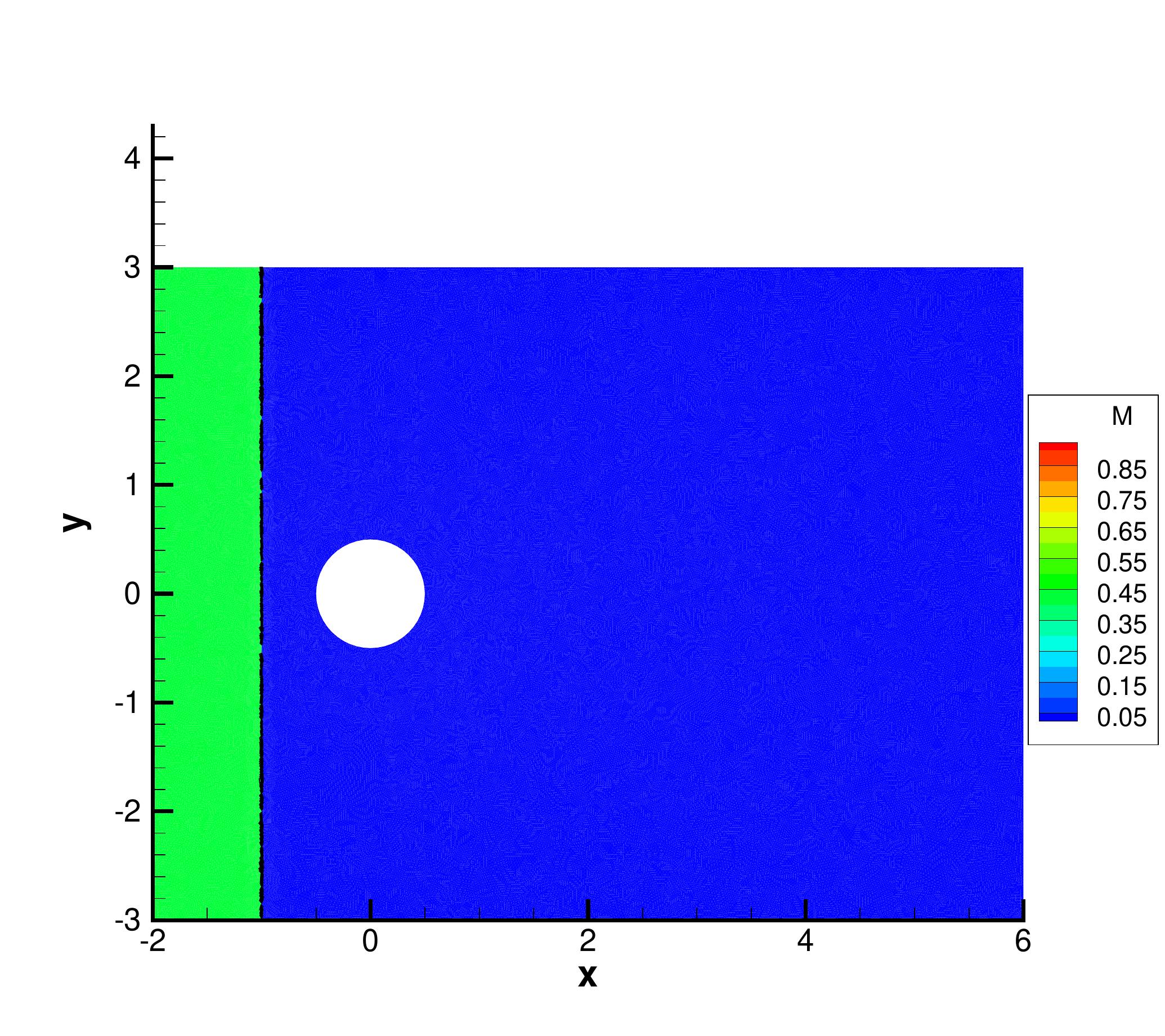}}
%{\includegraphics[width=0.45\textwidth,trim={0cm 0cm 0cm 4.5cm},clip]{ShockCylinder_Limiter_0.eps}}
%\caption{Test case setup for the shock-cylinder interaction of Section~\ref{ssec.shock-cyl2D interaction}. We show the initial density profile on the left, and on the right, in red, the cells crossed by the shock, thus detected as troubled, at time $t=0$.}\label{fig:ShockLimiter0}
\caption{Test case setup for the shock-cylinder interaction of Section~\ref{ssec.shock-cyl2D interaction}. We show the initial \Rall{Mach number} profile at time $t=0$.}\label{fig:ShockLimiter0}
\end{figure}

\begin{figure}
\centering
\subfigure[$t=0.5$]{\includegraphics[width=0.4\textwidth,trim={0cm 0cm 0cm 4.6cm},clip]{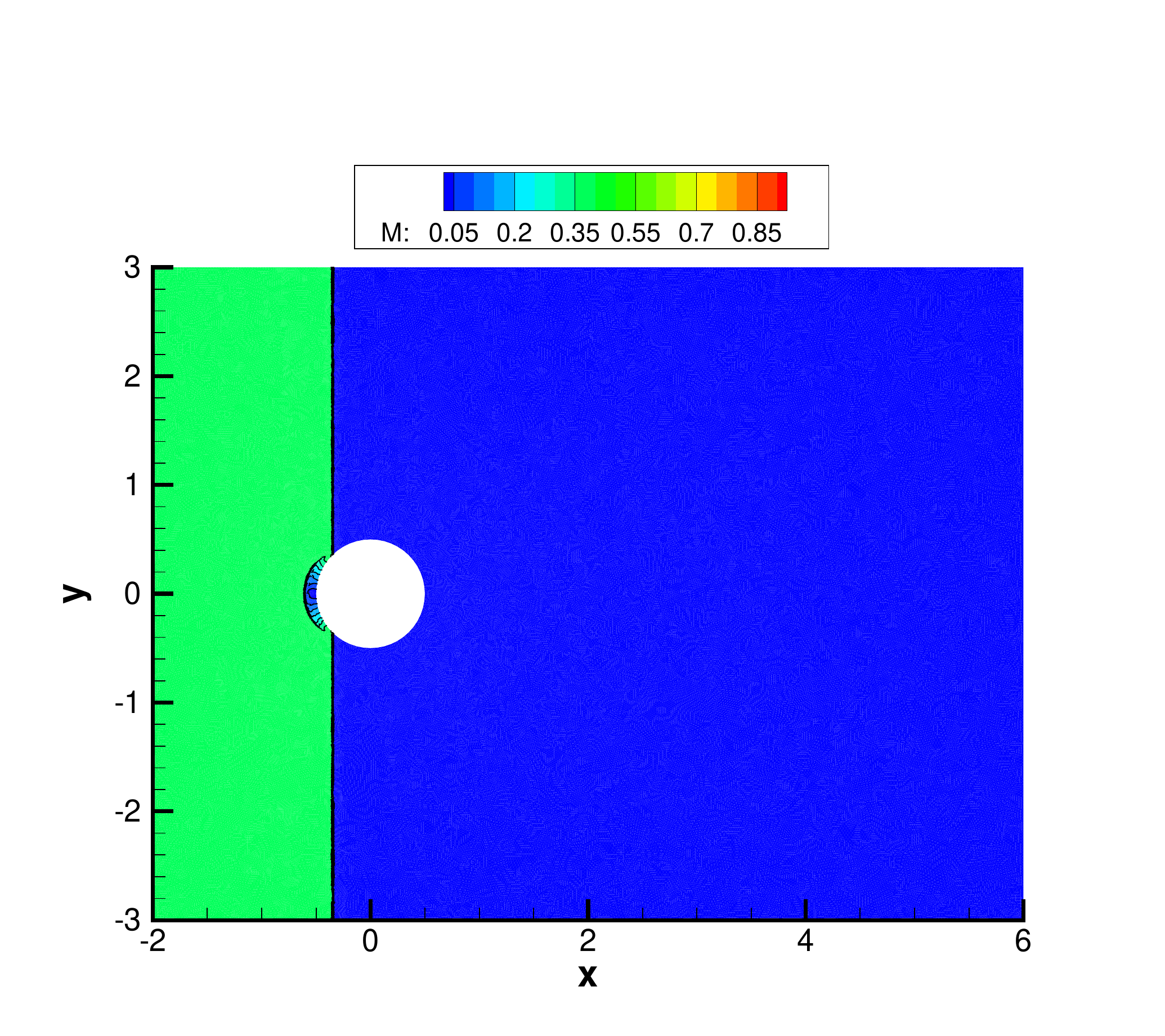}}
\subfigure[$t=0.5$]{\includegraphics[width=0.4\textwidth,trim={0cm 0cm 0cm 4.6cm},clip]{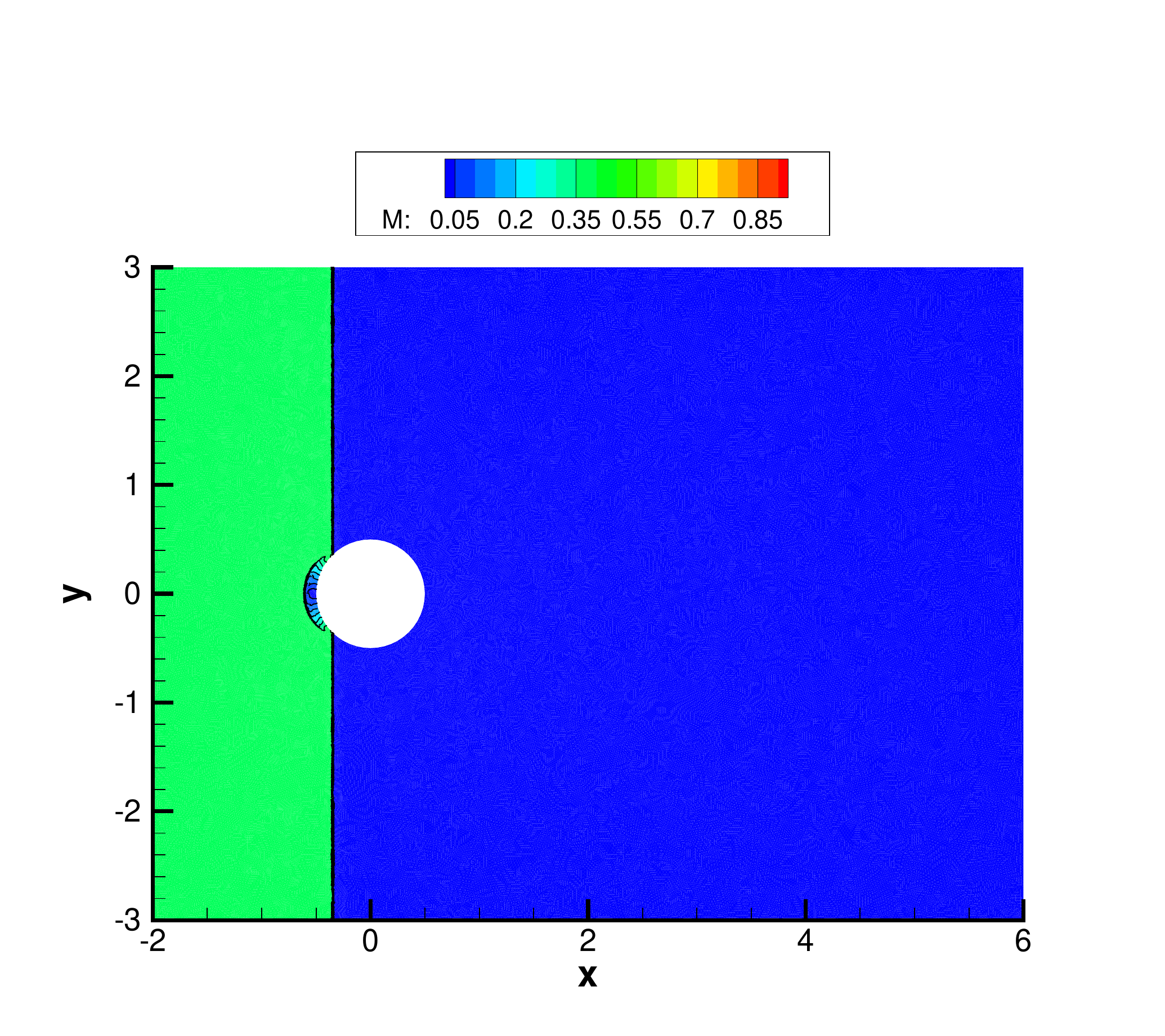}}
\qquad
\subfigure[$t=1$]  {\includegraphics[width=0.4\textwidth,trim={0cm 0cm 0cm 4.6cm},clip]{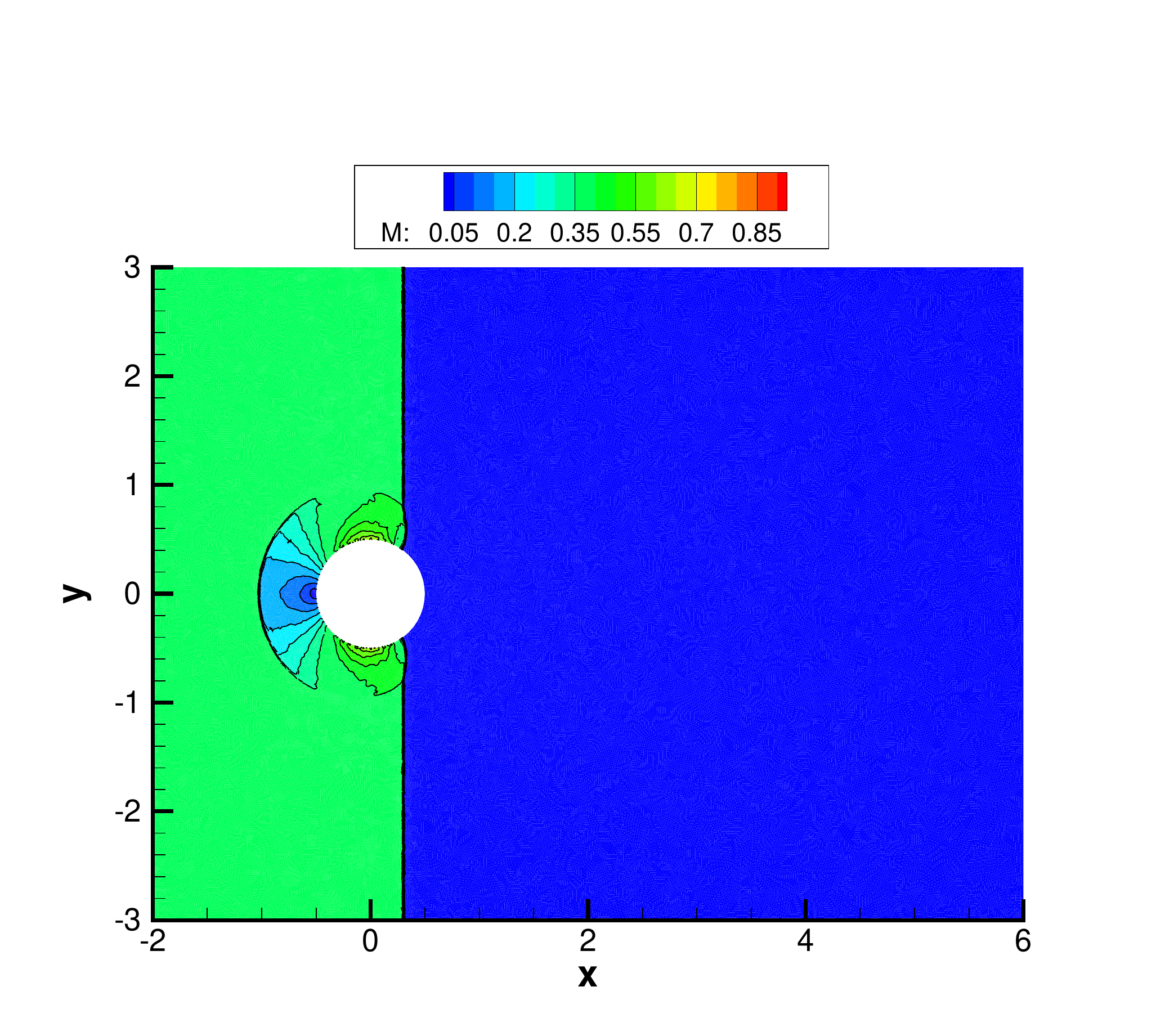}}
\subfigure[$t=1$]  {\includegraphics[width=0.4\textwidth,trim={0cm 0cm 0cm 4.6cm},clip]{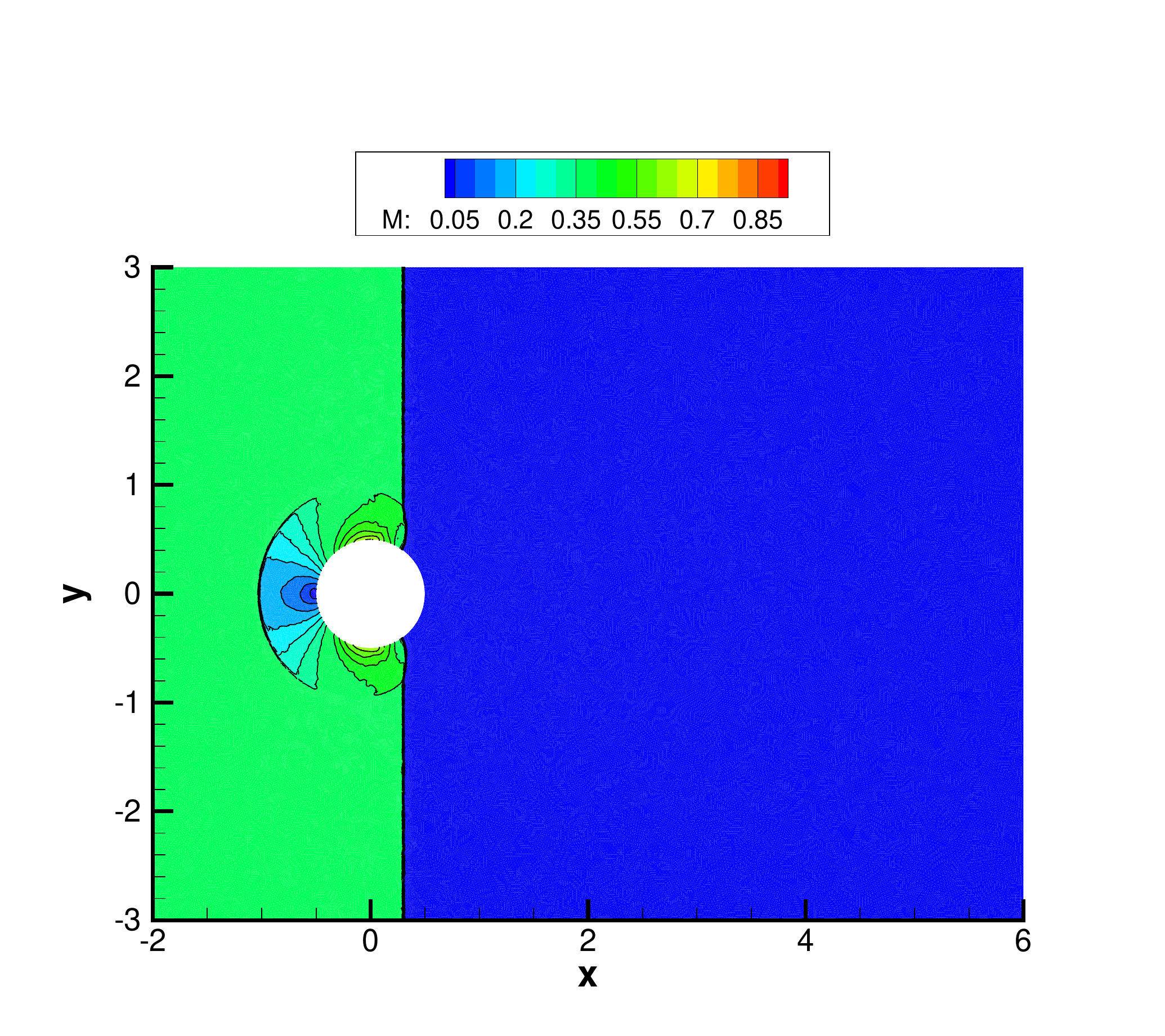}}
\qquad
\subfigure[$t=1.5$]{\includegraphics[width=0.4\textwidth,trim={0cm 0cm 0cm 4.6cm},clip]{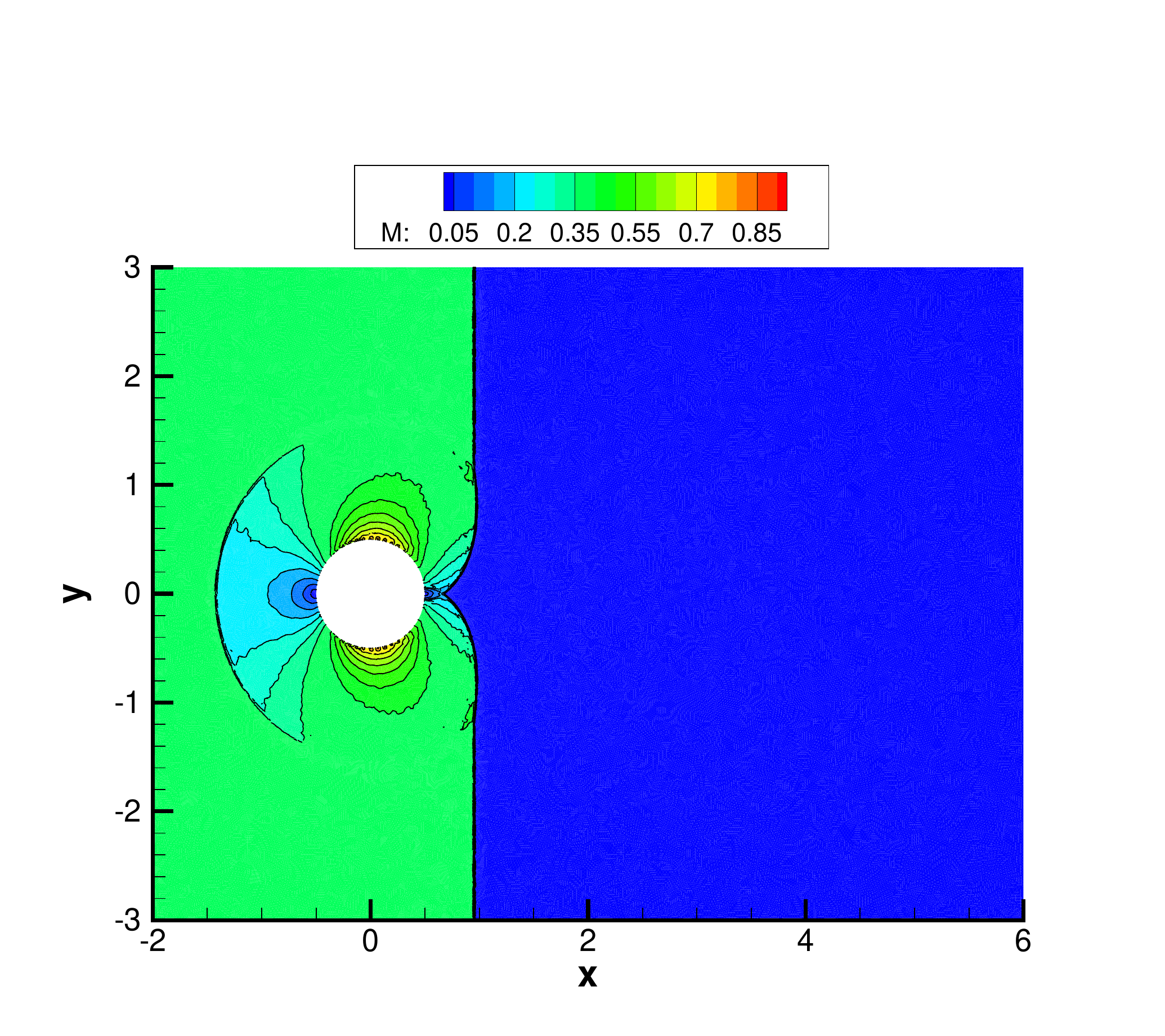}}
\subfigure[$t=1.5$]{\includegraphics[width=0.4\textwidth,trim={0cm 0cm 0cm 4.6cm},clip]{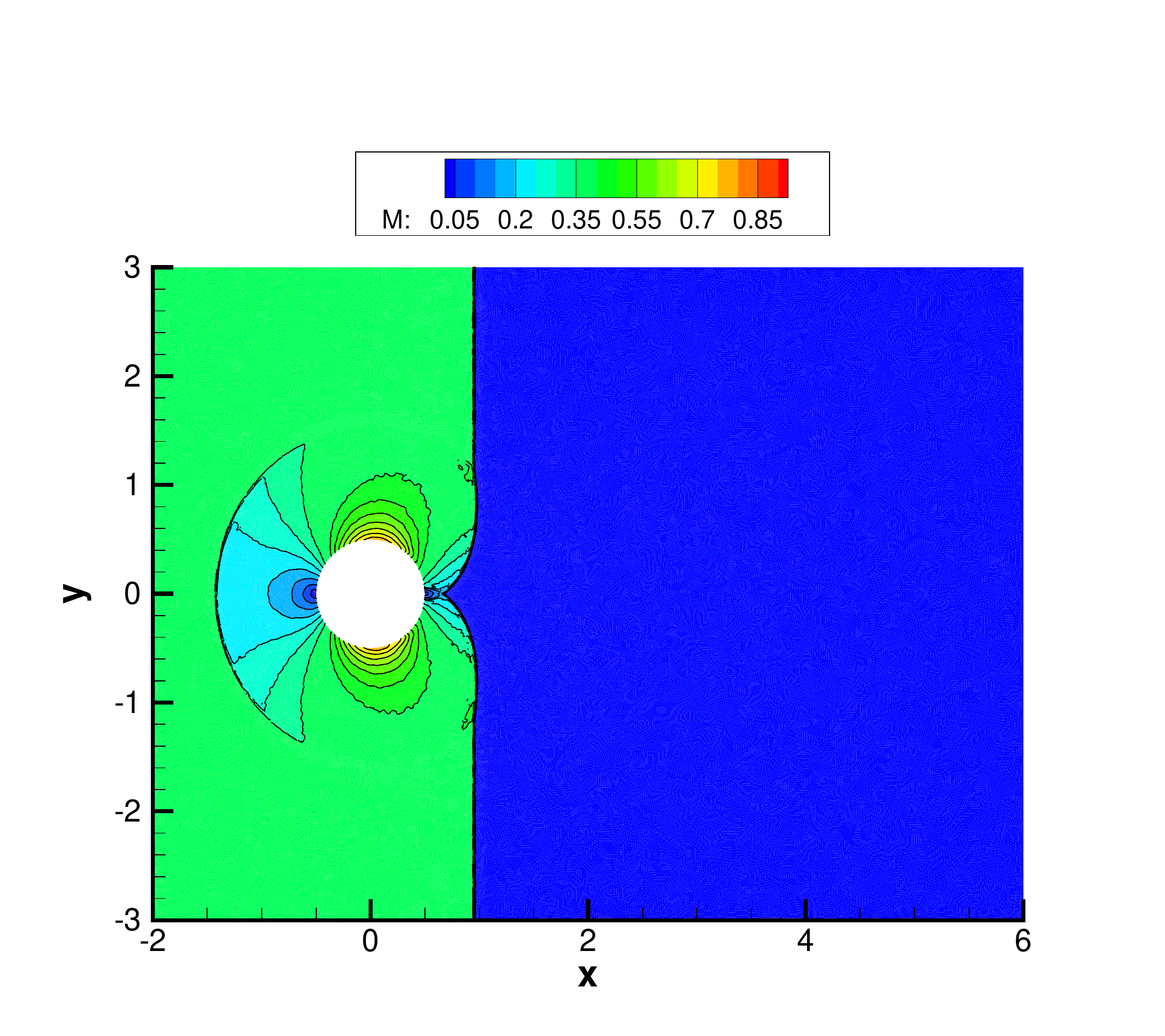}}
\qquad
\subfigure[$t=2$]  {\includegraphics[width=0.4\textwidth,trim={0cm 0cm 0cm 4.6cm},clip]{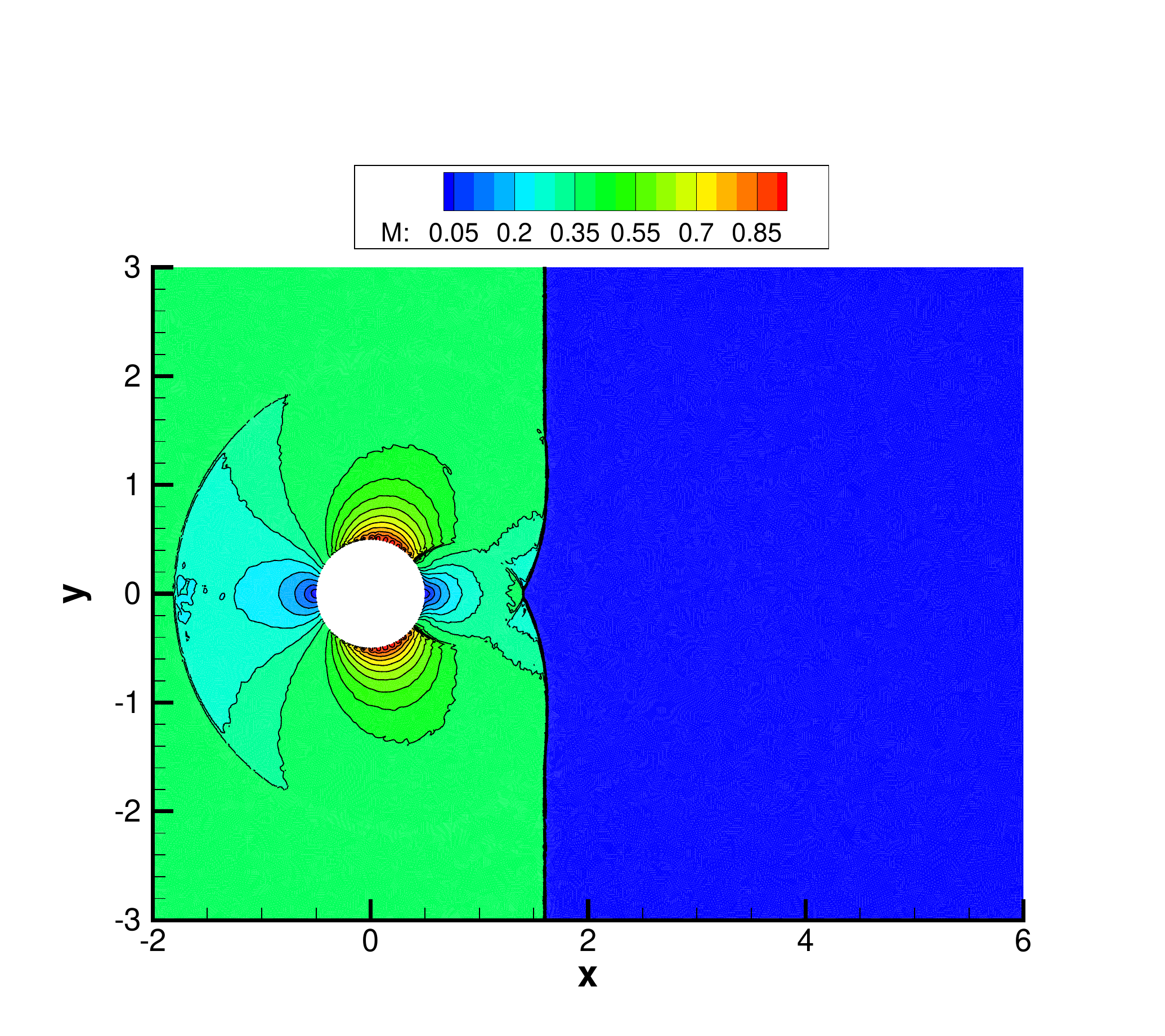}}
\subfigure[$t=2$]  {\includegraphics[width=0.4\textwidth,trim={0cm 0cm 0cm 4.6cm},clip]{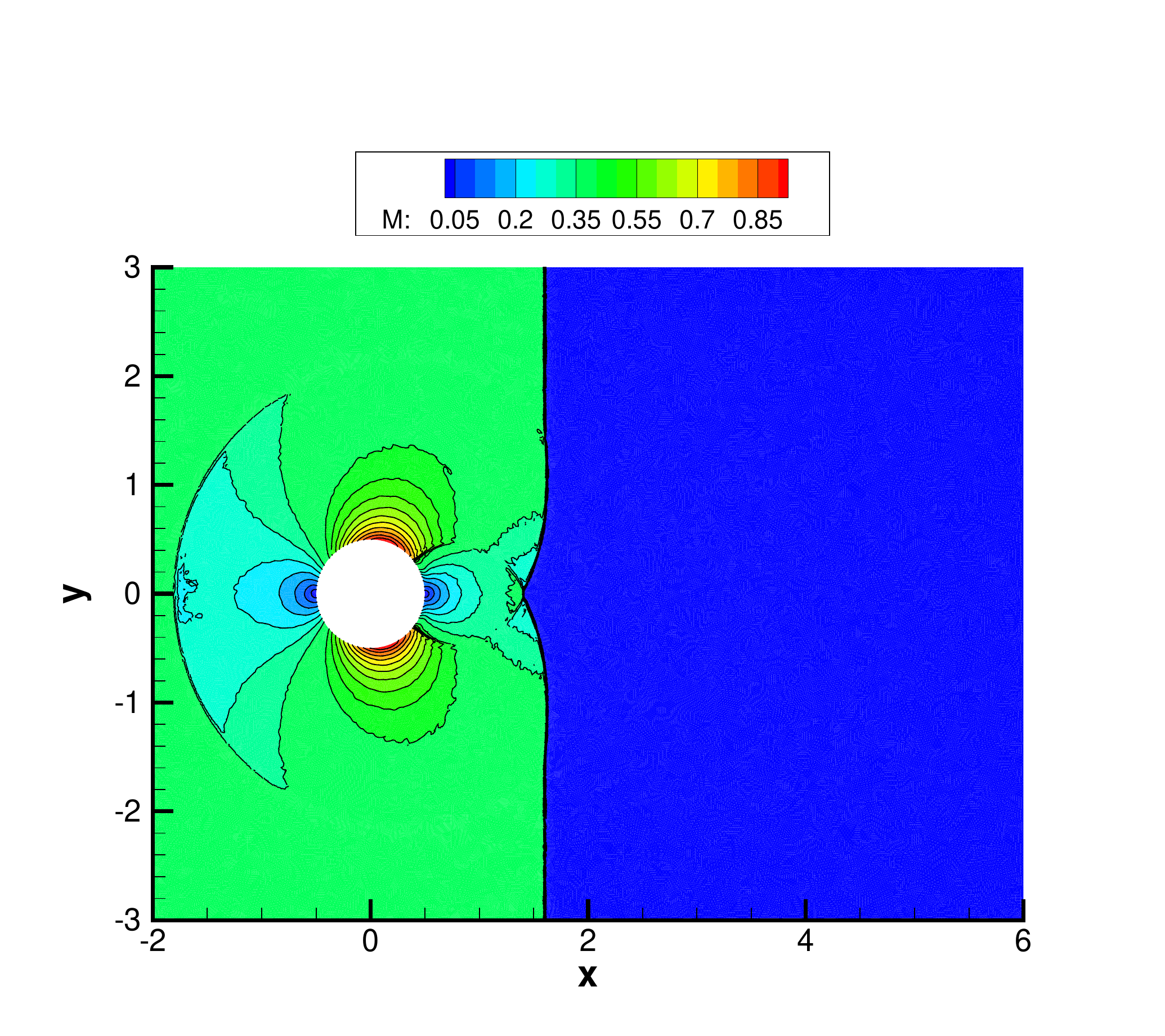}}
\caption{Shock-cylinder interaction of Section~\ref{ssec.shock-cyl2D interaction}. We show the Mach number iso-contours at different time steps \RIcolor{(simulation run with DG-$\mathbb{P}3$ on the left and DG-$\mathbb{P}3$/SBM-$\mathbb{P}3$ on the right)}.}\label{fig:ShockLimiter1}
\end{figure}

\begin{figure}
\centering
\subfigure[DG-$\mathbb{P}1$]                  {\includegraphics[width=0.4\textwidth,trim={0cm 0cm 0cm 0cm},clip]{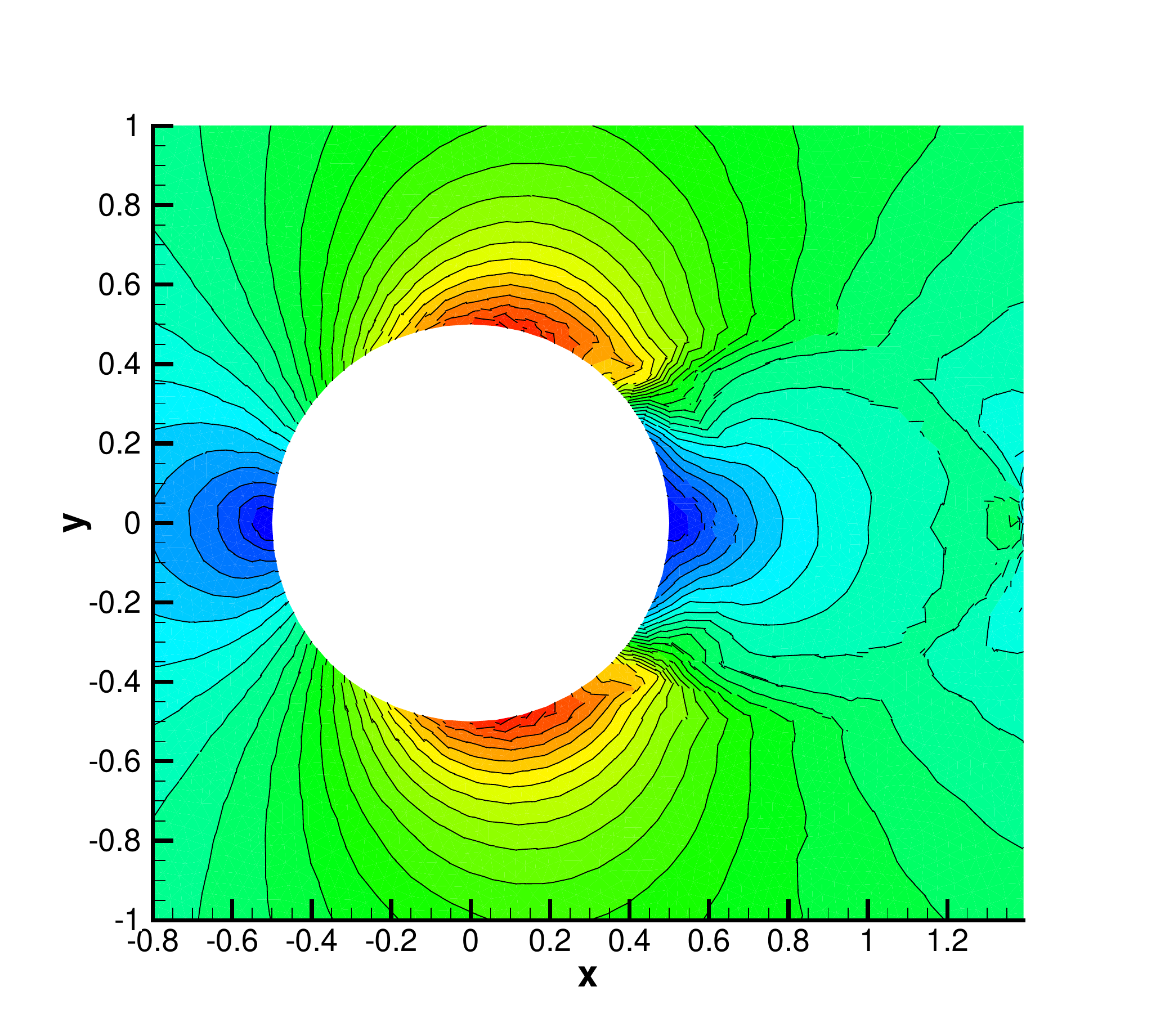}}
\subfigure[DG-$\mathbb{P}1$/SBM-$\mathbb{P}1$]{\includegraphics[width=0.4\textwidth,trim={0cm 0cm 0cm 0cm},clip]{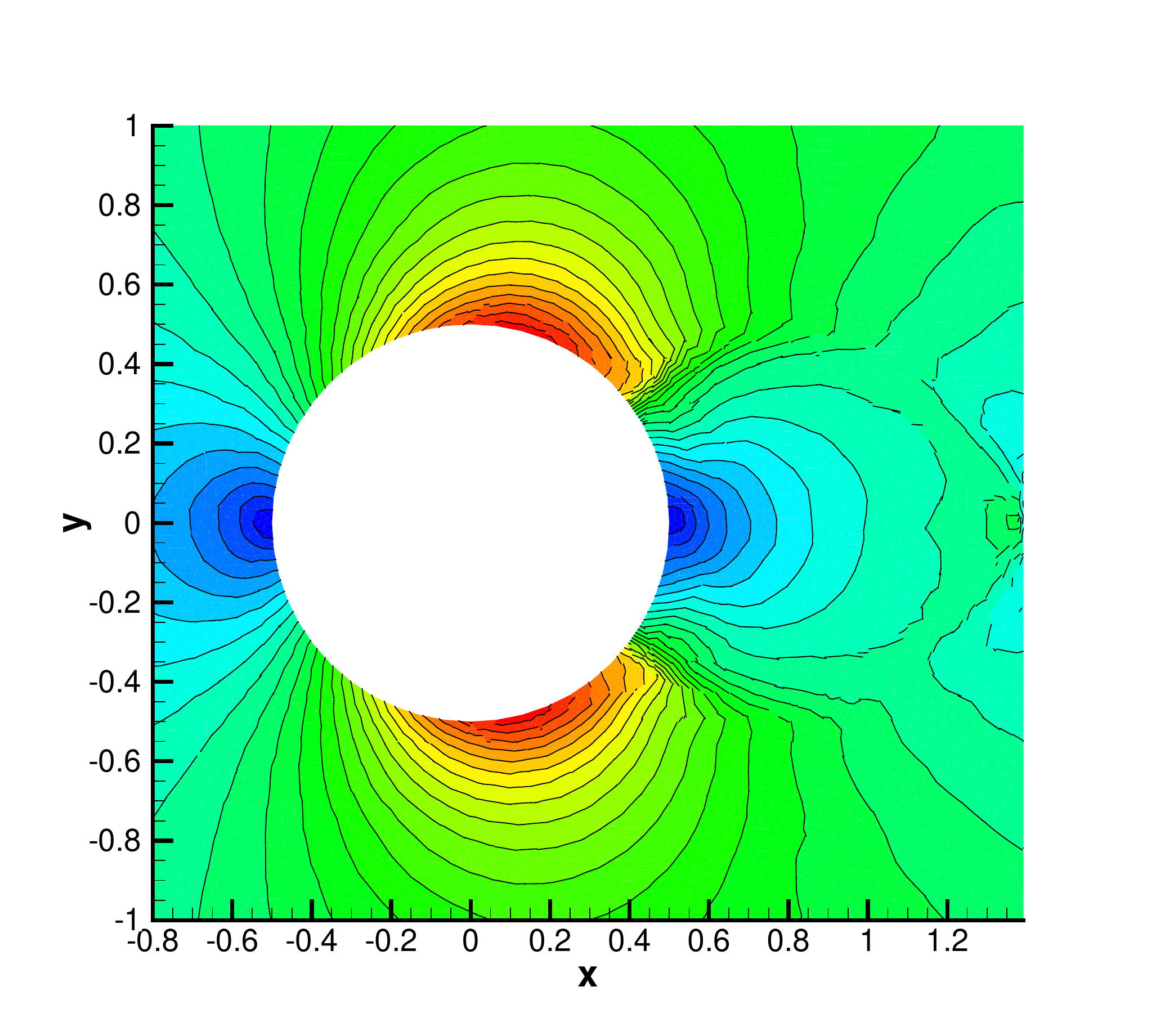}}
\qquad
\subfigure[DG-$\mathbb{P}2$]                  {\includegraphics[width=0.4\textwidth,trim={0cm 0cm 0cm 0cm},clip]{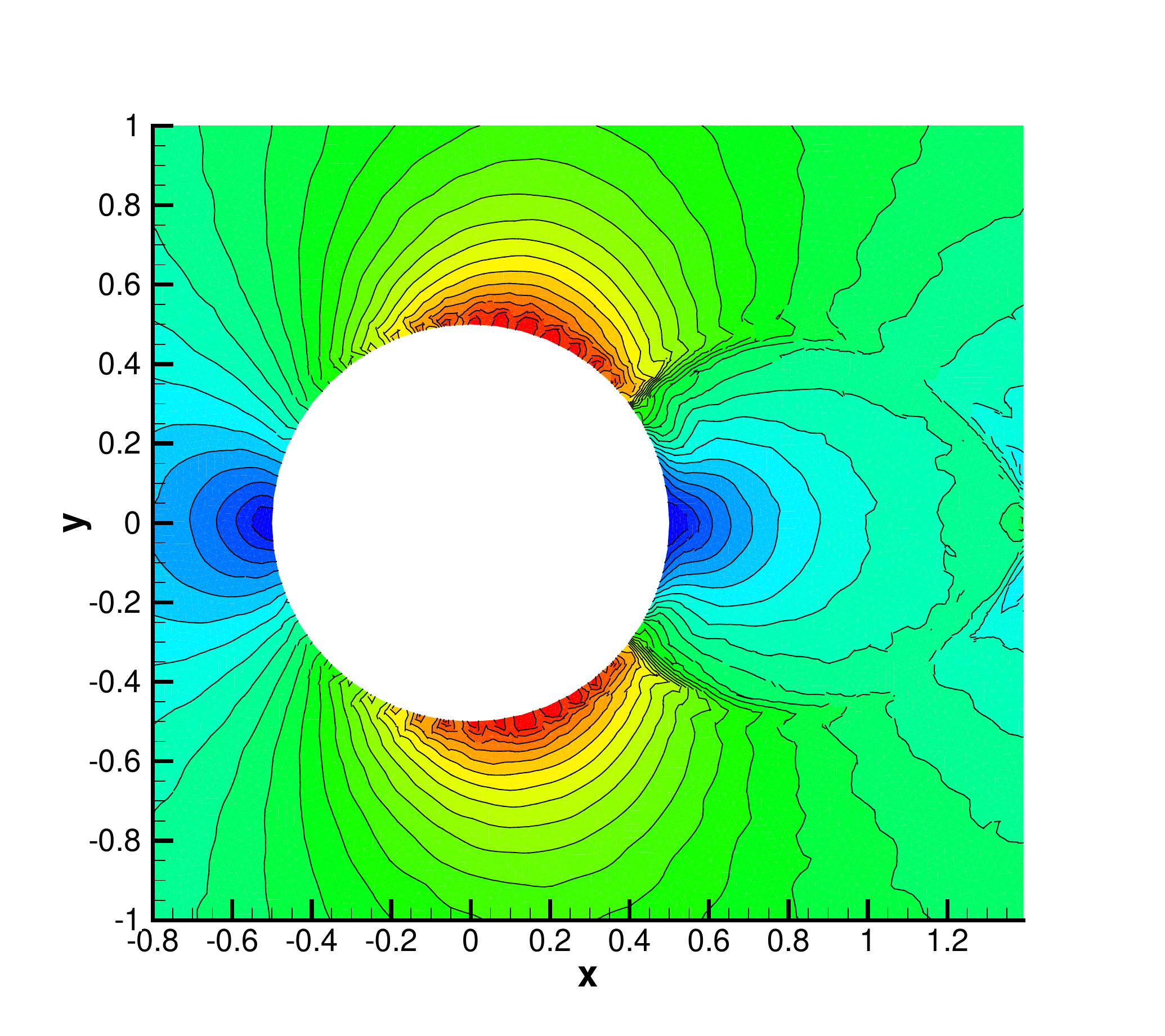}}
\subfigure[DG-$\mathbb{P}2$/SBM-$\mathbb{P}2$]{\includegraphics[width=0.4\textwidth,trim={0cm 0cm 0cm 0cm},clip]{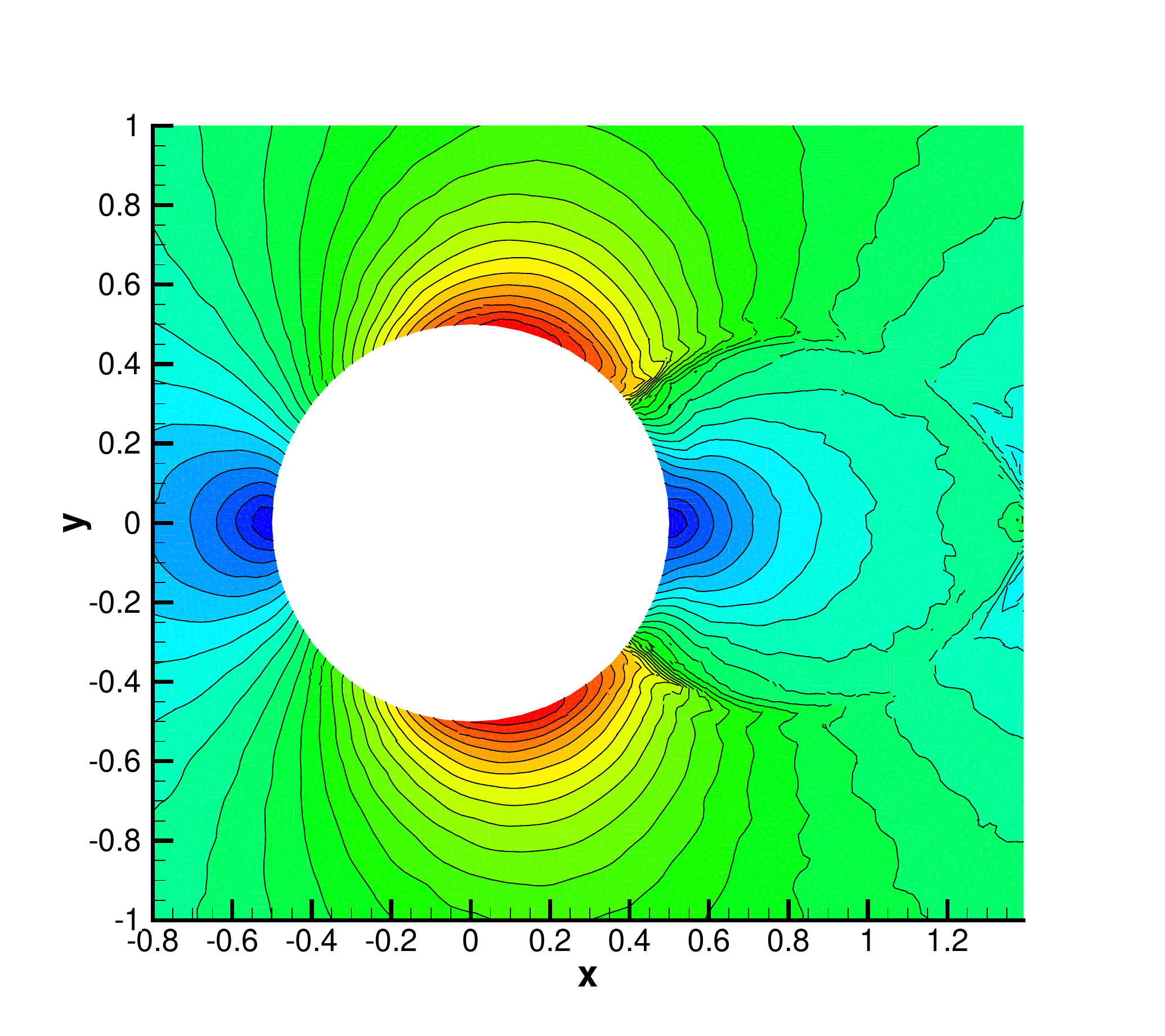}}
\qquad
\subfigure[DG-$\mathbb{P}3$]                  {\includegraphics[width=0.4\textwidth,trim={0cm 0cm 0cm 0cm},clip]{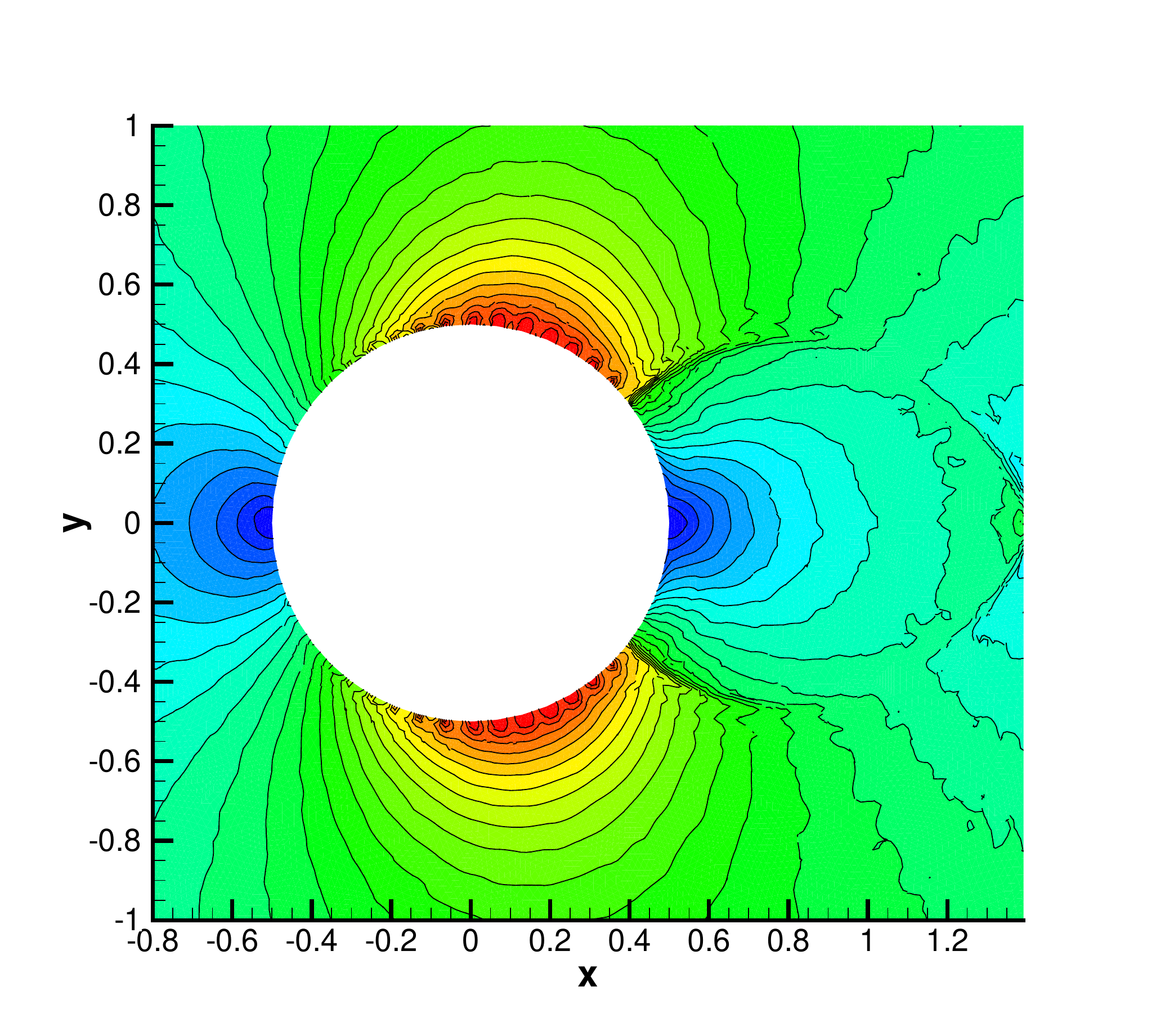}}
\subfigure[DG-$\mathbb{P}3$/SBM-$\mathbb{P}3$]{\includegraphics[width=0.4\textwidth,trim={0cm 0cm 0cm 0cm},clip]{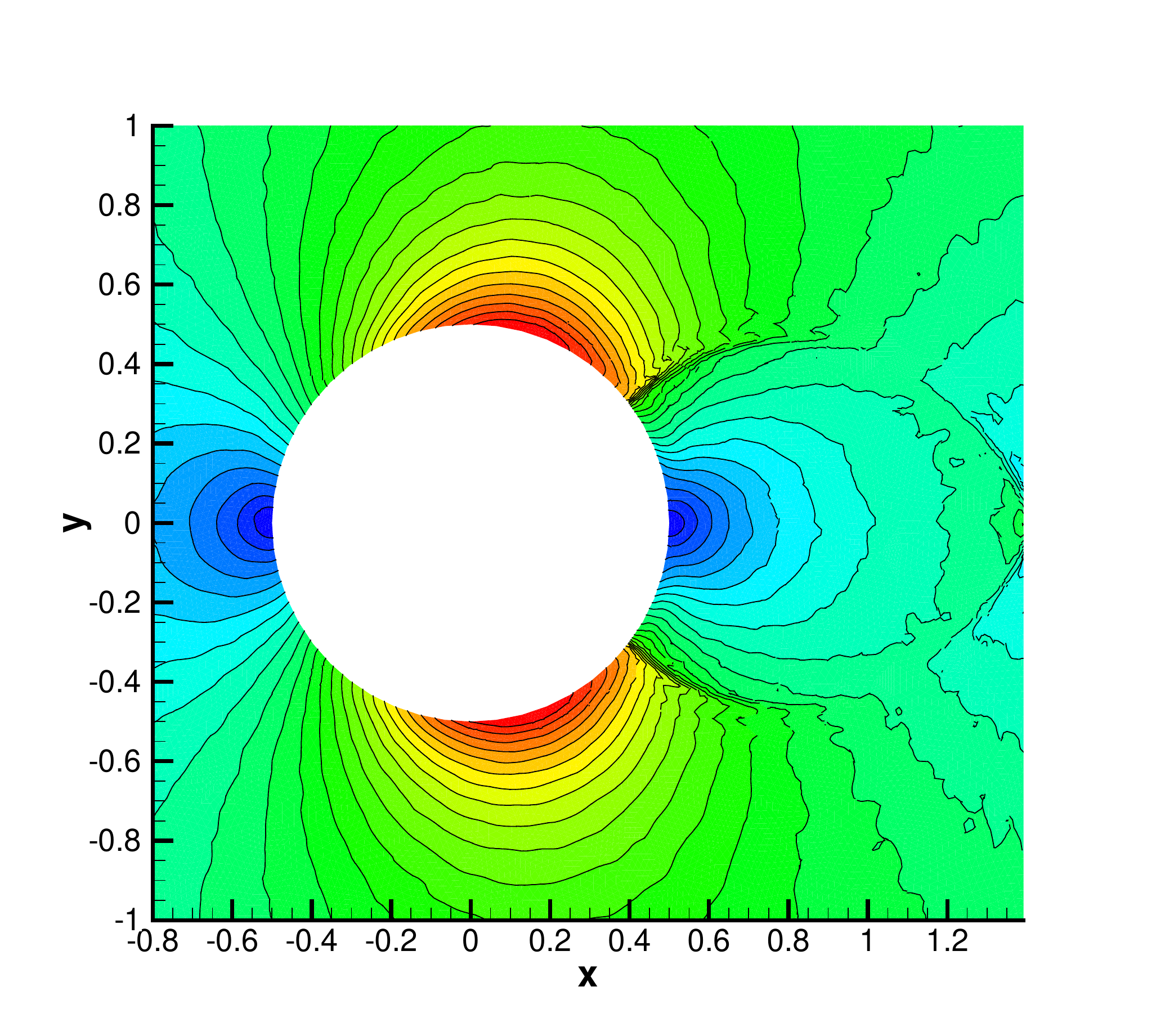}}
\caption{Shock-cylinder interaction of Section~\ref{ssec.shock-cyl2D interaction}. Mach number iso-contours at the final time $t=2$ (simulations run with several schemes).}\label{fig:ShockLimiter2}
\end{figure}

The last test case we want to address here consider the interaction of a shock wave with a two-dimensional cylinder. 
The computational domain is $[-2,6]\times[-3,3]$ discretized
with an unstructured triangulation made by 7,761 grid points and 15,198 elements. 
The employed ADER-DG scheme is supplemented with the \textit{a posteriori} sub-cell finite volume limiter.
The cylinder is centered in $(0,0)$ and has radius 0.5. 
The initial condition consists in a shock wave traveling at \textit{Mach number} $M_s = 1.3$ and is then given via the Rankine-Hugoniot conditions.
The flow upstream the shock is at rest and is characterized by density and pressure, respectively being $\rho=1.4$ and $p=1$.
The simulation has been run with polynomials $\mathbb{P}1$, $\mathbb{P}2$ and $\mathbb{P}3$ comparing the implemented wall boundary conditions, with and without the SBM flux correction.
Figure~\ref{fig:ShockLimiter0} shows the initial condition for \Rall{the Mach number distribution}. \RIcolor{The time evolution is then shown in Figure~\ref{fig:ShockLimiter1} where again we plot the Mach number distribution
computed with the DG-$\mathbb{P}3$ polynomials comparing the results obtained with and without the SBM correction. We already observe a notable improvement in the iso-contours close to the cylinder in Figure~\ref{fig:ShockLimiter1}.}
\Rall{Regarding the limiter activations, it should be noticed that no fundamental difference is observed for the classical and new wall boundary conditions.}
Finally, in Figure~\ref{fig:ShockLimiter2} we present a summary of the obtained solutions around the body with different methods.
We can observe that by increasing the order of the polynomials the results obtained with classical wall boundary conditions get worse and the SBM flux correction really introduces a notable improvement.

%--------- SECTION ---------------------------------------------------------------------------
\section{Conclusions} \label{sec.concl}
%---------------------------------------------------------------------------------------------

In this work we have provided a novel effective approach to handle boundary conditions with arbitrary high order of accuracy 
for curved domain discretized through DG schemes on simple linear meshes. 
The proposed strategy relies on the shifted boundary method that allows to overcome the second-order geometrical error due to the inconsistent treatment of curved boundaries thanks to a polynomial correction of the boundary condition and the boundary flux,
without being obliged to manage curvilinear meshes and the difficulties that go with \RIIcolor{that} (mesh generation process, isoparametric transformation and special quadrature formulas). 

\RIIcolor{In particular, as remarked in section 4.1, the original SBM approximation corresponds to 
a local change in basis functions at each quadrature point which is expensive when increasing the approximation order, and moreover unnecessary.
The approach used here is based on the use of the available  polynomial bases, whatever they are, at all quadrature points.} 
The new approach has been tested over a large set of benchmarks in 2D and 3D and with both steady and unsteady flows. 
The formal order of accuracy provided by the employed DG-$\mathbb{P}N$ schemes has been numerically retrieved in all the performed test cases,
and furthermore, the boundary corrections have also been coupled with the 
\textit{a posteriori} subcell FV limiter, thus allowing the effective simulation of shocks and discontinuities. 

Further extensions of the present work will concern first of all its application in the context of moving meshes and moving interfaces~\cite{gaburro2020high}, and then the additional development necessary for its usage in fully embedded computations~\cite{Scovazzi4,lesueur2022muct}, which do not require any conformal meshing of internal boundaries.
Finally, we also plan to use a similar approach in the context of more complex models as 
Navier-Stokes equations, for which the extension should be straightforward, 
the MHD equations~\cite{fambri2017space} where also the magnetic field should be correctly handled, 
up to the GPR unified model of continuum mechanics~\cite{dumbser2016high} for which also the distortion field should be taken into account.

\section*{Acknowledgments}
%---------------------------------------------------------------------------------------------

The authors are members of the CARDAMOM team in the Inria center at the University of Bordeaux.
E.~G. gratefully acknowledges the support received from the European Union's Horizon 2020 Research and Innovation Programme under the Marie Sk\l{}odowska-Curie Individual Fellowship \textit{SuPerMan}, grant agreement No. 101025563.  
The numerical simulations presented in this paper were carried out using the PlaFRIM experimental testbed, supported by Inria, CNRS (LABRI and IMB), Universit\'e de Bordeaux, Bordeaux INP and Conseil R\'egional d'Aquitaine~(see https://www.plafrim.fr/).

\bibliographystyle{abbrv}
\bibliography{biblio}

\end{document}